\documentclass[11pt,oneside,a4paper]{report}

\usepackage{amsmath}
\usepackage{amsthm}
\usepackage {amsfonts}
\usepackage{multirow}
\usepackage{tikz}

\newtheorem{d1}{Definition}[section]
\newtheorem{t1}{Theorem}[section]
\newtheorem{p1}{Proposition}[section]
\newtheorem{l1}{Lemma}[section]
\newtheorem*{t2}{Theorem}
\theoremstyle{definition}
\newtheorem{ex1}{Example}[section]

\newtheorem*{re1}{Remark}
\newtheorem*{re2}{Remarks}
\numberwithin{equation}{section}

\begin{document}

\author{Thanos Vasiliadis}
\title{An Introduction to Mean Field Games using probabilistic methods}
\date{Athens, 28/2/2019}
\maketitle

 \begin{abstract}
  
  This thesis is going to give a gentle introduction to Mean Field Games. It aims to produce a coherent text beginning for simple notions of deterministic control theory progressively to current Mean Field Games theory. The framework gradually extended form single agent stochastic control problems to multi agent stochastic differential mean field games. The concept of Nash Equilibrium is introduced to define a solution of the mean field game. To achieve considerable simplifications the number of agents goes to infinity and formulate this problem on the basis of McKean-Vlasov theory for interacting particle systems. Furthermore, the problem at infinity is being solved by a variation of the Stochastic Maximum Principle and Forward Backward Stochastic Differential Equations. To elaborate more the Aiyagari macroeconomic model in continuous time is presented using MFGs techniques.

 \end{abstract}

 \chapter*{Preface}
 
 This text started as my master thesis for completion of the M.Sc. program "Mathematical Modeling in New Technologies and Financial Engineering" offered by the National Technical University of Athens (NTUA) but soon exceeded its purpose and transformed into a text that I hope will provide the foundation for my future research in the Mean Field Games topic. 
 
 I aimed through the development of this text, to understand the topic in a sufficient depth that would allow we to review the most important literature, unify it and explain in a way that it can be helpful to anyone interested in studying the topic with minimum mathematical prerequisites. 
 
 Mean Field Games themselves are a heavy topic to discuss in any form. The requirements I identified as I engaged in this topic were Optimal Control Theory, Stochastic calculus, Stochastic Control theory and game theory. Even though one can deal with Mean Field Games without any prior knowledge of Game theory, it enhances a lot the intuition behind the models.
 
 In the same spirit I would like to advise the non-expert reader or generally anyone who wants to develop a feeling for MFGs and lack the background to start from the Appendixes and the move forward to the main text. They  have been designed to be able to stand on their own and provide a quick introduction and  review of the related subject without too much technical details. Nevertheless, I also have references to the appropriate sections of the appendix throughout the text. 
 
 Last but not least, I would like to thank my advisor professor Vassilis Papanicolaou for his support and encouragement throughout the course of this project. He has been an invaluable mentor and teacher. Our discussions has provided me with motivation and insight for various subjects broader than mathematics alone. 
 
 $\vspace{1mm} $
 
 \author{Thanos Vasiliadis}

 \tableofcontents
  \markright{}
\addcontentsline{toc}{chapter}{Preface}

\chapter{Introduction}

Mean field games (MFGs for short) are a relatively new branch of mathematics and more specifically they lie at the intersection of game theory with stochastic analysis and control theory. Since their first appearance in the seminal work of Lions and Lasry (2006) and independently by  Huang, Malhame and Caines(2006) two approaches have been proposed to study them, the coupled Hamilton-Jacobi-Bellman with Focker-Plank which comes from dynamic programming in control theory and PDEs and the Forward-Backward Stochastic Differential Equations (FBSDEs) of McKean-Vlasov type which comes from stochastic analysis. We will explain both of them in the appropriate section, however we will rely mostly on the second one for our analysis. Our program to deal with mean filed games is as follows: 
\begin{itemize}
 \item
 In the first chapter we will introduce MFGs intuitively, and the basic notions from game theory to improve the readability of the text for the non-expert reader. Moreover we will discuss briefly about large games to justify our probabilistic approach following in chapter 2 and 3
 \item
 In the second chapter we will develop our formal definitions about MFGs, we will explain the transition to the limit with motivation from statistical physics and introduce McKean-Vlasov stochastic differential equations to describe the dynamics of our system.
 \item
 In the third chapter we will be involved in solving MFGs and identifying equilibrium points. We will use a version of the Stochastic Maximum Principle along with Schauder's fixed point theorem to prove existence of a MFG equilibrium.
 \item
 In the fourth chapter in order to see the theory in action we will present and solve the Aiyagari model, a toy macroeconomic model. 
\end{itemize}

Let's start by addressing intuitively the question ``What is actually a mean field game?''

\section{What is a mean field game and why are they interesting?}

As the name proposes it is a strategic game dynamic and symmetric between a large number of agents  in which the interactions between agents are negligible but each agent's actions affect the mean of the population. In other words each agent acts according to his minimization or maximization problem taking into account other agents' decisions and because their population is large we can assume the number of agents goes to infinity and a representative agent exists (precise definitions for everything will be given later).

In traditional game theory we usually study a game with 2 players and using induction we extend to several, but with games in continuous time with continuous states (differential games or stochastic differential games) this strategy cannot be used because of the complexity that the dynamic interactions generate. On the other hand with MFGs we can handle large number of players through the mean representative agent and at the same time describe complex state dynamics.

MFG are becoming an increasingly popular research area because they can model a large variety of phenomena from large systems of particles in physics, to fish schooling in biology, but we will restrict ourselves here to economics and financial markets. We will devote the last part of this introductory section to MFG and economics to give further motivation.

We will now see one of the most common examples that accompanied mean field games since its early development \cite{GLL2010}.

\subsection{When does the meeting start?}

Suppose that we have $N$ university professors participating in a meeting, which is scheduled to begin at $\mathbf{t_0} $ (called $\mathit{initial \hspace{1mm} time}$ ).  All of them start from different locations to attend but are symmetric in a sense that they share the same characteristics (for example they have to cover the same distance to the venue or they are moving with the same speed, they need the same relative time etc). But because some of them are notorious for being late the organising committee decided to actually start the meeting only when the 75$\%$ of them gather to the venue. Each one given his or her preferences have a $\mathit{target} \hspace{1mm} time$ of arrival, $\mathbf{t_i}$ but due to non-anticipated events (weather conditions, traffic etc) they arrive at $\mathbf{X_i}$ ($\mathit{actual \hspace{1mm} time}$ of arrival). Each $X_i$ is the sum of professor's desired arrival time ($t_i$) which is completely under his control and random noise.
$$X_i=t_i+\sigma_i\epsilon_i \text{ for $i=1,2,...N$} $$

\begin{itemize}
 \item $(\epsilon_i)_{1\leq i \leq N}$ is an iid  sequence with $N(0,1)$
 \item $(\sigma_i)_{1\leq i \leq N}$ is also an iid  sequence with common distribution $\nu$
 \item $(\epsilon_i)_{1\leq i \leq N}$ is assumed to be independent of $(\sigma_i)_{1\leq i \leq N}$
\end{itemize}

So the actual time the meeting starts, $\mathbf{T}$ is a function of the empirical distribution $\bar{\mu}_X^N$ of the arrival times $X=(X_1,...,X^N)$.

The expected overall cost of professor i is defined as:
\begin{equation}
 J_i(t_1,...,t_N)=\mathbb{E}[\underbrace{A(X_i-t_0)^+}_{\text{reputation cost}}+\underbrace{B(X_i-T)^+}_\text{inconvenience cost}+\underbrace{C(T-X_i)^+}_\text{cost of early arrival} ]
\end{equation}

where $A,B,C>0$ constants

\begin{re1}

The fact that the choice of the start time T is a function of the empirical distribution $T=\tau(\bar{\mu}_X^N)$ is the source of the interactions between the agents who need to take into account the decisions of the other agents in order to make their own decision on how to minimize their cost.

\end{re1}

\section{Introduction to game theory}

\begin{d1}{Game}

 A (strategic) game is a model  of interacting agents (players), who take decisions.
\end{d1}

We can separate games in to four main categories:

\begin{itemize}
 \item \textit{Information}\\
  Regarding the structure of the available information:
 \begin{enumerate}
  \item Games of perfect or complete information, and 
  \item Games of partial or incomplete information usually called Bayesian games. 
 \end{enumerate}

 \item \textit{Time}\\
 Regarding time:
 \begin{enumerate}
  \item Static or one shot games in which the agents take only one decision regardless of the time horizon. 
  \item Dynamic games in which the agents take multiple decisions at discrete times. These games can be specified even more in:
  \begin{enumerate}
   \item Discrete time, or repeated games where the time is discrete and the dynamic game consists of multiple one shot games which are repeated in different time instances.
   \item Continuous time or differential games where the time is continuous and the agents take actions in a continuous manner i.e. use continuous functions to represent their decisions
  \end{enumerate}

 \end{enumerate}

\end{itemize}

\begin{d1}{Some terminology}

In order to define a game we need the following
 \begin{enumerate}
  \item
  $P$: the set of players (agents)\\
  $\#(P)=N$ the number of players
  \item
  $A^i$:  the set of actions for player i\\
  $A=A^1\times... \times A^N$\\
  $a=(a_1,...a_N)\in A$ is an action profile where $a_i$ is the action the individual players take and $a_{-i}$ the action profile including every player's action except $i$'s  $a_{-i}=(a_1,...a_{i-1},a_{i+1},...a_N)$ 
  \item
  \begin{enumerate}
   \item  $C$: the set of players' characteristics
   \item$\prec_i$ a preference relationship that partially orders $A^i$ and defines utility functions
  \end{enumerate}
  \item
  $\mathcal{U}$: the set of payoff functions $u_i:A_i\to \mathbb{R}$
  \item
  $\mathcal{M}$: the set of players' strategies
  $$\mathcal{M}:=\{f:[0,T]\to A| f \text{ arbitrary function} \} $$
  \item
  $\mathcal{P}(A)$ The set of probability measures on A
 \end{enumerate}

\end{d1}

$\vspace{1mm} $

\begin{re1}

 For this section we will assume that all players desire higher payoffs we will not go into details about utility functions or preferences or rational behavior of players since these concepts are broader than the scope of this text. We assume that $u_i:A_i \to \mathbb{R} $ is well defined and fulfills common assumptions which we will not mention.  We refer to the original work of Von Neumann and Morgensten "Game theory and Economic Behavior " and to almost any textbook in game theory for more information.
\end{re1}

$\vspace{1mm} $

Game theory is mostly concerned with the incentives of the agents. The main question is the existence of a strategic situation from which no one has incentive to deviate, the so called Nash equilibrium.

$\vspace{2mm} $

\begin{d1}{Nash Equilibrium}

An action profile $a^*\in A$ is called a Nash equilibrium if and only if for every player i
$$u_i(a^*)\geq u_i(a_i,a^*_{-i}) \hspace{1mm} \forall a_i\in A^i  $$
where $u_i(\cdot)$ is the payoff function of player i. 

\end{d1}

$\vspace{1mm} $

When we are solving a game we suppose that every player is acting according to his/her best interests, trying to respond optimally to other players actions. This is the concept of the best response function.

$\vspace{2mm} $

\begin{d1}{Best Response Function (BRF)}

The function $B_i:A \to A^i$ is called the best response of player $i$ to the actions of the other players denoted by $-i$.

$$B_i(a_i)=\{a_i\in A^i:u_i(a_i,a_{-i})\leq u_i(a^{'}_i,a_{-i}) \hspace{2mm} \forall a^{'}_i\in A^i \} $$

\end{d1}

$\vspace{2mm} $

Equivalently we can define Nash equilibrium in terms of Best Response Functions

\begin{d1}{Equivalent definition}

An action profile $a^*\in A$ is called a Nash equilibrium if and only if it is a fixed point of the best response function B, $B=B_1\times ... \times B_N$

\end{d1}

$\vspace{2mm} $

We will now give simple examples of games to elaborate more on the definitions and theory.

$\vspace{2mm} $

\begin{ex1}{\textbf{Prisoner's  dilemma}}

Suppose a robbery is committed and the police arrests the two suspects. Policemen decide to question them independently to increase their chances to unfold the truth. Each one can accuse the other but has also right to remain silent. If both of them accuse each other policemen will be sure that they are guilty and send them to prison for 5 years. In case both of them remain silent because of lack of details they will be sentenced only for one year and if one accuses and the other not then the one who accused the other will be free to leave and become a witness so the other will be sentenced to 10 years in jail. Policemen inform the suspects of the four possibilities, and they have to announce their decision simultaneously.

$\vspace{1mm} $

\textbf{Game formulation}\\
Players: 2\\
$A=\{\text{accuse, not accuse} \}$\\
The preferences relationship is defined as follows: The most preferable situation freedom is labeled 3, the next preferable situation 1 year in prison is labeled 2, the next 5 years is labeled 1 and the least 10 years in prison is labeled 0. This way we can use the usual order of natural numbers for the outcomes. Thus an action which yields a higher payoff is more preferable. We wrap everything in the next table of payoffs

\begin{table}[!h]
\caption{Prisoner's dilemma payoffs}
\centering
\begin{tabular}{r l|c|c|}
\multicolumn{1}{r}{}&\multicolumn{1}{r}{} & \multicolumn{2}{c}{suspect 2} \\
\cline{3-4}
{}& {} & accuse & not accuse \\
\cline{2-4}
\multicolumn{1}{ r|  }{\multirow{2}{*}{suspect 1} } & accuse & 1,1 & 3,0\\
\cline{2-4}
\multicolumn{1}{ r|  }{} & not accuse & 0,3 & 2,2\\
\cline{2-4}
\end{tabular}

\end{table}

\textbf{Solution}\\
We have the following strategic situations:
\begin{itemize}
 \item Assume player 2 plays \textbf{accuse}, then player 1 plays \textbf{accuse} since it has the greatest payoff.
 \item Assume player 2 plays \textbf{not accuse}, then player 1 players \textbf{accuse} since it has the greatest payoff.
\end{itemize}
The BRF for player 1 is accuse (whatever player 2 plays)
$$B_1=\{\textbf{accuse}\} $$
similarly
$$B_2=\{\textbf{accuse}\} $$

we conclude that the Nash equilibrium of the game is (\textbf{accuse, accuse})

\end{ex1}

$\vspace{2mm} $

\begin{ex1}{\textbf{Matching Pennies}}

Two players choose to show each other simultaneously the face of a coin if they choose the same player 2 pays player 1 1$\$$ if they choose different player 1 pays player 2 1$\$$.

$\vspace{1mm} $

$\mathbf{Game \hspace{1mm} formulation}$\\
Players: 2\\
$A=\{\text{Head, Tails} \}$\\
Payoffs:

\begin{table}[h]
\caption{Matching Pennies Payoffs}
\centering
\begin{tabular}{r l|c|c|}
\multicolumn{1}{r}{}&\multicolumn{1}{r}{} & \multicolumn{2}{c}{player 2} \\
\cline{3-4}
{}& {} & Head & Tails \\
\cline{2-4}
\multicolumn{1}{ r|  }{\multirow{2}{*}{player1} } & Head & 1,-1 & -1,1\\
\cline{2-4}
\multicolumn{1}{ r|  }{} & Tails & -1,1 & 1,-1\\
\cline{2-4}
\end{tabular}
\end{table}

\textbf{Solution}\\
This game has no Nash equilibrium. \\
Suppose $(Head,Head)$ is a Nash equilibrium, than player 2 will be in better position if he/she change his/her decision to Tails. So equilibrium moves to $(Head, Tails)$ but then again player 1 will be in better position if he/she change his/her decision to Tails and so on. There is no stable outcome, each player has incentive to deviate form any situation and so no Nash equilibrium exists.

\end{ex1}

$\vspace{2mm} $

\begin{ex1}{\textbf{Cournot Duopoly}}

Suppose we have a market with two firms producing the same product. Both of them face a common demand curve $Q=a-P$ where $a,b>0$ with $Q$ the total product ($Q=q_1+q_2$) and $P$ the price. Each firm has a linear cost function $C(q_i)=cq_i$ and try independently to maximize their profits $\Pi_i(q_i)=q_iP-cq_i$, for  $i=1,2 $

$\vspace{1mm} $

$\mathbf{Game \hspace{1mm} formulation}$\\
Players: 2\\
$A=\{q_1, q_2 \}$ (here the actions are continuous variables)\\
Payoffs: $\Pi_1(q_1), \Pi_2(q_2)$

\textbf{Solution}\\
For each firm the profit function can be expressed as:
$$\Pi_i(q_i;q_{-i})=q_i(a-(q_1+q_2))-cq_i \text{ for $i=1,2$}$$
using the inverse demand curve.\\
The first order condition for profit maximization of firm 1 yield:
$$\frac{d\Pi_1}{dq_1}(q_1;q_2)=a-2q_1-q_2-c=0$$

\begin{equation}
q_1=B_1(q_2)=\frac{a-c-q_2}{2}
\end{equation}

This is BRF of firm 1. \\
Similarly the BRF of firm 2 is
\begin{equation}
q_2=B_2(q_1)=\frac{a-c-q_1}{2}
\end{equation}

So for the Nash equilibrium we are looking for an intersection point in the system (1.2.1)(1.2.2) which yields

$$
\left\{
\begin{split}
 &q_1^*=\frac{a-c}{3}\\
 &q_2^*=\frac{a-c}{3}
\end{split}
\right.
$$

$(q_1^*,q_2^*) $ is the unique Nash equilibrium, as the following figure shows.

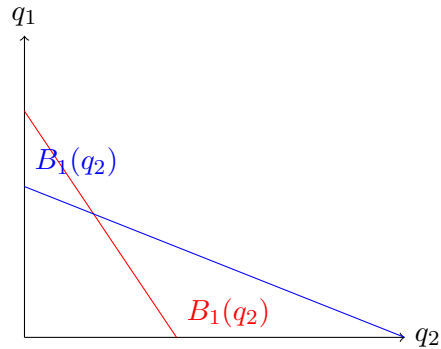
\begin{figure}[!h]
\centering
 \begin{tikzpicture}
    \draw[->] (0,0) -- (5,0) node[right] {$q_2$};
    \draw[->] (0,0) -- (0,4) node[above] {$q_1$};
    \draw[color=red][domain=0:2]     plot (\x,{1.5*(2-\x)})             node[above right]{$B_1(q_2) $};
    \draw[color=blue][domain=0:5] (0,2)  node[above right] {$B_1(q_2)  $}    plot (\x,{0.4*(5-\x)})  ;
  \end{tikzpicture}
 \caption{BRF functions}
\end{figure}

\end{ex1}

\begin{d1}{Symmetric game}

 A game is called symmetric if
 \begin{enumerate}
  \item Each player has the same action set
  $$A^1=...=A^N $$
  \item And his/her preferences can be represented by utility functions $u_i,u_j$ such that $u_i(a_1,a_2;a^{-i-j})=u_j(a_2,a_1;a^{-i-j}) \hspace{1mm} \forall (a_1, a_2)\in A $
 \end{enumerate}

\end{d1}

All of the previous examples, including "When does the meeting star?" are static, symmetric and as we are going to discuss in section 2, symmetry is one of the core characteristics of mean filed games.

\section{Nash equilibrium in Mixed Strategies}

As we saw in the example (1.2) a game does not always have a Nash equilibrium. But what would happen if we allow the players to randomize their behavior? Let look again at "Matching Pennies" while we allow players to choose their actions based on probability.

$\vspace{2mm}$

\begin{ex1}{\textbf{Matching Pennies with randomized behavior}}

Assume \textbf{player 2} chooses \textbf{Head} with probability $\mathbf{q}$ (and Tails with $1-q$) then \textbf{player 1} chooses \textbf{Head} with probability $\mathbf{p}$ (and Tails with $1-p$) and keeping in mind table 2 each outcome ((H,H),(H,T),(T,H),(T,T)) has probability $qp$. Now let look at the following situations:
  \begin{itemize}
   \item If player 1 chooses \textbf{Head} with probability 1 his expected payoff would be:
   $$q1+(1-q)(-1)=2q-1 $$
   \item If player 1 chooses \textbf{Tails} with probability 1 his expected payoff would be:
   $$q(-1)+(1-q)(1)=1-2q $$
  \end{itemize}
So if $q<\frac{1}{2}$ then he/she is in better position playing Tails and vice versa for $q>\frac{1}{2}$. For $q=\frac{1}{2}$ then $p=\frac{1}{2}$(each strategy gives the same expected payoff).
The best response of player 1 is:
$$B_1(q)=\begin{cases}
          \{0\} &\text{ if $q<\frac{1}{2}$}\\
          p=\frac{1}{2} &\text{ if $q=\frac{1}{2}$}\\
          \{1 \} &\text{ if $q>\frac{1}{2}$}\\
         \end{cases}
 $$

And similarly we can construct the BRF for player 2. Combining them and noticing the fixed point we conclude that the unique Nash equilibrium is when each one is randomizing with $p=q=\frac{1}{2}$.

This equilibrium has a special name called Nash equilibrium in mixed strategies as the following definitions indicates.

\end{ex1}

\begin{d1}{Mixed strategy}

 A mixed strategy for a player in a strategic game is a probability distribution, $\mu \in \mathcal{P}(A^i) $ for his/her actions given the actions of the other players.

\end{d1}

In mixed strategies each player randomize his/her actions according to the distribution $\mu_i$ which is a probability measure defined on $A_i$ the set of player's actions. While, $\mu$ is the product measure defined on $A$ the Cartesian product of $A_i$s.

\begin{d1}{Nash equilibrium in mixed strategies}

 A mixed strategy profile $\mu^* \in \mathcal{P}(A) $ is called a Nash equilibrium in mixed strategies if and only if for every player $i$
 $$ u_i(\mu^*)\geq u_i(\mu) \hspace{2mm} \forall \mu \in \mathcal{P} (A)   $$
 where $u_i$ the utility that player $i$ gets from probability distribution $\mu$
\end{d1}

It is rather obvious that a pure-strategies Nash equilibrium is an equilibrium in degenerate mixed strategies.

\begin{re1}
 Again we assume $u_i:\mathcal{P}(A^i) \to \mathbb{R}$ is well defined satisfying certain assumptions.
\end{re1}

Again as we saw in the previous example we can define Best Response Functions in terms of mixed strategies in the same way as we did with pure strategies. And of course we have the equivalent definition of Nash equilibrium in mixed strategies as a fixed point of the BRFs.

\begin{t1}{Nash}

 Every strategic game with a finite action set, has a Nash equilibrium in mixed strategies
\end{t1}

\section{Games with a large number of players}

Nash's theorem is quite general and was proved by himself in a very elegant way but in order to stress the importance of finiteness in the theorem we will use an example were we violate this finiteness and highlight the need for measure-theoretic tools to analyse games with a large number of players. 

We consider a game where the number of players is infinite and set-up a rule to introduce this infiniteness in the strategy profiles. This counterexample is by Peleg (1969) \cite{Peleg1969}

\begin{ex1}{A game with infinite number of players}

This game has no physical interpretation, we suppose the Positive natural numbers is the set of players, each player is given to options to player 0 or 1 and the payoff he gets depends upon the choice he made and the sum of the actions of the other players as follows.

\textbf{Game formulation}\\
$P\equiv \mathbb{Z} $ the set of players\\
$A_i=\{\ 0,1 \}$ individual's action set \\
$u_i(a_i,\displaystyle\sum_{i\in P} a_i)=\begin{cases}
 1  & \text{if $a_i=1$ and $\displaystyle\sum_{i\in P} a_i<\infty$}\\
 0  & \text{if $a_i=0$} \\
 -1  & \text{if $a_i=1$ and $\displaystyle\sum_{i\in P} a_i=\infty$}\\
\end{cases}
$ payoff function.

The game does not have a Nash equilibrium in pure strategies nor in mixed. For \textbf{pure strategies}, suppose that that \textbf{1} is a Nash equilibrium, i.e. all players choose to play 1 then $\displaystyle\sum_{i\in P} a_i=\infty$ and $u_i=-1 \forall i \in P$. But then everyone would be in better position if they choose \textbf{0}, to get a payoff of 0 and  $\displaystyle\sum_{i\in P} a_i<\infty$ but then again everyone would have incentive to deviate and play \textbf{1} to get a payoff of 1 but then $\displaystyle\sum_{i\in P} a_i=\infty$ e.t.c. 

For \textbf{mixed strategies} suppose that $\mu\in \mathcal{M}(A)$ is an equilibrium probability distribution and $(p, 1-p)$ the corresponding mixed strategy equilibrium profile for player $i$, then his expected payoff under $\mu_{-i}^*$ over $A_{-i}$ would be:
$$u_(\mu_i,\mu_{-i}^*)= \mathbb{E}_{\mu_{-i}^*} \big[\mu(\{0\})u_i(0,\displaystyle\sum_{k\neq i} a_k)+\mu(\{1\})u_i(1,\displaystyle\sum_{k\neq i} a_k)\big]$$
$$=p\mathbb{E}_{\mu_{-i}^*}\big[u_i(1,\displaystyle\sum_{k\neq i} a_k)\big]=-p $$

Then he would would get a payoff  $-p$ depending on the convergence of the sum $\displaystyle\sum_{k\neq i} a_k$ which is given by $\mu(\underset{i\to \infty}{liminf}(\{a_i=0\}))=0 $. This way  player $i$ would gain if he play a pure strategy (0) and hence $\mu$ could not be a mixed strategy and we come back to the pure strategies case. 

\end{ex1}

As seen in the previous example solving games with a large number of players can be difficult (remember also eq. (1.1.1) from "When does the meeting start?") if not impossible and the strategy behind MFGs, to overcome this difficulty, is to search for simplifications in the limit $N\to\infty$ of large games. Of course simplifications are not for free we need to rely on assumptions to achieve them, but we will come to them in the second part where we will give our formal definitions. 

For now we will describe interactions based on empirical distributions of individual responses intuitively using the example "When does the meeting start?" 

\subsection{Revisiting "When does the meeting start?"}

We restate the key elements of the model in a compact way. 

$\vspace{1 mm}$

\textbf{Game formulation}
\begin{itemize}
 \item $P=\{1,2,...N \}$
 \item $A_i=[0,E]$ the players can choose any positive time $t_i$, with $0$ representing maybe the start of the day and $E$ the end of the event but that is not important for our analysis\\
 $A=A_1\times... \times A_N=[0,E]^N$

 \item 
$J_i(t_i,\tau(\bar{\mu}_{X_{-i}}^{N-1}))$ the payoff functions 

 \item $\mathcal{P}(A)$ the set of probability measures on $A$ here the nature randomize the behavior of the players.

\end{itemize}

As we already mentioned the interactions in the model happen based on a function of the empirical distribution $\tau(\bar{\mu}_X^N)$ .

\begin{d1}{Empirical distribution}
 
 Let $X_1,...,X_N$ be independent and identically distributed random variables, with distribution $F(x)=\mu(X_1<x)$. The empirical distribution is defined as
 $$\mu_X^N=F_N(x)=\frac{1}{N}\sum_{i=1}^N \mathit{1}_{X_i\leq x} $$
\end{d1}

\begin{re1}
 We are going to use $\mu$ for the empirical distribution instead of $F$ and save $F$ for other uses as commonly used in the MFGs literature. We will explain more about empirical distributions in MFGs in section 2 and the purpose of this awkward notation will be clear.  
\end{re1}

And as $N\to \infty$ we would like $\mu_X^N$ to  converge to a distribution $\mu$ by a law of large numbers. Indeed it is true by the next theorem

\begin{t1}{Glivenko-Cantelli Lemma}
 
 The empirical distribution converges uniformly to $\mu$ i.e.
 $$ \underset{x\in\mathbb{R}}{sup}|\mu_X^N-\mu| \overset{a.s.}{\to}0  $$
 as $N\to\infty$
\end{t1}

In order to define convergence formally we need to equip $\mathcal{P}(A)$ with a topology, namely the topology of weak convergence ($\mathcal{W}^*$) i.e 

\begin{d1}{Weak convergence}

 Given a sequence $\{\mu_n\}_{n\in \mathbb{N}}$ of probability measures in $\mathcal{P}(A)$ we say that
 $$\{\mu_n\}_{n\in \mathbb{N}} \underset{w}{\to} \mu \in \mathcal{P}(A)$$
  if and only if
 $$\int_Afd\mu_n \to \int_Afd\mu \hspace{2mm} \forall f \in C(A)$$
 as $n\to \infty$
\end{d1}

\begin{t1}
$\hspace{1mm} $

 If A is compact then $(\mathcal{P}(A),\mathcal{W}^*)$ is compact and can be metrized by the Kantorowich-Rubinstein distance 
 \begin{equation}
  d_1(\mu,\nu)=sup\{\int_A f d\mu - \int_A f d \nu) \}
 \end{equation}
 where $f:A\to \mathbb{R}$ bounded Lipschitz continuous
\end{t1}

\begin{proof}
 We start with a $A$ being a compact metric space or a compact subset of a metric space and $C_0^\infty(A)$ the space of continuous functions on $A$ that vanish at infinity equipped with the infinity norm. By the following extension of Riesz representation theorem (tailored for measures) we have that $(\mathcal{P}(A)$ is isometric to $C_0^\infty(A)^*$
 
 \begin{description}
  \item[Theorem (Riesz-Markov-Kakutani)]

  For any positive linear functional $I$ on $C_0^\infty(A)$  there is a unique regular Borel measure $\mu$ on $A$ such that
 $$I(f)=\int_Af(x)d\mu(x)$$
 for all $f$ in $C_0^\infty(A)$
 \end{description}

 From the compactness of $A$ we get that $C_0^\infty(A)=C_b(A)$ and the weak* topology on $C_0^\infty(A)^*$ induced in $(\mathcal{P}(A)$ by the isometry coincide with the weak topology of measures. 
 
 From Alaoglu-Banach the unit ball in $C_0^\infty(A)^*$ is weak* compact and so compact in weak topology of measures. 
 
\end{proof}

\section{Solution of "When does the meeting start?"}
 We are now ready to solve "When does the meeting start?", using the notions from the previous subsections. 
 
 We assume that as $N$ the number of agents approaches infinity a number of simplifications kick in:
 \begin{itemize}
  \item $J_1=...J_N=J$
  \item $X_1=...=X_N=X$
  \item $t_1=...t_N=t $
  \item $(\epsilon_i)_{1\leq i\leq N} \to \epsilon \sim \mathit{N}(0,1) $
  \item $(\sigma_i)_{1\leq i\leq N} \to \sigma \sim \nu$
 \end{itemize}

 The law of large numbers together with the symmetry of the model, provide us a way to reformulate the problem in terms of the representative agent. The first three bullets are just for notational convenience we could also write whatever follows in terms of agent i.   
 
 $$X=t+\sigma \epsilon= t+Z $$
 
 The core of the problem is the distribution of $ \sigma_i \epsilon_i $ (the idiosyncratic shocks) which generate the uncertainty in the model. Since they are independent their distribution $F(\cdot)$ is going to be:
 
 \begin{equation}
  F_Z(z)=\mathbb{P}[Z<z]=\mathbb{P}[\sigma \epsilon <z]=\int_{-\infty}^\infty \nu(x)\Phi(\frac{z}{x})\frac{1}{|x|}dx=\int_{-\infty}^\infty \Phi(\frac{z}{x})\nu(dx)
 \end{equation}

 The next step is to compose the best response of the representative agent to the distribution of actions of the other players. For that reason we notice that the empirical distribution $\bar{\mu}_X^N$ approaches a distribution $\mu$ and $T=\tau(\bar{\mu}_X^N)$ approach $T^*=\tau(\mu)$ as the number of agents goes to infinity and the $BRF(T^*)$ is the solution of the minimization problem:
 
 $$\underset{t\in A}{inf}J(t;T^*) $$

Which comes from the first order condition of (1.1.1) i.e.

$$
J(t;T^*)=\mathbb{E}\big[A(X-t_0)^+ +B(X-T^*)^++C(T^*-X)^+\big] $$
$$=\mathbb{E}\big[A(X-t_0)^++B(X-T^*)^++C(T^*-X)(1-\mathit{1}_{\{T^*<X  \}}) \big] 
$$
$$=A\mathbb{E}\big[(t+Z-t_0)^+\big] +(B+C)\mathbb{E}\big[t+Z-T^*\big]+C(t-T^*)$$
$$
=A\mathbb{E}\big[(t+Z-t_0)\mathit{1}_{\{t+Z-t_0>0  \}}\big] +(B+C)\mathbb{E}\big[(t+Z-T)\mathit{1}_{\{t+Z-T^*>0  \}}\big]+C(t-T^*)$$

now to find the minimum from the first order condition

$$\frac{d}{dt}J(t;T^*)=0$$
$$ A\mathbb{E}\big[\mathit{1}_{\{t+Z-t_0>0  \}}\big] +(B+C)\mathbb{E}\big[\mathit{1}_{\{t+Z-T^*>0  \}}\big]-C=0$$
$$A\mathbb{P}\big[t+Z-t_0>0\big] +(B+C)\mathbb{P}\big[t+Z-T^*>0\big]=C $$
$$A\mathbb{P}\big[Z<t-t_0\big] +(B+C)\mathbb{P}\big[Z<t-T^*\big]=C $$

using (1.5.1) we get an implicit equation of $t$
\begin{equation}
 AF(t-t_0)+(B+C)F(t-T^*)=C
\end{equation}

and this way we have proven the following proposition

\begin{p1}

 If $A,B,C$ are positive constants and $X=t+\sigma \epsilon$ with $\sigma, \epsilon$ as described by (1.5.1) then there exists $\hat{t}$, a unique minimizer of (1.1.1) given by (1.5.2) with $T^*=\tau(\mu)$ being fixed.  
\end{p1}

\begin{proof}
 For uniqueness we have to notice that $F(\cdot)$ given by (1.5.1) is strictly monotone, suppose there are two minimizers $t_1,t_2$ and show that they are identical. 
\end{proof}

The next step to identify a Nash equilibrium is to search for a fixed point in the BRFs. Here we need to be careful because the players interact through the distribution of the states (the time the event begins is a function of the arrival times $T^*=\tau(\mu) $ in the limit).
We are going to define an operator and then use Banach's fixed point theorem.

\begin{p1}
 
 Let 
 $$ \hat{t}:=G(T^*)  $$
 then $G:A\to A$ has a unique fixed point, i.e. 
 \begin{equation}
  G(T^*)=T^* 
 \end{equation}

\end{p1}

\begin{proof}
 $A$ is closed by definition, and map $A$ to itself, all that remain is to show that $G$ is contractive to apply Banach's fixed point theorem. 
 
 Let's use the implicit function theorem on (1.5.2) with respect to $T^*$
 
 \begin{align*}
  \frac{d\hat{t}}{dT^*}=& \frac{(B+C)F'(\hat{t}-T^*)}{AF'(\hat{t}-t_0)+(B+C)F'(\hat{t}-T^*)}\\
  &=\frac{1}{\frac{AF'(\hat{t}-t_0)}{(B+C)F'(\hat{t}-T^*)}+1}\\
  &= \lambda\leq 1
 \end{align*}
because $A,B,C>0$ and $F'$ is nonnegative for rules $\tau(\cdot)$ that satisfy the following properties:
\begin{itemize}
 \item $\forall \mu \hspace{2mm} \tau(\mu)\leq t_0$ the meeting never starts before $t_0$
 \item \textbf{Monotonicity} If $\mu([0,t])\leq \mu'([0,t])$ for all $t\leq 0$ then $\tau(\mu)\geq \tau(\mu')$
 \item \textbf{Sub-additivity} For all $t\leq 0 \hspace{3mm} \tau(\mu(\cdot-t)\geq \tau(\mu)+t$
\end{itemize}

\end{proof}

\section{Differential games and Optimal control}

\section{MFG and Economics}

We end this introductory section by presenting some of the most important ideas in economics  that led to the development of the MFGs theory. 

Technically speaking MFGs are the result of the advances that happened in stochastic control and stochastic differential games during the last thirty years. However, the ideas behind modelling a large number of symmetric agents which independently try to optimize are at least two hundred old. The first well known author  that spoke about a large number of agents that collectively appear one representative agent is Adam Smith\footnote{actually there were other minor authors before Adam Smith that introduced some of the ideas he synthesized in his theory} who used the notion of the "invisible hand" that brings the market into an equilibrium. There is a famous quote that is attributed to him "We don't eat meat by the kindness of the butcher nor bread by the kindness of the baker, it is their personal interest to earn money that guide them to sell us meat or bread." 

Later by the beginning of 20th century the Marginal school appeared in economics and differential calculus they introduced notions like marginal benefit and marginal cost to study agents behavior (firms or consumers) in a single market (partial equilibrium models). This approach is what is usually called microeconomic  where the center the analysis is the single agent. Meanwhile, a different approach was also developed where we could study the aggregate variables of the economy, like total product, inflation, labour etc neglecting the single agent, usually called macroeconomic. 

A first attempt by Leon Walras to provide an explicit analytical model which could combine microeconomic elements and produce laws that govern the whole economy, born the general equilibrium economic models. His approach was not satisfactory and the problem remained partially open until the famous proof by Arrow-Debreu which inaugurated a new era in mathematical economics using abstract analysis techniques. It is worth to mention also Aumann's famous article about a market with a continuum of traders that pushed this line of thinking even further. 

In the same spirit game theory- born out of von Neumann's collaboration with Morgensten- was a mathematical attempt to study in a consistent way human incentives, in situations where they have to take actions. John Forbes Nash initiated the study of games with many players with his famous theorem about existence of equilibrium in mixed strategies and unified game theory with current economic theory.

The history continuous with Rufus Isaac who first studied games in continuous time (differential games) using optimal control methods around fifties, to end with the development of stochastic differential games and finally Mean field games. 
 
In the fourth section where we present a MFG version of a macroeconomic model and implement exactly the way of thinking that mentioned above, to start from the agent's level and end up with a general equilibrium for the whole economy. 

\chapter{Mean Field Game Formulation}

In this chapter we are going to develop our formal definitions about MFGs and Nash equilibrium. We discuss MFGs in a continuous time interval $[0,T]$ with continuous states so that our analysis borrows elements for the theory of stochastic differential games rather than traditional game theory approach. 

We aim to provide functional and conceptual definitions helpful in understanding mean field games modeling in stead of achieving the greatest mathematical generality. Starting from a fairly general setting of a stochastic N-player differential game we motivate the need and usefulness of the mean field games assumptions. (For more information about differential games and stochastic differential games we refer to \cite{Isaacs1965}, \cite{Card2010} and  \cite{BF1984})

\section{General model set-up}

Suppose we start with a finite set of players $P$ with $\#(P)=N$ and each player $i$ can choose an action $a^i_t$ from an action set $A^i$ this action is a functional which can take continuous or discrete values is space usually called state space $X$ for our state variables $X^i_t$s. These processes $X^i_t$ characterize the position of each player $i$ ( the meaning of "position" can vary according to the context of each specific game we study). Each player has also the choice to randomize his behaviour playing a mixed strategy but we will not consider this case here.  

Furthermore, each agents has a functional $J^i$ as his cost or benefit criterion which he is interested to optimize. 

 To conclude we assume also $(\Omega,\mathcal{F},\{\mathcal{F}_t\}_{t\geq 0},\mathbb{P}) $ to be a standard filtered probability space on which we can defined an N-dim (same dimension as our state process) Brownian Motion $W_t)$ such that $\{\mathcal{F}_t\}_{t\geq 0}$ is generated by $W$ augmented with all the $\mathbb{P}$-null sets in $\mathcal{F}$.

We warp everything as our terminology:

\begin{d1}{Terminology}
\begin{enumerate}
  \item
  $P$: the set of players (agents)\\
  $\#(P)=N$ the number of players
  \item
  $X$: the state space, can be a metric space or a subspace of a metric space (usually assumed compact)
  
  \item
  $A^i$:  the set of actions for player i\\
  $A=A^1\times... \times A^N$\\
  $a=(a_1,...a_N)\in A$ is an action profile where $a_i$ is the action the individual players take and $a_{-i}$ the action profile including every player's action except $i$'s  $a_{-i}=(a_1,...a_{i-1},a_{i+1},...a_N)$ 
  
  $A^i_{adm}:=\{a^i_t:[0,T]\to X| a^i_{(\cdot)} \text{ arbitrary, admissible function} \} $: The set of admissible strategies for player $i$
  
  \item
  $\mathcal{U}$: the set of payoff functions $J_i:A_i\to \mathbb{R}$ (Cost/Benefit criterion) for player $i$
  \item
  $\mathcal{P}(A^i)$ The set of probability measures on $A^i$ 
  \item
  $\mathcal{M}=\mathcal{P}(A^1)\times ... \times \mathcal{P}(A^N)   $ The Cartesian product of sets of probability measures  on individual action sets, are called mixed strategies. 
 \end{enumerate}

\end{d1}

Here because we talk about MFGs in continuous time i.e.  mean field stochastic differential games we does not distinguish between actions and strategies and use the terms interchangeably. An action functional $a_t$ represent the decisions of the agents at each moment. Appendix A and B explain more about their role in optimal control problem. 

\subsection{Agents playing pure strategies with noise}

In the case where each agent plays a pure strategy then we can formulate the game as an N-player optimal control problem where everyone interacts with each other through their controls and/or their states. We model the state variables $X_t$ to evolve according to a system of coupled SDEs, where the decisions of the agents and the distribution of the controls and/or their states provide the coupling of the SDEs. This yields a stochastic optimal control problem for each agent, given the distribution of the actions and/or the states of the rest of the players. To make everything more precise we assume the following problem for each agent $i$:

\newtheorem*{prob1}{Individual agent's problem playing pure strategies}

\begin{description}
 \item[Individual agent's problem playing pure strategies] 
 
 \begin{equation}
  \underset{a^i_{(\cdot)}\in A^i_{adm}}{inf} J^i(a_{(\cdot)};\nu^{-i}_{(\cdot)})=\underset{a^i_{(\cdot)}\in A^i_{adm}}{inf} \int_0^T \mathbb{E}\big[ f(X^i_t,a^i_t,\nu^{-i}_t) dt + g(X^i_T,\nu^{-i}_t) \big]
 \end{equation}
subject to
\begin{equation}
\left\{
\begin{split}
 &dX^i_t=b(t,X^i_t,a^i_t,\nu^{-i}_t)dt+\sigma(t,X^i_t,a^i_t,\nu^{-i}_t)dW^i_t\\
 &X^i_0=\xi^i\in L^2\\
 \end{split}
 \right.
\end{equation}
where
\begin{equation}
 \nu_t=(\nu^1_t,...,\nu^N_t)=
 \begin{cases}
                              \nu^1_t=\mathcal{L}(X^1_t,a^1_t)\\
                              ...\\
                              \nu^N_t=\mathcal{L}(X^N_t,a^N_t)\\
                              \end{cases}
\end{equation}

\end{description}

$\{\nu_t\}_{0\leq t \leq T}$ represent the flow of probability measure, i.e each coordinate of the N-tuple $(\nu^1_t,...,\nu^N_t)$ is a flow of distributions for each player to interact with each other.

\begin{re2}

\begin{itemize}
$\hspace{0,5mm}$
\item 
Action sets can be \textbf{finite} or \textbf{infinite}. 
\item 
The class of models described by (2.1.1), (2.1.2) are called Second-order MFGs in the literature because the dynamics are described SDEs and using dynamic programming principle (DPP) we end up in a second order Hamilton-Jacobi-Bellman (HJB) equation. If instead we use ODEs for (2.1.2) we end up with First-order HJB equation and so they are called First order MFGs. 
\item
From a game theoretic point of view it is very natural that the agents interact through their controls (actions) and their states $\nu^i_t=\mathcal{L}(X^i_t,c^i_t)$ with $\mathcal{L}$ representing the common law of the state and the control. If the players are interacting only according to states (as in "When does the meeting start?") then $\nu^i_t=\mathcal{L}(X^i_t)$. Historically the first models that were developed (McKean 1968) were of interacting states both because of their simplicity and their connection with statistical physics. In the next few subsections we are going to follow this line of thinking indeed and also explain more about the flow of probability measures and why it is a natural concept to describe large scale models as mentioned already in the introduction. 
\item
 The random noise $W^i_t$ is independent for each player $i$, there is a possible extension in our modelling by adding a $W^0_t$ common for all the players, these models are known as Mean Field Games with common noise and are considerably more difficult and require different treatment than what we are going to present here. 

\end{itemize}
\end{re2}

Solving his optimization problem each agent can construct his best response function, given the distributions of the rest of the players. The intersection of all BRFs is the Nash equilibrium of the game. 

As we have already seen in the introduction this problem is very difficult to solve and this is where MFGs kick in. We can achieve considerable simplifications if we assume, symmetry and that the number of agents goes to infinity.  

%\section{Agents playing mixed strategies}

\section{Limiting behaviour of large systems}

In this subsection we present the ideas that opened the way for development of MFGs. The situation, when the number of players goes to infinity, is of great importance for MFGs and we are going to borrow the so called propagation of (molecular) chaos from statistical physics to describe it.

We will start by describing, the simplest case, about what is called as a hard sphere gas, where everything is deterministic and governed by ODEs and gradually extend the framework to interacting diffusions which will be described be SDEs and then draw an analogy with a game of interacting players. This way we will give some intuition for the complex system (2.1.1)-(2.1.3). 

\subsection{Boltzmann's theory of hard sphere gases}

The simplest way we can imagine the molecules of a dilute gas is, as small hard spheres of some radius $r$ and mass $m$ that are moving randomly and can collide. Let's assume that they live in a position-velocity space $S\subset \mathbb{R}^6 $ (generally it can be any finite-dimensional separable metric space) and a $N$-particle system is a point in $S^N$ (Cartesian product). Moreover their dynamics are Markovian in a sense that the future position of the system only depends on its current position, this way we can define transition functions.   Let $N$ be the number of molecules of the gas and define the density (i.e. the number of molecules per unit volume of $S$) as  $ f(x,u,t) $ where $x$ is the position and $u$ is the velocity. While
$$ \frac{1}{N} \int_U\int_X f(x,u,t)dxdu$$
is the proportion of molecules which, at time $t$ are located in a region of space $X$ and have velocities in $U$.

Now we are ready to state Boltzmann's equation for the evolution of $f(x,u,t)$ (derivation of the equation escapes the scope of this text but we refer the reader to the original work of Boltzmann \cite{Boltzmann1995}

\begin{equation}
\frac{\partial f}{\partial t}=\frac{u}{m}\nabla_xf+C[f]
\end{equation}

where $\frac{u}{m}\nabla_xf $ gives the rate of change due to streaming and $C[f]$ is the collision operator applied on $f$. which gives the rate of change of the density due to collisions which are governed by principles of momentum and energy conservation. We need further assumptions to describe the collision operator $C[f]$, but intuitively speaking  we can say that depends upon the rate at which collisions are happening and the post-collision velocities of the molecules.

\subsubsection{Vlasov's theory of plasmas}
Anatoly Vlasov proposed his theory about plasmas in 1938 and published it as a monograph "Theory of Vibrational Properties of an Electron Gas and Its Applications" in 1945. The primary focus of Vlasov was to describe plasmas where the ions never collide and instead have long range interactions which Boltzmann's equation cannot describe properly.

We adopt the same setting as before with the extra assumption that all the particles are of the same kind (for example electrons). Let $F(x)$ be the force that a particle at the origin would exert at a particle at x. Since the interactions cover the whole space S they generate a force field $F_f(x)$ (again the technical details about particles escape the scope of this text)
\begin{equation}
 F_f(x)=\int_SF(x-x')f(x',u',t)dx'du'
\end{equation}
The particle density changes through the motion of particles subject to the force filed $F_f(x)$ and Vlasov's equation for the evolution of density is

\begin{equation}
\frac{\partial f}{\partial t}=\frac{u}{m}\nabla_xf+\frac{1}{m} F_f(x)\nabla_uf
\end{equation}

\subsection{Propagation of chaos in Boltzmann's and Vlasov's theory}
So far,  we have presented the basic kinetic theories for gases and plasma, now we would like to introduce also the idea of molecular chaos propagation and use it to better understand the continuum limit of MFGs.

 Suppose we have an $N$-particles system and a probability measure is assigned to each particle so we get a sequence of probability measures $\{\mu_i \}_{i=1}^N $. We think of the measures as giving the joint probability distributions of the first i particles, for example $\mu_3$ gives the joint distribution of particles 1,2,3.

\begin{d1}{Propagation of chaos}

We say that a sequence of probability measures $\{\mu_i \}_{i=1}^N $ is $\mu$-chaotic if $k$ coordinates, become independent and tend to $\rho$ as $N$ goes to infinity i.e. for any $k\in \mathbb{N}$ and $g_1(s),...,g_k(s)\in C_b(S)$
$$\underset{N\to \infty}{lim}\int_S g_1(s_1)...g_k(s_k)\mu_N(ds_1...ds_k)=\displaystyle\prod_{i=1}^k\int_Sg_i(s)\mu(ds) $$

\end{d1}

To elaborate more on the idea of molecular chaos propagation we will discuss the case of the Vlasov equation and we will show she propagates chaos.

We assume the same setting as previous section with the extra assumptions that $F: S \to S $ be bounded and Lipschitz and  we define a deterministic $N$-particle process in $S$ for each $N$

\begin{equation}
 \left\{
 \begin{split}
  &\frac{d}{dt}x_i^N(t)=u_i^N(t)\\
  &\frac{d}{dt}u_i^N(t)=\frac{1}{N}\displaystyle\sum_{i=1}^NF(x_i^N-x_j^N)
 \end{split}
\text{for $i=1,...N$}\right.
\end{equation}

as shown in \cite{BH1977} if the initial conditions $x_i^N(0),u_i^N(0)$ for $i=1,...,N$ are such that:
$$\frac{1}{N} \sum_{i=1}^N\delta_{(x_i^N(0),u_i^N(0))}\to \mu_0\in \mathcal{P}(S) $$
then for $t>0$
$$\frac{1}{N}\sum_{i=1}^N\delta_{(x_i^N(t),u_i^N(t))}\to \mu_t\in \mathcal{P}(S) $$
where $\mu_t$ is the weak solution at time t of the Vlasov equation

\begin{equation}
\left\{
\begin{split}
 &\frac{\partial f}{\partial t}+\frac{u}{m}\nabla_xf+\frac{1}{m} F_f(x)\nabla_uf \\
 &F_f(x)=\int_SF(x-x')f(x',u',t)dx'du'\\
 &\mu_0=f(x,u,0)dxdu
\end{split}
\right.
\end{equation}

Thus this $N$-particle system propagates chaos.

\subsection{Interacting Diffusions}

The deterministic particle system can be generalized to interacting diffusions, McKean in his article "Propagation of Chaos and a class of nonlinear parabolic equations" \cite{McK1967} initiated the study of those systems. %This class of diffusions are Markov processes, solutions of specific SDEs which are now called McKean-Vlasov.%

Suppose we have $N$ particles, each one is making a diffusion in a $d$-dim space, the drift and the volatility of their movement are affected by the empirical distribution of the rest $N-1$ particles, but are common for every particle. 

\begin{equation}
 \left\{
 \begin{split}
  &dX_t^i=\{\frac{1}{N}\displaystyle\sum_{j=1}^nb(X_t^i,\bar{\mu}_t^N)\}dt+\{\frac{1}{N}\sum_{j=1}^N\sigma(X_t^i,\bar{\mu}_t^N)\}dW_i \text{ for $i=1,...,N$}\\
  &\bar{\mu}_t^N=\frac{1}{N}\sum_{j=1}^N\delta_{X_t^j<x}
 \end{split}
\right.
\end{equation}

where $b:\mathbb{R}^d \times\mathbb{R}^d \to \mathbb{R}^d $ and $\sigma:\mathbb{R}^d \times\mathbb{R}^d \to \mathbb{R} $ bounded and Lipschitz and $X_i^n$ with values in $\mathbb{R}^d$. The Wiener processes $W_i$ are taken to be independent of each other and of the initial conditions $X_1^N(0),...,X_N^N(0)$

McKean in his article assumes that volatility is constant and equals $1$ and that drift term is given by:

\begin{equation}
 b(X_t^i,\mu_t^N)=\int \bar{b}(X_t^i,y)\mu_t^N(dy)=\frac{1}{N}\sum_{j=1}^N b(X_t^i,X_t^j) 
\end{equation}

with the last equality given by the fact that $\bar{\mu}_t^N $ is an empirical distribution. This way we arrive in:

\begin{equation}
 dX_t^i=\{\frac{1}{N}\displaystyle\sum_{j=1}^N\bar{b}(X_t^i,X_t^j)\}dt+dW_i \text{ for $i=1,...,N$}
\end{equation}

 where $\bar{b}:\mathbb{R}^d \times\mathbb{R}^d \to \mathbb{R}^d $ and $\bar{\sigma}:\mathbb{R}^d \times\mathbb{R}^d \to \mathbb{R} $ bounded and Lipschitz and the rest as before. 

In \cite{McK1967} he proves the following theorem:

\begin{t1}{Propagation of chaos for diffusions}
 
If the particles are initially stochastically independent but with common distribution $\mu_0$, then the sequence of n-particle joint distributions at time $t$ is $\mu_t$-chaotic, $\mu_t$ being the (weak) solution at time t of the nonlinear McKean-Vlasov equation

\begin{equation}
\left\{
\begin{split}
 &\frac{\partial}{\partial t}f_t=-\nabla[V_ff_t]+\frac{1}{2}\Delta f_t\\
 &V_f(x)=\int_{\mathbb{R}^d}\bar{b}(x,x')f_t(x')dx'\\
 &f_0(x)dx=\mu_0
\end{split}
\right.
\end{equation}
where the subscript $t$ in $f_t$ is used to stress the connection with $\mu_t$, not to be mistaken by a time derivative. 
\end{t1}

From Sznitman \cite{Sznitman1989} we get also an alternative statement of the theorem. As $N \to \infty$, $X_t^i$ has a natural limit $\bar{X}_t^i$. Each $\bar{X}_t^i$ will be an independent copy of the nonlinear process $\bar{X}_t$.

\begin{t1}{(Sznitman)}

 There is existence and uniqueness both trajectorial and in law for the nonlinear process $\bar{X}_t$:
 
 \begin{equation}
 \left\{
 \begin{split}
  &d\bar{X}_t=\{\int_{\mathbb{R}^d}\bar{b}(\bar{X}_t,y)d\mu_t\}dt+d\bar{W}\\
  &\bar{X}_0=x_0 \text{ $\mu_0$-distributed, $\mathcal{F}_0$-measurable random variable}\\
  &\mu_t \text{ the law of $\bar{X}_t$}\\
 \end{split}
\right.
 \end{equation}
\end{t1}

\begin{proof}
 To begin let us assume as usual that $C:=C([0,T];\mathbb{R}^d)$ is the space of continuous functions on $[0,T]$ with values in $\mathbb{R}^d$
 and $\mathcal{P}(C)$ the space of probability on $C$. We equip $\mathcal{P}(C)$ with the Kantorowich-Rubinstein metric then $\mathcal{P}(C)$ is complete as we have already discussed in the introduction. We take $T>0$ and define $\Phi$ as the map that associates to $\mu\in\mathcal{P}(C)$ the law of the solution of 
 \begin{equation}
 \left\{
 \begin{split}
  &dX_t=\{\int_Cb(X_t,w_t)d\mu(w)\}dt+dW_t \text{ t$\leq$T}\\
  &X_0=x_0
 \end{split}
\right.
 \end{equation}
The law does not depend on the specific choice of the space $\Omega$.
If $\{X_t\}_{t\leq T}$, is a solution of (2.2.10), then its law on $C$ is a fixed point of $\Phi$ , and conversely if $\mu $ is such a fixed point of $\Phi$ (2.4.11) defines a solution of (2.2.10) up to time $T$. 

For the fixed point argument using Banach's theorem we refer to \cite{Sznitman1989}.  
\end{proof}

To connect the  nonlinear process with the nonlinear PDE (2.2.9) we use Ito's formula for $f \in C_b^2(\mathbb{R}^d) $. 
$$ f(\bar{X}_t)=f(\bar{X}_0)+f'(\bar{X}_t)dW_t+\{(\frac{1}{2}\Delta f+\int_{\mathbb{R}^d}\bar{b}(\bar{X}_t,y)u_t(dy)\nabla f(\bar{X}_t)\}dt$$

and assuming it is a true martingale we set $dt$ part equal to zero and get (2.4.9)

\begin{re1}
$\hspace{0,5mm}$
 If we set $W_t\equiv 0$ in McKean's model we get Vlasov's equation for plasmas.
\end{re1}

\subsection{Propagation of chaos and MFGs}

Interacting diffusions can can extend to stochastic differential games if we grant the freedom of choice to every particle $i$ and rename the particles as agents or players. In the particular case we study each player can decide about his drift, which can affect his position $X^i_t$ and as a consequence the empirical distribution of the states, $\mu^N_t$, and these decisions are thought to be measurable functions, which we are going to define in detail later. Furthermore we introduce a criterion for the decisions.

For $i=1,...,N$
\begin{equation}
 J^i(a^1_t,...a^N_t)=\mathbb{E}\bigg[\int_0^Tf(X^i_t,\mu^N_t,a^i_t)dt+g(X^i_T,\mu^N_T)\bigg]
 \end{equation}
 
subject to

 \begin{equation}
  dX_t^i=\{\frac{1}{N}\displaystyle\sum_{j=1}^N\bar{b}(X_t^i,a^i_t,X_t^j)\}dt+dW_i  
 \end{equation}

It is only natural to extend the previous theorems under our current set-up, since we keep the Lipschitz assumption about $\bar{b}$, we can repeat the proof with no changes. 

We let $N\to\infty$ and  $X_t^i\to \bar{X}_t^i$. Again each $\bar{X}_t^i$ will be an independent copy of the nonlinear process $\bar{X}_t$ as before and also we want to provide a limit for $J^i$. For this reason we are going to investigate in the next subsection symmetric functions of many variables.

\section{Symmetric functions of many variables}

Considering (2.2.12) it is not precise the way it is written, $J^i(a^1_t,...,a^N_t)$ is a function that depends on $N$ variables and on the other hand we have functions of empirical measures. The definition would hold true only for $J^i(a^i_{(\cdot)};\mu^N_{(\cdot)})$ a function that depends upon individual's $i$ decision and the empirical distribution of the states $\mu^N_t$ of the rest of the players but what is special about this dependence is that we have identical players and diffusions. So this functional should enjoy some properties which can lead us to define a limit when $N\to\infty$, in addition, we would like to be able to approximate functions as $J^i(a^1_t,...,a^N_t)$ by functions of measures to make (2.2.12) precise.

\begin{d1}{Symmetric functional}

A function $J^n:Q^n\to \mathbb{R}$ with $Q^n$ being compact is called symmetric iff 
\begin{equation}
 J^n(c_1,...,c_n)=J^n(c_{\pi(1)},...,c_{\pi(n)}) \text{ for every permutation $\pi$ on $\{1,...,n\}$ }
\end{equation}
\end{d1}

\begin{t1}
 
 For each $n\in \mathbb{N} $, let $u^n:Q^n\to \mathbb{R}$ be a symmetric function of its $n$ variables. \textbf{We assume}:
 
 \begin{enumerate}
 
  \item (\textit{Uniform boundedness}) There some $C>0$ such that
   \begin{equation}
   ||u^N||_{L^\infty(Q^n)}\leq C 
   \end{equation}
  \item (\textit{Uniform continuity})
  This a modulus of continuity $\omega$ independent of $n$ such that 
  \begin{equation}
   |u^n(X)-u^n(Y)|\leq \omega d_1(\mu_X^n,\mu_Y^n) \hspace{2mm} \forall X,Y\in Q^n,\forall n\in \mathbb{N}
  \end{equation}
where $\mu_X^n=\frac{1}{n}\sum_{i=1}^n \mathit{1}_{x_i}$ and $\mu_X^n=\frac{1}{n}\sum_{i=1}^n \mathit{1}_{x_i}$ if $X=(x_1,...,x_n)$ and $Y=(y_1,...,y_n)$ by $d_1$ we mean the Kantorowich-Rubinstein distance in $\mathcal{P}(Q^n)$
  
 \end{enumerate}
 
\textbf{Then} there is a subsequence $u^{n_k}$ of $u^n$ and a continuous map $U:\mathcal{P}(Q)\to \mathbb{R}$ such that: 
\begin{equation}
 \underset{k\to \infty}{lim}\underset{X\in Q^n}{sup} |u^{n_k}(X)-U(\mu_X^{n_k})|=0
\end{equation}

\end{t1}

Before we give the proof some remarks are one the way. 

\begin{re2}

 \begin{enumerate}
  \item 
  This map $U$ is going to play the role of our payoff functional in what follows, we need its domain to be the space of probability measures since the payoff functional depends on the decisions of all players. i.e. the empirical distribution $\mu_t^n$
  \item
  The assumptions (1),(2) are essential to make use of the Arzela-Ascoli theorem which is going to give us the uniformly convergent subsequence
 \end{enumerate}

\end{re2}

\begin{proof}
 Since we know that $\mathcal{P}(Q^n)$ is compact and complete we want to exploit this and construct the map $U:\mathcal{P}(Q)\to \mathbb{R}$ which is going to satisfy the assumptions for Ascoli-Arzela theorem. 
 
 Let us begin by defining $(U^n)_{n\geq 1}$ on $\mathcal{P}(Q^n)$ by:
  \begin{equation}
   U^n(\mu)=\underset{X\in Q^n}{inf} \big[u^n(X)+ \omega d_1(\mu_X^n,\mu) \big] \hspace{2mm} \mu\in\mathcal{P}(Q)
  \end{equation}

We need to prove that these functions qualify for Ascoli-Arzela. Boundedness is checked easily since $(u^N)_{n\geq 1}$ are bounded from assumption 1 and $\mathcal{P}(Q)$ is compact so $\forall \mu,\nu\in \mathcal{P}(Q)$ distance $d_1(\mu,\nu)$ is bounded. 

Also easily from def (2.4.5) together with assumption 2, we can show that these functions $(U^n)_{n\geq 1}$ extend the original $(u^n)_{n\geq 1}$ to $\mathcal{P}(Q)$ meaning that 

$$U^n(\mu_X^n)=u^n(X) \text{ for any $X\in Q^n$} $$

Furthermore, let us show that $(U^n)_{n\geq 1}$ have $\omega$ for modulus of continuity on $\mathcal{P}(Q)$ i.e. are equicontinuous. 
Indeed if $\mu,\nu\in\mathcal{P}(Q)$ and if $X\in Q^n$ is $\epsilon$-optimal in the definition of $U^n(\nu)$, then

 $$U^n(\mu) \leq u^n(X)+\omega d_1(\mu_X^n,\mu)$$
 $$= U^n(\nu)-\omega d_1(\mu_X^n,\nu)+\epsilon+\omega d_1(\mu_X^n,\mu)$$
 $$\leq U^n(\nu)+\epsilon -\omega d_1(\mu_X^n,\nu) +\omega \big(d_1(\mu_X^n,\nu)+d_1(\nu,\mu)\big)$$
 $$=U^n(\nu)+\omega d_1(\mu,\nu)  +\epsilon$$

Now since $\mathcal{P}(Q)$ is compact Arzela-Ascoli gives the existence of a subsequence $(n_k)_{k \geq 1}$ for which $(U^{n_k})_{k \geq 1}$ converges uniformly to a limit $U$ and since $u^{n_k}(X)=U^{n_k}(\mu_X^{n_k}) $ for any $X\in Q^n$ we get (2.3.4)

\end{proof}

By the means of the above theorem we can approximate $J^i(a^1_t,...,a^N_t)$ by functions of measures and  construct limits, in the sense that when $N\to \infty$ then $(J^i)_{i=1}^N\to \bar{J}$ with $\bar{J}$ given by:

$$
 \bar{J}(\bar{a}_{(\cdot)};\mu_{(\cdot)})=\mathbb{E}\bigg[\int_0^Tf(\bar{X}_t,\mu_t,\bar{a}_t)+g(\bar{X}_T,\mu_T)\bigg]
$$

To make a brief summery up until now, we would like to think of the players as if they were ions in a plasma where the particles never collide. The key observation is that we cannot apply directly Boltzmann's theory of hard sphere gases because the gravity force usually is modelled as an inverse square potential with a singularity at zero and this would be unrealistic for systems of interacting players who never collide in a physical sense. So in case of deterministic games we would stick with Vlasov's theory and for stochastic differential games we would go with McKean's theory for interacting diffusions. Then we need also to define a limit for the sequence of criterion functions $(J^i)_{i=1}^N$ to describe fully  the situation at infinity. In other words we are interested in the finite game:

\begin{description}
 \item[Individual agent's problem playing pure strategies (simplified version)] 
 For each $i=\{1,...,N\}$ 
 
 \begin{equation}
  \underset{a^i_{(\cdot)}\in A^i_{adm}}{inf} J^i(a_{(\cdot)};\bar{\mu}^{-i}_{(\cdot)})=\underset{a^i_{(\cdot)}\in A^i_{adm}}{inf} \int_0^T \mathbb{E}\big[ f(X^i_t,a^i_t,\bar{\mu}^{-i}_t) dt + g(X^i_T,\bar{\mu}^{-i}_t) \big]
 \end{equation}
subject to
\begin{equation}
\left\{
\begin{split}
 &dX^i_t=b(t,X^i_t,a^i_t,\bar{\mu}_t^N)dt+dW^i_t\\
 &X^i_0=\xi^i\in L^2\\
 \end{split}
 \right.
\end{equation}
where
\begin{equation}
 \bar{\mu}_t^N=\frac{1}{N}\sum_{j=1}^N\mathit{1}_{X_t^j<x}
\end{equation}

\end{description}

As we let $N\to\infty$ we end up with: 

\begin{description}
 
\item[Representative agent's problem playing pure strategies]

 \begin{equation}
 \underset{\bar{a}_{(\cdot)}\in \bar{A}_{adm}}{inf}\bar{J}(\bar{a}_{(\cdot)};\mu_{(\cdot)})=\underset{\bar{a}_{(\cdot)}\in \bar{A}_{adm}}{inf}\mathbb{E}\bigg[\int_0^Tf(\bar{X}_t,\mu_t,\bar{a}_t)+g(\bar{X}_T,\mu_T)\bigg]
 \end{equation}

 subject to
 
 \begin{equation}
 \left\{
 \begin{split}
  &d\bar{X}_t=\{\int\bar{b}(\bar{X}_t,\bar{c}_t,y)d\mu_t\}dt+d\bar{W}\\
  &\bar{X}_0=x_0 \text{ $\mu_0$-distributed, $\mathcal{F}_0$-measurable random variable}\\
  &\mu_t=\mathcal{L}(\bar{X}_t)\\
 \end{split}
\right.
 \end{equation}

 where $\bar{X}_t$ is the nonlinear defined earlier. 
\end{description}

Now we ready introduce the concept of Nash equilibrium for MFGs.

\section{Nash equilibriums for MFGs}
 
 Suppose now that we have a well defined Mean Filed Game in the sense of subsection 2.1, with $N$ players, for each $i=\{1,...,N\}$ we have a criterion function $(J^i)N$ that enjoy the properties of theorem 2.3. We are interested to study the strategic situation that each agent optimize his payoff functional, because as we discussed in the introduction it is not feasible for all the agents to get their global maximum or minimum we are searching for a situation that each one plays his best response to other players actions, and the system reaches an equilibrium where no one has the incentive to deviate. This is the Nash equilibrium as we introduced it in the previous section for our simple games. 
 
 \begin{d1}{Nash equilibrium}
 
 An action profile $a^*_t\in A$ is called a Nash equilibrium for a fixed time $t\in [0,T]$ if and only if for every player i
$$J^i(a^*_t)\leq J^i(a^i_t,a^{*,-i}_t) \hspace{1mm} \forall a^i_t\in A^i  $$
where $J^i$ is the payoff functional of player i. 

 \end{d1}
 
 The above definition is for fixed time and the action profile $a^*$ is a vector of functions $a^{*,i}_t$. For the various forms this control function can take and the different information structures the players can depend upon to adapt their strategies, we distinguish between the following cases:
 
 \begin{itemize}
  \item \textbf{Open loop equilibrium}\\ 
  Apart from the initial data $X_0$ player $i$ cannot make any observation of the space of the system. There is no feedback, and thus this is an open loop  process. 
  $$a^i_t=c^i(t,X_0,W_{[0,t]},\mu^N_t) $$
  For measurable deterministic functions $c^i$, $i=\{1,...,N\}$
  
  \begin{d1}{Open loop Nash equilibrium}
  
  An action profile $a^*\in A$ is an open loop Nash equilibrium, if whenever a player $i\in\{1,...,N\}$ uses a different strategy $a^i_t=c^i(t,X^i_0,W^i_{[0,t]},\mu^N_t)$ from $a^{*,i}_t$  while the other players keep using the same, then $J^i(a^*_t)\leq J^i(a^i_t,a^{*,-i}_t)$ 

  \end{d1}

  \item \textbf{Closed loop equilibrium}\\ 
  Here the player $i$ can also observe the state space and can use this information to update his strategy (in form of feedback) and thus this is a closed loop process. However he has no additional information about the strategy of the other players.  
  $$a^i_t=c^i(t,X^i_{[0,t]},\mu^N_t) $$
  For measurable deterministic functions $c^i$, $i=\{1,...,N\}$ and $\{X^i_t\}_{0\leq t\leq T}$ the solution of the state dynamics SDE. Also we notice that this is an implicit, path-dependent form since $X^i_t$ depends also on the controls. 
  
  In case the player can observe the whole space we say that he has perfect or complete observability and in the case he can observe only the states of some players e.g. being close to his position, we say that he has partial observability and we take this in to account in the definition of $c^i$ 
  
  \begin{d1}{Closed loop Nash equilibrium}
   
   Suppose $\{X^*_t\}_{0\leq t\leq T}$ is the solution of the state dynamics SDE (2.3.13) when we use the actions $a^*=(c^{*,1}(t,X^{*,i}_{[0,t]},\mu^N_t),...,c^{*,N}(t,X^{*,i}_{[0,t]},\mu^N_t)$.\\
   An action profile $a^*_t\in A$ is a closed loop Nash equilibrium, if whenever a player $i\in\{1,...,N\}$ uses a different strategy $a^i_t=c^i(t,X^i_{[0,t]},\mu^N_t)$ while the rest continue to use $b_t=(c^{*,j}(t,X^i_{[0,t]},\mu^N_t)$, $\forall j\neq i$ but with  $\{X^i_t\}_{0\leq t\leq T}$ the solution of the state dynamics SDE (2.3.13) when we use the actions $a_t=(a^i_t,b_t)$. Then $J^i(a^*_t)\leq J^i(a_t)$ 
  \end{d1}

  In a closed loop Nash equilibrium a change in a player's strategy  will result in an change in the state process. The rest of the players will also adjust their actions because their payoff functional changes. They will keep using the same $c^{*,j}$, $\forall j\neq i$ to compute their controls but according to the new path of the state process $X_t$.
  
  \item \textbf{Markovian Nash equilibrium}\\ 
Markovian Nash equilibrium is a special class of close loop controls which depend only in the current value of $X_t$ instead of the whole path $X^i_{[0,t]}$. 
$$a^i_t=c^i(t,X^i_t,\mu^N_t) $$
For measurable deterministic functions $c^i$, $i=\{1,...,N\}$

\begin{d1}{Markovian Nash equilibrium}
   
   Suppose $\{X^*_t\}_{0\leq t\leq T}$ is the solution of the state dynamics SDE (2.3.13) when we use the actions $a^*=(c^{*,1}(t,X^*_t,\mu^N_t),...,c^{*,N}(t,X^*_t,\mu^N_t)$.\\
   An action profile $a^*_t\in A$ is a close loop Nash equilibrium, if whenever a player $i\in\{1,...,N\}$ uses a different strategy $a^i_t=c^i(t,X^i_t)$ while the rest continue to use $b_t=(c^{*,j}(t,X^i_t,\mu^N)$, $\forall j\neq i$ but with  $\{X^i_t\}_{0\leq t\leq T}$ the solution of the state dynamics SDE (2.3.13) when we use the actions $a_t=(a^i_t,b_t)$. Then $J^i(a^*_t)\leq J^i(a_t)$ 
  \end{d1}
  
  Using Makrovian action profiles instead of state insensitive adapted processes (open loop controls) will affect the dependence upon the state variable in the third section where we are going to search for solutions for MFGs.
  
 \end{itemize}
 
\subsection{Limits of Nash equilibrium}
 
 Now we turn to study the situation when we have a Nash equilibrium and send $N$ the number of agents at infinity, for this subsection we assume that the time is frozen at $t$ and everything refers to this particular moment and for this reason we will drop $t$ from our notation for now. Once we have reach the Nash equilibrium, time is not significant any more as we will argue later and in the third section where we are going to solve MFGs we will examine the time frame until we reach Nash equilibrium and see how the various forms of the control functions and Nash equilibriums described earlier affect the solution of the game.
 
 A natural question to ask is:
 "\textbf{Whenever we have a Nash equilibrium for the $N$-player game, is this a Nash equilibrium for the infinite game also?}"
 
 The answer is positive and was given by Lions in \cite{LL2007}, in the form of the following theorem:

\begin{t1}
 
 Assume that $a^*$ is a Nash equilibrium for the game (2.3.6-2.3.8). Then up to a subsequence, the sequence of empirical measures of actions $\{\nu^i\}_i^N$ (different than $\{\mu^i\}_i^N$ the empirical measure of $(X^i)_{i=1}^N$) converges to a measure $\hat{\nu}\in\mathcal{P}(A)$ such that:
 
 \begin{equation}
  \int_A J(y,\hat{\nu})d\hat{\nu}(y)=\underset{\nu^i\in\mathcal{P}(A^i)}{inf} \int_A J(y,\hat{\nu})d\nu^i(y) 
 \end{equation}

\end{t1}

 \begin{proof}
From our game definition we have a sequence $\{J^i(a^1,...,a^N)\}_{i=1}^N$ which fulfils the requirements of Theorem 2.3 and can be approximated (up to a subsequence) by an empirical distribution of the original arguments i.e. the actions.
 $$J^i(a^1,...,a^N)=J^i(a^i,\nu^{N-1}) \text{ for any } i\in\{1,...N\}$$
 where
 $$\nu^{N-1}=\frac{1}{N-1}\sum_{1\leq i\neq j\leq N}\delta_{a^j}$$
  and remembering (2.4.5) for $\nu^N\in\mathcal{P}(A^N)$ and fixed $a^i$ we write:
  \begin{equation}
 J^i(a^i,\nu^N)=\underset{a^{-i}\in A^{N-1}}{inf}\big[J^i(a^i,a^{-i})+\omega d_1(\nu^{N-1},\nu^N)\big]  
  \end{equation}
  
  where
  
  \begin{align*}
   d_1(\nu^{N-1},\nu^N)&=d_1(\nu^{N-1},\frac{1}{N}\sum_{1\leq i\leq N}\delta_{a_i})\\
   &=d_1(\nu^{N-1},\frac{1}{N}(\sum_{1\leq i\neq j\leq N-1}\delta_{a_j}+\delta_{a_i}))\\
   &=d_1(\nu^{N-1},\frac{N-1}{N}\nu^{N-1}+\frac{1}{N}\delta_{a_i})
  \end{align*}
  
  by definition of the Kantorowich-Rubinstein 
  
  \begin{equation}
  \begin{split}
   &d_1(\nu^{N-1},\frac{N-1}{N}\nu^{N-1}+\frac{1}{N}\delta_{a_j})=\frac{f(\delta_{a_i})}{N}\leq \frac{C}{N} \\
   &\text{ for any $f$ bounded and Lipschitz}
  \end{split}
  \end{equation}

  using the definition 2.4 of Nash equilibrium 
  
  $$ J^i(a^{*,i},\hat{\nu}^{N-1})\leq J^i(a^i,\hat{\nu}^{N-1})$$
  or equivalently 
  \begin{equation}
 J^i(a^{*,i},\hat{\nu}^{N-1})=\underset{\nu^i\in\mathcal{P}(A^i)}{inf}J^i(\nu^i,\hat{\nu}^{N-1})= \underset{\nu^i\in\mathcal{P}(A^i)}{inf} \int_A J^i(a,\hat{\nu}^{N-1})d\nu^i(a)   
  \end{equation}
 
 It is obvious that the r.h.s. of (2.4.3) has its minimum at $\delta_{a^{*,i}}$. We can rewrite (2.4.2), for fixed $a^{*,i}$ in Nash equilibrium, as:
 
 $$J^i(a^{*,i},\hat{\nu}^N)= \underset{a^{-i}\in A^{N-1}}{inf}\big[J^i(a^{*,i},a^{-i})+\omega d_1(\hat{\nu}^{N-1},\hat{\nu}^N)\big] $$
 
  $$J^i(a^{*,i},\hat{\nu}^N) \leq J^i(a^{*,i},\hat{\nu}^{N-1}) + \frac{C}{N}$$

  so $\delta_{a^{*,i}}$ is $\epsilon$-optimal also for the problem: 
  $$  \underset{\nu^i\in\mathcal{P}(A^i)}{inf} \int_A J^i(a,\hat{\nu}^{N})d\nu^i(a)  $$
  assuming %$\epsilon_N=\frac{C}{N}$ and 
  $N$ is sufficiently large. 
  
  The empirical measure $\hat{\nu}^N=\frac{1}{N}\sum_{i=1}^N \delta_{a^{*,i}}$ is also optimal since it is a linear combination and $C$ is independent of $N$ so 
  
  $$\int_A J^i(a,\hat{\nu}^N)d\hat{\nu}^N(a)\leq  \underset{\nu^i\in\mathcal{P}(A^i)}{inf} \int_A J^i(a,\hat{\nu}^N)d\nu^i(a)  + \epsilon_N $$

Letting $N$ go to infinity and the result follows

 \end{proof}

 \subsection{MFGs equilibrium}
 
 As we just mentioned a Nash equilibrium for the finite game can be extended to infinite agents and we are going to ask the opposite question in the last subsection of Section 3. Here we would like to discuss intuitively about the concept of a MFGs equilibrium and the possible similarities and differences between Nash equilibriums and MFGs equilibriums.

 \begin{d1}{MFG equilibrium}
 
  A deterministic measure-function $t\mapsto\mu_t$ is called a MFG equilibrium if $\mu_t=\mathcal{L}(\bar{X}_t)$ for each $t\in[0,T]$, for some admissible control $a^*_t$ which is optimal for representative agent's problem. With $\bar{X}_t$ the nonlinear process defined in section 2.2.4
 \end{d1}

The representative agent cannot influence $\mu_t$ (the distribution of an infinity of agents' state processes) and thus considers it as fixed when solving the optimization problem. If each agent among the infinity is identical and acts in the same way, then the law of large numbers suggests that the statistical distribution of the representative's optimally controlled state process at time $t$ must agree with $\mu_t$. 

 The aforementioned definition is very close to the concept of a Nash equilibrium apart from the fact that when the number of agents is infinite we cannot distinguish agents as we cannot distribution points in a continuous line. In a MFGs equilibrium no one has incentive to deviate since everyone plays his best response. 
 
 The biggest difference is that in a MFGs equilibrium little if anything at all can be said about the actual positions of the agents (states-actions). In the end it is just a distribution. 
 
\subsection{Games with countably infinite players versus a continuum of players  and Approximate Nash Equilibrium}

We would like to end this chapter with a small intuitive explaination about games with countably infinite players as oposed to games with a continuum of players based on the  pioneering work of Aumann \cite{Aumann1964}, Mas-Colell \cite{Mas-Colell1983}  and Aproximate Nash equilibrium by G. Carmona \cite{GCarmona2004}

As mentioned before transition to the limit in the case of interacting particles comes naturally using the idea of propagation of molecular chaos. But in the case of real humans it is not trivial how we should understand a continuum of players and how to intrprete it in a model. 

The key observation (made by G.Carmona in \cite{GCarmona2004} ) is that a game with a finite number of players is similar to a game with infinite players (countable or uncountable) if it can approximately describe the same strategic situation as the infinite. We say that a sequence of finite games approximates the strategic situation described by the given strategy in the infinite game if both the number of agents and the distribution of states and/or actions of the finite game converges to that of the infinite. We will make these statements precise in the last section of the next chapter where we are going to discuss about solutions of the finite game given we solved the infinite game. 

\begin{re1}
 In most of the economic literature discussing games with an uncountable infinite number of players \cite{Aumann1964} \cite{GCarmona2004} \cite{Mas-Colell1983}  (producing a continuum ) atomless probability spaces (or generally measure spaces) are used which is quite different from the measuretheoretic structure adopted by MFGs. 
\end{re1}

\chapter{Solution of Infinite players Mean Field Games}

\section{Revision of the infinite MFG problem}
 
 We will depart from section's 2 notation for the infinite agents problem, to make notation more compact since we are devoting the whole section to MFGs at infinity.

\begin{description}
\item[Representative agent's problem playing pure strategies]

For each fixed deterministic flow $\{\mu_t\}_{0\leq t \leq T}$ on $\mathbb{R}$ solve:
 \begin{equation}
  \underset{a_{(\cdot)}\in \mathcal{A}}{inf} J(a_{(\cdot)};\mu_{(\cdot)})=\underset{a_{(\cdot)}\in A}{inf} \int_0^T \mathbb{E}\big[ f(X_t,a_t,\mu_t) dt + g(X_T,\mu_t) \big]
 \end{equation}
subject to
\begin{equation}
\left\{
\begin{split}
 &dX^{a}_t=b(t,X_t,a_t,\mu_t)dt+\sigma(t,X_t,a_t,\mu_t)dW_t\\
 &X_0^a=\xi\in L^2\\
 \end{split}
 \right.
\end{equation}

 \item[Equilibrium of the MFG]

 Find a flow $\hat{\mu}=\{\mu_t\}_{0\leq t \leq T}$ such that 
 \begin{equation}
  \mu_t=\mathcal{L}(\hat{X}^{\hat{\mu}}_t,\hat{a}_t) \text{ for all $t\in [0,T]$}
 \end{equation}
 
 where $(\hat{X}^{\hat{\mu}}_t,\hat{a}_t)$ is the optimal pair of state and control processes solution of the Representative agent's control problem 
\end{description}

We placed the over-scripts $a,\mu$ to keep track with respect to what we optimize at each step because if we assume the control $a_t$ has some special (feedback) form as articulated in section 2 then the equilibrium of the MFG comes as a fixed point of implicit functions and can be tricky to keep track of the notations.

To make thinks easier we will work in the same set-up as in the previous case where the volatility is constant and the agents interact only through states, this keeps the presentation simpler and at the same time wealthy enough to understand the ideas better. 
 
\section{Preliminaries}
We begin our effort to solve the MFG problem by reviewing some notions about spaces of probability measures which will come in handy.

\begin{d1}{The space $\mathcal{P}_p(A)$}

 Let $A$ be a compact subset of a metric space $(\mathcal{S},d)$ (or the space itself) and $\mathcal{P}(A)$ the space of probability measures on $A$. We define $\mathcal{P}_p(A)\subset\mathcal{P}(A)$ as the space of probability measures of order $p$, $p\in\mathbb{R}^+$ with the $p$-th power of the distance to a fixed point $a_0\in A$ integrable i.e. 
 $$\big( \int_A d(a_0,a)^pd\mu(a)\big)^{\frac{1}{p}}<\infty$$
\end{d1}

In the next subsection where we discuss the maximum principle approach to the control problem we are going to work with flows in $\mathcal{P}_2(A)$ because $\mathcal{P}(A)$ is rather big to achieve our result. For $\mathcal{P}_p(A)$ we have available also the theorems from the introduction. 

We would like now to give a representation of probability measures in terms of random variables which will be our guiding intuitions in the next sections. Our aim is to find a random variable with a given law on any space. First we will work an example in a subset of $\mathbb{R}$ and then to a general metric space. 

\begin{ex1}
 Suppose we have $\big((0,1),\mathcal{B},\mathbb{P}\big)$ a fixed probability space, with $\mathcal{B}$ the Borel $\sigma$-algebra, for any distribution function $F$ on $\mathbb{R}$ let:
 $$X(t)=X_F(t)=inf\{x:F(x)>t\} \hspace{1mm} 0<t<1 $$

\begin{p1}
 For any distribution function $F$, $X_F$ is a random variable with distribution function $F$. 
\end{p1}
\begin{proof}
 
\end{proof}

\end{ex1}

\begin{t1}{Skorokhod's representation theorem}
 
Let $(\mu_n)_{n\in\mathbb{N}}$ be a sequence of probability measures on a metric space $\mathcal{S}$ such that $\mu_n\to\bar{\mu}$ on $\mathcal{S}$ when $n\to \infty$ and the support of $\bar{\mu}$ is separable. Then there exist random variables $X_n$ defined on a common probability space $\big(\Omega,\mathcal{F},\mathbb{P}\big)$ such that 
\begin{enumerate}
 \item $ \mathcal{L}(X_n)=\mu_n \hspace{1mm} \forall n\in \mathbb{N}$
 \item $X_n\to\bar{X}\hspace{1mm} \mathbb{P}$-a.s. as $n\to \infty$
\end{enumerate}

\end{t1}

The proof far exceed the purpose of this text and is omitted.

\begin{d1}
 A random variable $X$ with values in $(\mathcal{S},d)$ is said to be of order $p$, $p\in\mathbb{R}^+$ if $\mathbb{E} [d(x_0,X)^p]<\infty \hspace{1mm} \forall x_0\in E$. Moreover
 \begin{itemize}
  \item For $X,Y$ with values in $(\mathcal{S},d)$ and the Kantorowich-Rubinstein distance we have 
  $$d_1(\mathcal{L}(X),\mathcal{L}(Y))\leq \mathbb{E} [d(X,Y)] $$
  \item For $\mu,\nu$ of order 1
  $$d_1(\mu,\nu)=inf\{\mathbb{E} [d(X,Y)]|\mathcal{L}(X)=\mu,\mathcal{L}(Y)=\nu\} $$
 \end{itemize}

\end{d1}

\section{The Representative Agent's problem}

Now we are ready to solve the Representative agent's problem. The traditional way to study these problems is first to check for existence of a minimizer of (3.1.1) an optimal control process as we may call it and then  identify it using the stochastic maximum principle. Working this way one has to start with an appropriate space of continuous then using Girsanov's theorem translate the canonical process of the space into the state process and using convexity and/or compactness arguments retrieve the optimal control as a weak limit. Indeed a similar method can be used to cope with MFGs problems from a probabilistic point of view for example [citation] 

However, we will use a different approach. The strategy we present here to find the MFG equilibrium is to use a form of the Stochastic Maximum Principle (SMP) to connect existence of a minimizer for the control problem when the $\mu_t$ is fixed with a system of Forward-Backward Stochastic Differential Equations (FBSDEs). Once the flow of probability measures giving the fixed point (3.1.3) is injected in the FBSDE, then the equilibrium of the MFG comes as the solution of a FBSDE system of McKean-Vlasov type.

As usual in control problems first we have to define the Hamiltonian and the associated adjoint process.

\begin{d1}{Hamiltonian}

 Let $H:[0,T]\times\mathbb{R}\times \mathcal{P}(A)\times\mathbb{R}\times A\to \mathbb{R}$ with 
 $$H(t,x,\mu,y,a)=\rangle b(t,x,\mu,a),y \langle +f(t,x,\mu,a)$$
 be the classical Hamiltonian associated with control problem (3.1.1-3.1.2) and y the costate variable
\end{d1}

\begin{re1}
If we allow for control in volatility we have to make a series of "corrections" in our SMP approach
\begin{itemize}
 \item We have to "correct" the Hamiltonian with a risk adjustment term since the decisions on the control can affect the volatility of the state and so increase the uncertainty of the controller for future costs. 
 \item In the same spirit we have to introduce a second adjoint process to reflect this inter-temporal risk optimization.  
\end{itemize}
\end{re1}

\begin{d1}{First order adjoint process}

 We call the solution of 
 \begin{align}
  &dY_t=-\frac{\partial H}{\partial x} (t,X_t,\mu_t,Y_t,\hat{a}(t,X_t,\mu_t,Y_t))dt+Z_tdW_t\\
  &Y_T=\frac{\partial g}{\partial x}(X_T,\mu_T)
 \end{align}
 first order adjoint process associated with control problem (3.1.1-3.1.2)
 
\end{d1}

\begin{re1}
 We use $\hat{a}$ for solutions of optimal control problems while we same $*$ for Nash equilibriums, in a later section where we are going to discuss their connection we will revise our notation. 
\end{re1}

\subsection{Assumptions}

Unfortunately we cannot continue the discussion at full generality and we have to impose some assumptions to achieve our existence theorem. Usually to minimize the Hamiltonian we require some convexity and in this particular case we will demand an affine structure in $b(t,x,\mu,a)$ along  with convexity in $f(t,x,\mu,a)$

We list the complete set of our assumptions here and refer to them whenever needed providing also some motivation. We denote $S_0-S_3$ the assumptions required for retrieving the stochastic maximum principle and $F_4-F_6$ for the fixed point problem, even thought all of them have to be fulfilled for the MFG equilibrium to exists. 

$\vspace{1mm}$

\textbf{Assumptions}
%$\vspace{1mm}$
\begin{enumerate}
\item[$S0$.] $A\subset \mathbb{R}$ compact, convex and the flow of probability measures $\mu_t$ is deterministic.

 \item[$S1$.]
 $$ b(t,x,\mu,a)=b_0(t,\mu)+b_1(t)x+b_2(t)a $$
 where $b_2$ is measurable and bounded and $b_0,b_1$ measurable and bounded on bounded subsets of $[0,T]\times\mathcal{P}_2(A)$ and $\mathbb{R}$ respectively.
 
 \item[$S2$.] $f(t,\cdot,\mu,\cdot)$ is $C^{1,1}$, bounded with bounded derivatives in $x,a$ and satisfies the convexity assumption
 $$f(t,x',\mu,a')-f(t,x,\mu,a)-(x'-x)\frac{\partial f}{\partial x}-(a'-a)\frac{\partial f}{\partial a}\leq \lambda |a'-a|^2\hspace{1mm} \lambda\in\mathbb{R}^+ $$
 
 \item[$S3$.] $g$ is bounded, for any $\mu\in\mathcal{P}_2(A)$ the function $x\to g(x,\mu)$ is $C^1$ and convex. 
 
 \item[$F4$] $b_0,b_1,b_2$ are bounded by $c_L$. Moreover for any $\mu,\mu'\in\mathcal{P}_2(\mathbb{R})$ we have $|b_0(t,\mu')-b_0(t,\mu)|\leq c_L d_1(\mu,\mu')$
 
 \item[$F5$] $|g_x(x,\mu)|\leq c_B$ and $|f_x(t,x,\mu,a)|\leq c_B$ for all $t\in[0,T], x\in\mathbb{R},\mu\in\mathcal{P}_2(\mathbb{R}),a\in\mathbb{R}$

\end{enumerate}

\subsection{Stochastic Maximum Principle}

We begin with a version of the SMP from Pham to achieve existence and motivate our strategy. 

\begin{t1}{(Pham)}

 Let $\hat{a}\in A$ and $\hat{X}$ the associated state process.  Assume
 \begin{enumerate}
 \item $g$ is bounded $C^1$ and convex
  \item There exist a solution $(Y_t,Z_t)_{0\leq t\leq T}$ of the BSDE: 
  \begin{equation}
 \left\{
 \begin{split}
  &dY_t=-\frac{\partial H}{\partial x} (t,X_t,\mu_t,Y_t,\hat{a}(t,X_t,\mu_t,Y_t))dt+Z_tdW_t\\
  &Y_T=\frac{\partial g}{\partial x}(X(T),\mu_T)
 \end{split}
\right.
\end{equation}
such that 
$$\mathcal{H}(t,\hat{X}_t,\mu_t,\hat{a}_t,\hat{Y}_t,\hat{Z}_t)=\underset{a_t\in A}{min}H(t,\hat{X}_t,\mu_t,a,\hat{Y}_t,\hat{Z}_t) \hspace{1mm} 0\leq t\leq T \hspace{1mm} a.s.$$
\item $(x,a)\to H(t,x,\mu,a,\hat{Y}_t,\hat{Z}_t)$ is a convex function for all $t\in[0,T]$. 
  \end{enumerate}
Then $\hat{a}$ is an optimal control.
\end{t1}

\begin{re1}
 Since $\mu_t$ is fixed, it does not affect our calculations, which are pretty standard for stochastic control problems and we can drop it from our notation without any harm. 
\end{re1}

\begin{proof}
 For any $a_t\in A$ we need to calculate:
 $$J(\hat{a}_t)-J(a_t)=\mathbb{E}\big[ \int_0^T f(t,\hat{X}_t,\hat{a}_t)-f(t,\hat{X}_t,a_t)dt+g(\hat{X}_T)-g(X_T)\big] $$
 
 So we need to get estimates for $ f(t,\hat{X}_t,\hat{a}_t)-f(t,X_t,a_t)$ and $g(\hat{X}_T)-g(X_T)$ 
 
 \begin{description}
  \item For $g(\hat{X}_T)-g(X_T)$ we use the assumption 1, along with the BSDE and Ito's rule
  \item For $ f(t,\hat{X}_t,\hat{a}_t)-f(t,X_t,a_t)$ we change $f$ for $H$ and then use assumption 3. 
  \item In the end we combine the estimates with the definition of  $\mathcal{H}$ to prove that 
  $$J(\hat{a}_t)-J(a_t)<0 \text{ for $0\leq t\leq T$}$$ and taking the inf the desired relationship comes. 
 \end{description}

\end{proof}

We now state and prove a variation of the classical Stochastic Maximum Principle in the spirit of Pham's theorem, tailored for our needs. Apart from the FBSDE system which we are going to use in the next step, we can find a way to compare control as the inequality (3.3.5) show in the theorem which is also going to be helpful in the next step.

\begin{t1}
 Assume $(S_0-S_1)$ in addition, if the map $t\to \mu_t\in\mathcal{P}_2(\mathbb{R})$ is measurable and bounded and if the FBSDE system:
 \begin{equation}
 \left\{
 \begin{split}
  &dX_t=b(t,X_t,\mu_t,\hat{a}(t,X_t,\mu_t,Y_t))dt+\sigma dW_t\\
  &dY_t=-\frac{\partial H}{\partial x} (t,X_t,\mu_t,Y_t,\hat{a}(t,X_t,\mu_t,Y_t))dt+Z_tdW_t\\
  &X_0=x_0\\
  &Y_T=\frac{\partial g}{\partial x}(X(T),\mu_T)
 \end{split}
\right.
\end{equation}
has a solution $(X_t,Y_t,Z_t)_{0\leq t\leq T}$ such that 
$$\mathbb{E}\big[\underset{0\leq t\leq T}{sup}(|X_t|^2+|Y_t|^2)+\int_0^T|Z_t|^2dt \big]<\infty $$ 
Then for any admissible $(a_t)_{0\leq t\leq T}$ the variational inequality:

\begin{equation}
 J(\hat{a}_t;\mu)+\lambda \mathbb{E} \int_0^T |a_t-\hat{a}_t|^2dt \leq J(a_t;\mu)
\end{equation}

 where $\hat{a}_t=(\hat{a}_t)_{0\leq t\leq T}$ is the minimizer of the Hamiltonian, holds
\end{t1}

\begin{proof}
 As before we need to calculate $J(\hat{a}_t)-J(a_t)$ with $\hat{a}_t$ the minimizer of the Hamiltonian and $a_t$ admissible. As before we drop $\mu$ to lighten notation since it is fixed and doesn't affect calculations.
 $$J(\hat{a}_t)-J(a_t)=\mathbb{E}\big[ \int_0^T f(t,\hat{X}_t,\hat{a}_t)-f(t,\hat{X}_t,a_t)dt+g(\hat{X}_T)-g(X_T)\big] $$
 So we need to get estimates for $ f(t,\hat{X}_t,\hat{a}_t)-f(t,\hat{X}_t,a_t)$ and $g(\hat{X}_T)-g(X_T)$ 
 
 \begin{description}
  \item[For $g(\hat{X}_T)-g(X_T)$] we use the $(S3)$, and so we get
  $$g(X_T)\geq g(\hat{X}_T)+g_x(\hat{X}_T)(X_T-\hat{X}_T) $$
  $$g(\hat{X}_T)-g(X_T)\leq Y_T(\hat{X}_T)(X_T-\hat{X}_T) $$
  using Ito's rule we end up with
 $$
   \mathbb{E}\big[g(\hat{X}_T)-g(X_T)\big]\leq \mathbb{E}\big[\int_0^T-(X_t-\hat{X}_t)H_{\hat{x}}dt+\int_0^TY_t(b(t,\hat{x},\hat{a})-b(t,x,a))dt\big]
  $$
and imposing $(S_1)$
\begin{equation}
 \mathbb{E}\bigg[g(\hat{X}_T)-g(X_T)\bigg]\leq \mathbb{E}\bigg[\int_0^T-(X_T-\hat{X}_T)H_{\hat{x}}dt+\int_0^T\hat{Y}_t\bigg(b_1(t)(\hat{X}_t-X_t)+b_2(t)(\hat{a}_t-a_t)\bigg)dt\bigg]
\end{equation}

  \item[For $ f(t,\hat{X}_t,\hat{a}_t)-f(t,X_t,a_t)$]we use $(S2)$. 
  $$\mathbb{E}\bigg[\int_0^Tf(t,\hat{X}_t,\hat{a}_t)-f(t,X_t,a_t)dt\bigg]\leq\mathbb{E}\bigg[\int_0^T(\hat{X}_t-X_t)f_x-(\hat{a}_t-a_t)f_a+ \lambda |\hat{a}_t-a_t|^2 dt\bigg] $$
  using the definition of the Hamiltonian
  
  \begin{align*}
   \mathbb{E}\bigg[\int_0^Tf(t,\hat{X}_t,\hat{a}_t)-f(t,X_t,a_t)dt\bigg]\leq\mathbb{E}\bigg[&\int_0^T(\hat{X}_t-X_t)(H_{\hat{x}}-b_{\hat{x}}Y_t)-(\hat{a}_t-a_t)(H_{\hat{a}}-b_{\hat{a}}Y_t)\\ 
   &+ \lambda |\hat{a}_t-a_t|^2 dt\bigg]  
   \end{align*}
   
   \begin{align*}
   =\mathbb{E}\bigg[\int_0^T(\hat{X}_t-X_t)H_{\hat{x}}-(\hat{a}_t-a_t)H_{\hat{a}}+Y_t\big((\hat{a}_t-a_t)b_{\hat{a}}-(\hat{X}_t-X_t)b_{\hat{x}}\big)+\lambda |\hat{a}_t-a_t|^2 dt\bigg] 
  \end{align*}
  
   \begin{equation}
    =\mathbb{E}\bigg[\int_0^T(\hat{X}_t-X_t)H_{\hat{x}}-(\hat{a}_t-a_t)H_{\hat{a}}+Y\big((\hat{a}_t-a_t)b_2(t)-(\hat{X}_t-X_t)b_1(t)\big)+\lambda |\hat{a}_t-a_t|^2 dt\bigg] 
  \end{equation}
  
 \end{description}

 We sum (3.3.6) and (3.3.7) to get:
 \begin{equation}
  \begin{split}
   &J(\hat{a}_t)-J(a_t)\leq \mathbb{E}\bigg[\int_0^T -(\hat{a}_t-a_t)H_{\hat{a}}+\lambda |\hat{a}_t-a_t|^2 dt\bigg]\\
  \end{split}
 \end{equation}
 
All that is left is to prove existence and uniqueness of $\hat{a}_t$. We take care of that with the following lemma. 

\begin{l1}{Minimization of the Hamiltonian}
  
Assume $(S_0-S_2)$ then for each $(t,x,\mu,y)$ is the appropriate domain there exists a unique minimizer of the Hamiltonian $H$  
\end{l1}

And so $H_{\hat{a}}=0$ and the proof is complete.

\end{proof}

 \section{The fixed point problem}

 The second step for the MFG equilibrium is now to find a family of probability distributions $(\mu_t)_{0\leq t\leq T}$ such that the process $\{\hat{X}_t\}_{0\leq t\leq T}$ solving (3.2.1) admits $(\mu_t)_{0\leq t\leq T}$ as flow of marginal distributions i.e. 
 $$\mu_t=\mathcal{L}(\hat{X}_t,) \text{ for all $t\in [0,T]$} $$
  so we rewrite (3.2.1) as a McKean-Vlasov FBSDE system:
 \begin{equation}
 \left\{
 \begin{split}
  &d\hat{X}_t=b(t,\hat{X}_t,\mathcal{L}(\hat{X}_t),\hat{a}(t,\hat{X}_t,\mathcal{L}(\hat{X}_t),Y_t))dt+\sigma dW_t\\
  &dY_t=-\frac{\partial H}{\partial x} (t,\hat{X}_t,\mathcal{L}(\hat{X}_t),Y_t,\hat{a}(t,\hat{X}_t,\mathcal{L}(\hat{X}_t),Y_t))dt+Z_tdW_t\\
  &\hat{X}_0=x_0\\
  &Y_T=\frac{\partial g}{\partial x}(X(T),\mu_T)
 \end{split}
\right.
\end{equation}
 
 and we have the following theorem for solving MKV-FBSDEs. 
 
 \begin{t1}
  Under Assumptions ($S0-F4$) the FBSDE system (3.4.1) has a solution $(X_t,Y_t,Z_t)_{0\leq t\leq T}$ . Moreover for any solution, there exists a function $u:[0,T]\times \mathbb{R}^d$ such that $Y_t=u(t,X_t) \hspace{1mm} a.s. \forall t\in[0,T]$ and satisfies the growth and Lipschitz properties:
  \begin{enumerate}
   \item $ |u(t,x)|\leq c(1+|x|)\hspace{1mm}  \forall t\in[0,T] $ 
   \item $|u(t,x)-u(t,x')|\leq c|x-x'| \hspace{1mm}  \forall t\in[0,T], \hspace{1mm}  \forall x,x'\in \mathbb{R}^d $ 
  \end{enumerate}

 \end{t1}

 The theorem itself is difficult to prove and the strategy behind it is far from trivial. Let's start articulating the steps we need to follow to solve this problem. 
 
 \begin{enumerate}
  \item Given a flow of probability measures $\mu\in\mathcal{P}_2(C)$ with $\mu=(\mu_t)_{0\leq t\leq T}\mathcal{L}(\hat{X}_t)$ and $C$ the space of real continuous functions, we prove that the FBSDE system is uniquely solvable. 
  \item We set $\Phi:\mu\mapsto\mathcal{L}(\hat{X}^{x_0;\mu})$ the map that associates each $\mu$ with $\mathcal{L}(\hat{X}^{x_0;\mu})$ the probability distribution of state process, solution of the previous step.
  \item MFG equilibrium comes as a fixed point of $\Phi$ 
 \end{enumerate}

  In the first step we fix $\mu$ and for each particular fixed $\mu$ we solve the FBSDE to achieve this we use a similar approach as in the previous subsection where we rely on existing theory of FBSDE system for existence and uniqueness. In the second step we associate each solution of the FBSDE with a flow of probability measures and in the third step we use the compactness of $\mathcal{P}(C)$ to apply the Schauder's fixed point theorem.  

The complete proof of Theorem 3.4 is long and cumbersome so we will prove only the most important points that are going to help us gain a better understanding of the subject. 

For the first step we have the following lemma. 

\begin{l1}
 Given $\mu\in\mathcal{P}_2(C)$ with marginal distributions $(\mu_t)_{0\leq t\leq T}=\mathcal{L}(\hat{X}_t)$ and $C$ the space of real continuous functions the FBSDE:
 \begin{equation}
 \left\{
 \begin{split}
  &d\hat{X}_t=b(t,\hat{X}_t,\mathcal{L}(\hat{X}_t),\hat{a}(t,\hat{X}_t,\mathcal{L}(\hat{X}_t),Y_t))dt+\sigma dW_t\\
  &dY_t=-\frac{\partial H}{\partial x} (t,\hat{X}_t,\mathcal{L}(\hat{X}_t),Y_t,\hat{a}(t,\hat{X}_t,\mathcal{L}(\hat{X}_t),Y_t))dt+Z_tdW_t\\
  &\hat{X}_0=x_0\\
  &Y_T=\frac{\partial g}{\partial x}(X(T),\mu_T)
 \end{split}
\right.
\end{equation}
has a unique solution 
\end{l1}

We will not attempt a complete proof here but we would rather sketch some arguments. First we notice that the assumptions $(F-F)$ with $\mu$ fixed and bounded and properties of the driver o the BSDE gives as existence according to the fairly straight Forward 4-step scheme we developed in the appendix. Giving just a short reminder here, we suppose a deterministic function $\theta$ such that $Y_t=\theta(t,X_t)$ and applying Ito's rule we end up with a quasiliniar non-degenerate parabolic PDE.  However this is only a local result for small time, the idea of the extension for arbitrary time as proposed by Delarue in \cite{Delarue2002} is the following: we assume it holds in an interval of the form $[T-\delta,T]$ with $\delta$ sufficiently small including also $t_0$ of $x_0$. There $\theta(T-\delta,\cdot)$ plays the role of the terminal data in the BDSE namely $g_x$. 

$\vspace{1mm}$

For step 2 we give the following definition

\begin{d1}
 For any continuous flow of probability measures $\mu=(\mu_t)_{0\leq t\leq T}\in \mathcal{P}_2(C)$ and $(\hat{X}_t^{x_0;\mu},\hat{Y}_t^{x_0;\mu},\hat{Z}_t^{x_0;\mu})$ solution of FBSDE (3.4.2) we define the map $\Phi:\mu\mapsto\mathcal{L}(X^{x_0;\mu})$ and we call MFG equilibrium or solution, any fixed point of $\Phi$.
\end{d1}

We take care of the fixed point in the next lemma

\begin{l1}
There exists a closed convex subset $E$ of $\mathcal{P}_2(C)$ which is stable for $\Phi$ with a relatively compact range, $\Phi$ is continuous on $E$. $\Phi$ has a fixed point.  
\end{l1}

\begin{proof}
In the proof we will make use of
\begin{t2}{Schauder's fixed point theorem}

 Let $E$ be a nonempty, compact, convex subset of a Banach space and $\Phi$ a continuous (compact) map from $E$ to itself then $\Phi$ has a fixed point. 
\end{t2}

We need to identify a compact convex subset of $\mathcal{P}_2(C)$ and prove that $\Phi$ defined as earlier maps this subset to itself and is continuous. 

To do this we start by looking for  bounds of our solutions. For assumptions $(F4-F5)$ we get that 
$$|H_x|=|b_1(t)y+f_x|\leq c_L|y|+c_B$$
and so if we write $Y_t=g_x(X_T,\mu_T)-\int_t^T-H_x dt-\int_t^TZ_tdW_t$ and under our assumption and a comparison principle for SDEs [citation] it is straightforward that:
$$\text{for any $\mu\in\mathcal{P}_2(C)$ and } t\in [0,T]\hspace{3mm} |Y_t^{x_0;\mu}|\leq c \hspace{2mm} \text{ a.s.} $$ 
where $c$ depends upon $c_B,c_L $ and $T$. 

Remembering Lemma 3.3.1 with our assumptions yields 

\begin{align}
 &|\hat{a}(t,x,\mu,y)|\leq \lambda^{-1}(c_L+c_L|y|)  \notag \\
 &|\hat{a}(t,X_t^{x_0;\mu},\mu,Y_t^{x_0;\mu})|\leq \frac{c_L}{\lambda}(1+c)=c'
\end{align}

By Theorem 5.4 in \cite{CD2013} we have 
$$\mathbb{E}(\underset{0\leq t\leq T}{sup} |X_t^{x_0;\mu}|^2)\leq c' (1+\mathbb{E}(|x_0|^2))=K$$

So we consider the set:

$$E:=\{\mu\in\mathcal{P}_4(C): \underset{0\leq t\leq T}{sup}\int_\mathbb{R}|x|^4d\mu_t(x)\leq K  \}$$

$E$ is convex and closed in the $d_1$. 

Now we have to show also that it is relatively compact, and to do so we use a tightness argument with Prohorov's theorem for the family of processes $((X_t^{x_0;\mu})_{0\leq t \leq T})_{\mu\in E}$ and the corresponding laws. 

For continuity we need to show that: 

\begin{description}
 \item $\Phi$ is continuous in $\mu' \in E$ $\iff$ For all $\epsilon>0$ there exists $\delta>0$ such that if $\mu\in E$ and $d_1(\mu,\mu')<\delta$ then $d_1(\Phi(\mu),\Phi(\mu'))<\epsilon$
\end{description}

 From Definition 3.2 in the Preliminaries subsection of this section we write:
 $$d_1(\Phi(\mu),\Phi(\mu'))=d_1(\mathcal{L}(X_t^{x_0;\mu}),\mathcal{L}(X_t^{x_0;\mu'}))\leq \mathbb{E}\underset{0\leq t\leq T}{sup}|X_t^{x_0;\mu}-X_t^{x_0;\mu'}| $$
 
 Now to get an estimate for $|X_t^{x_0;\mu}-X_t^{x_0;\mu'}|$ we need to use (3.3.5) along with the state process under the "environment" $\mu'$.  
\end{proof}

As Carmona and Delarue showed in their original paper \cite{CD2013} we can relax assumption $(F6)$ by approximating cost functions $f,g$ by sequences of functions $f_n,g_n$ that satisfy $(F5)$ uniformly.
 
 \section{Analytic approach, Connection of SMP with dynamic programming}
 
 Here we are going to describe the so called analytical method. Since the initial appearance of the MFGs in the mathematical literature it has served as the primary solution method and has been intensively studied. It has its roots in the Dynamic Programming Principle as introduced by Bellman and classical analytical mechanics.
 
 In a nutshell the method uses a value function which solves a special kind of PDE called Hamilton-Jacobi-Bellman and under appropriate assumptions can produce an optimal control in feedback form which coupled with the state dynamics (stochastic or not) can give us solution to optimal control problems. 
 
 However in our case we have many agents that optimize, each one with his optimal control problem which are coupled since they interact through their states and/or controls. So when each agent takes decisions have to take into account the empirical distribution of the states and/or controls of the other players (in our case only the states for simplicity). This interaction indicates that when we solve the MFG problem i.e. search for a distribution of states that no one has intention to deviate, in addition to solving the optimal control for each agent, we have also to describe the evolution of the state distribution. This idea was first introduced in subsection 2.3.1 - 2.3.5 where we ended with a Fokker-Plank equation for the evolution of particle distribution in the case of particles and a states in our MFGs case. Here we are going to introduce the HJB equation and couple it with the FP equation to derive a MFG equilibrium as defined earlier. 
 
 \subsection{Hamilton Jacobi Bellman}
 
 We use the same definition of the Hamiltonian as before in addition we let $\mathcal{H}$ the Legendre transform of H i.e. $\mathcal{H}=\underset{a\in A_{adm}}{inf}H$
 
 We define the value function of representative's agent problem u as:
 $$u(x,t)=\underset{a\in A_{adm}}{inf}J(a_{(\cdot)};\mu_{(\cdot)}) $$ 
 
As pointed in the appendix if $u$ is sufficiently regular it solves the HJB equation:

\begin{equation}
\left\{
\begin{split}
u_t+\mathcal{H}(x,-u_x,-u_{xx}, \mu_t)=0\\
u(x,T)=g(x)
\end{split}
\right.
\end{equation}

\begin{re1}
The usual way to relax the assumptions about (3.5.1) is to look for viscosity solutions but this is a concept that we are not going to discuss.  
\end{re1}

Once we have a solution of the HJB, using also the SDE we can compute $\hat{a}_t  \in A_{adm}$ the optimal control function were the specific form depends also on the modelling we use namely open or close loop controls. To make everything rigorous we need a verification theorem. 

As explained in the previous sections the in order to solve the Representative's agent problem we fix the flow of probability measures, $\mu_t$. We use the fixed point condition as earlier: 
$$\mu_t=\mathcal{L}(\hat{X}_t) \text{ for all $t\in [0,T]$} $$
and this implies coupling with the Fokker-Plank equation.

 \subsection{Fokker Plank}
 
\begin{d1}{Infinitesimal Generator }
\begin{itemize}
 \item For a general Markov process $\{X_t\}_{0\leq t\leq T}$  starting from $x_0$ we define the Infinitesimal generator of the process as: 
 $$\lim_{t\downarrow t_0} \frac{\mathbb{E}_{x_0}[f(X_t,t)]-f(x_0,t_0)}{t-t_0}$$
 \item For a process that satisfies our state SDE we have the following definition
 
$$A[f(t,x)]=\frac{\partial}{\partial t}f(t,x)+b(t,x,\mu)\frac{\partial}{\partial x}f(t,x)+\frac{1}{2}\sigma^2\frac{\partial^2}{\partial x^2}f(t,x) $$
\end{itemize}

\end{d1} 
 
\begin{d1}{Adjoint operator}

 Let $A$ be an operator we define the adjoint operator of $A$ as:
 $$\int_\mathbb{R}\phi(x)A[f(x)]dx=\int_\mathbb{R}f(x)A^*[\phi(x)]dx \hspace{3mm}  \forall \phi\in C_0^\infty (\mathbb{R})$$
 
\end{d1}

 As already seen in section 2 but with alternative notation now the evolution of the population's distribution, given an initial distribution $\mu_0$ is given by:
 
 \begin{equation}
  \left\{
  \begin{split}
  & A^*[\mu_t]=0\\
  & \mu_0=\mathcal{L}(X_0)
  \end{split}
\right.
 \end{equation}

 While in our case (3.5.2) becomes
 
 \begin{equation}
  \left\{
  \begin{split}
  & \partial_t\mu_t+b(t,x,\mu_t)\partial_x\mu_t-\frac{1}{2}\sigma \partial_{xx} \mu_t=0\\
   &\mu_0=\mathcal{L}(X_0)
  \end{split}
\right.
 \end{equation}
 
 Combining (3.5.1) with (3.5.3) we have a system of coupled PDEs the solution of which provide us with a MFG equilibrium distribution as in Definition 3.5. 
 
 \begin{equation}
\left\{
\begin{split}
&u_t+\mathcal{H}(x,-u_x,-u_{xx}, \mu_t)=0\\
& \partial_t\mu_t+b(t,x,\mu_t)\partial_x\mu_t-\frac{1}{2}\sigma \partial_{xx} \mu_t=0\\
&\mu_0=\mathcal{L}(X_0)\\
&u(x,T)=g(x) 
\end{split}
\right.
\end{equation}
 
We can see an analogy with the MKV-FBSDE system (3.4.1) since here also we have the HJB equation backward in time and the FP forward in time. But here we are dealing with infinite dimensions problem while the MKV-FBSDE problem is in finite, which also the big advantage of the probabilistic method, apart from the interpretation.

 \section{From infinite game to finite}

 We have already given, in the introduction, some motivation about our strategy to pass to the limit of infinite players. Now that we have study the infinite game enough, we would like to ask the question, \textbf{"What can we say from the situation at infinity (infinite game) about the $N$-player game?"} or in other words can we reconstruct the finite game from the infinite?
 
 First let us recall the finite game from section 2.

For $i=1,...,N$
 \begin{equation}
  \underset{a^i_{(\cdot)}\in A^i_{adm}}{inf} J^i(a_{(\cdot)};\bar{\mu}^N_{(\cdot)})=\underset{a^i_{(\cdot)}\in A^i_{adm}}{inf} \int_0^T \mathbb{E}\big[ f(X^i_t,a^i_t,\bar{\mu}^N_t) dt + g(X^i_T,\bar{\mu}^N_t) \big]
 \end{equation}
subject to
\begin{equation}
\left\{
\begin{split}
 &dX^i_t=b(t,X^i_t,a^i_t,\bar{\mu}^N_t)dt+dW^i_t\\
 &X^i_0=\xi^i\in L^2\\
 \end{split}
 \right.
\end{equation}
where
$$
 \bar{\mu}^N_t=\frac{1}{N}\sum_{i=1}^N\mathit{1}_{X_i<x}
 $$

Now suppose that we have solved the infinite game () using the SMP then we have a MFG equilibrium $\mathbf{\mu}_t $ and a value function $\theta(t,X_t)=Y_t$ from theorem () as pointed also in the appendix for FBSDE systems. Then we can define players control strategies $\hat{a}(t,X_t,\mu_t,\theta(t,X_t)$ for the infinite game in feedback form. 

We would like to set each player's strategy in the finite game as
\begin{equation}
 a^i_t=\hat{a}(t,X^i_t,\mu_t,\theta(t,X^i_t))
\end{equation}

and prove that this collection of strategies $a_t=(a^1_t,...,a_t^N)$ is indeed a Nash equilibrium for the finite game. As we have already mentioned $\hat{a}(t,X^i_t,\mu_t,\theta(t,X^i_t))$ is a feedback function for control $a^i_t$ and we are will restrict ourselves into closed loop Nash equilibriums and specifically Markovian in the spirit of Definition 2.7. As pointed in the introduction the most natural way to approach this problem is through the idea of an Approximate Nash Equilibrium (ANE). We combine the notions of Markovian and approximate Nash equilibrium in the following definition. 

\begin{d1}{Markovian Approximate Nash Equilibrium}

Let 
$$ \phi^i:[0,T]\times C([0,T];\mathbb{R}) $$
be a measurable feedback function and $\{X_t^{*,i}\}_{0\leq t\leq T}$ the solution of the state SDE (3.6.2) for $i=1,...N$ if we use the admissible controls:
$$a_t^{i,*}=\phi^i(t,X_t^i) \hspace{3mm} 0\leq t\leq T $$

Given $\epsilon>0$ the action profile $a^*_t=(a^{*,1}_t,...,a_t^{*,N})$ is an Markovian $\epsilon$-approximate Nash equilibrium if, whenever a player $i$ uses a different strategy $a_t^i=\psi(t,X_t^i)$ while the rest continue to use $b_t^j=\phi^j(t,X_t^j)\hspace{1mm}\forall j\neq i$ but with $\{X_t^i\}_{0\leq t\leq T}$ the solution of the state SDE (3.6.2) when we use the actions $a_t=(a^i_t,b^{-i}_t)$. Then 

$$J^i(a^*_t)\leq J^i(a_t)+\epsilon$$
 for each $i\in\{1,...,N\}$
\end{d1}

\begin{t1}{Existence of MANE}
 
 Under assumptions $(S0-F6)$ there exists a sequence $(\epsilon_N)_{N\geq 1}$ with $\epsilon_N\to 0$ as $N\to\infty$ such that the strategy profile $a^*_t=(a^{*,1}_t,...,a_t^{*,N})$ with $a^{*,i}_t$ defined in (3.6.3) is a Markovian $\epsilon_N$-approximate Nash equilibrium for the finite game. 
\end{t1}

The proof of this theorem is rather long and we will omit it but can be found in [citation]. Instead some remarks are on the way to elaborate more on this interesting result. 

\begin{re2}
 \begin{itemize}
  \item Theorem holds true also for open loop controls and for closed loop controls with light modification. 
  \item .....
  
  \end{itemize}

\end{re2}

\chapter{The Aiyagari Model }

 In this section we would like to present and "solve" an actual macroeconomic model \cite{Aiyagari1994}, as an application to the theory developed so far in this text. As already discussed in the introduction macroeconomic models are typical examples of MFGs since they incorporate a large number of symmetric agents usually separated in sectors who share the same incentives. 
 
 This chapter is organised as follows. At first we will make a brief introduction to the class of macro-models with a large number of agents subject to idiosyncratic shocks and present shortly some stylized facts to motivate the key features of the model that follow. Then we will present the model itself and derive existence according to chapter 3. 
 
 \section{Introduction to models with a large number of agents subject to idiosyncratic shocks }
 
 We are going to discuss about an economy with two sectors commonly refereed to as households and firms. Households provide factors of productions to firms in order to produce and gain back income as compensation. In our particular case we will assume we are dealing only with labour ($l_t$) and capital ($k_t$). We have $N$ workers working in $N$ perfectly competitive, identical firms that we will not distinguish and consider them as one representative firm (This is a common practice in macroeconomic literature). The representative firm produces only one product ($Y_t$) consumed by households and pay $w_t$ as wage and $r_t$ as compensation for capital. Since all workers are identical they get paid by the same amount, same goes for capital. We will describe the position (state) of each agent with a vector $X^i_t=(k^i_t,l_t^i)$ and use an empirical distribution for the states
 $$\mu_t^N=\frac{1}{N} \sum_{i=1}^N \delta_{X_t^i\leq x} $$
 
 The mean capital and labour are defined as:
 
 \begin{equation*}
  \left\{
  \begin{split}
  K^N_t=\frac{1}{N}\sum_{i=1}^Nk_t^i=\int kd\mu_t^N(k,l) \\
  L^N_t=\frac{1}{N}\sum_{i=1}^Nl_t^i=\int ld\mu_t^N(k,l)
  \end{split}
  \right.
 \end{equation*}

 Households can control their consumption ($c_t^i$) and as a general rule aim to maximize their discounted utility i.e. each unit of product they consume offers them a certain satisfaction and they have specific preferences regarding the time horizon of the satisfaction. We represent their preferences with a utility function $U(c^i_t)$ satisfying certain assumptions which we are going to specify later. They face a budget constraint: 
 
 $$ w_tl^i_t+r_tk^i_t=c^i_t+\frac{dk^i_t}{dt} $$
 
 where on the l.h.s we have the income of agent $i$ and on r.h.s. we have the expenditure namely consumption and rate of capital accumulation or decrease. 
 
 Firms on the other hand aim to maximize their profit while they control the mean capital ($K_t$) and  labour ($L_t$) they enter in their production. This is a reasonable assumption since all of them are identical and appear as one representative firm.  
 
 Some stylized facts provided by \cite{Aiyagari1994} will guide us to specify our model further. 
 
\begin{description}
 \item[Stylized Facts] 
\begin{enumerate}
 \item Barsky, Mankiw, Zeldes 1986  and Deaton 1991 
 
 Individual consumptions are much more variable than aggregate. This indicates that heterogeneity may be important due to incomplete markets.
 
 \item C. Carroll 1991 
 
 Individual wealth holdings are highly volatile which is hard to explain in absence of temporary idiosyncratic shocks. 
 
 \item Mankiw Zeldes 1991 
 
 Considerable diversity in portfolio compositions for households with different wealth levels. 
 
 \item Avery, Elliehusen Kennickell 1988
 
 The top end of the wealth scale owns stocks while low end owns liquid assets, which is hard to explain when the markets are frictionless. 
 
\end{enumerate}
 \begin{re1}
  The situation described by the above empirical research represents U.S. economy during the late 20th century and currently  might seem obsolete but we are presenting it for completeness and educational reasons. 
 \end{re1}

 \item[Model's Key features ]
 
 \begin{enumerate}
  \item Endogenous heterogeneity
  \item Aggregation through mean field interactions
  \item Finite horizon
  \item Borrowing constraint
  \item General equilibrium (endogenously determined interest rate)
 \end{enumerate}

 \end{description}

 In a nutshell we are dealing with an income fluctuation problem: households face uncertain earnings.
 
\section{Model formulation}

In order to solve this income fluctuation problem households take decisions on consumption or alternatively assets accumulation or decrease in order to maximize expected value of the discounted utility of consumption.

Let's introduce some notation to specify our model.
\begin{d1}

 The following will be useful
 \begin{enumerate}
  \item $c^i_t$: agent's $i$ consumption  
  \item $k^i_t$: agent's $i$ capital
  \item $l^i_t$: agent's $i$ labor
  \item $w^i_t$: agent's $i$ wage
  \item $U^i(c^i_t)$: agent's $i$ utility function 
  $$U:A\to \mathbb{R} $$
  \item $\beta$: utility discount factor
  \item $\lambda=\frac{(1-\beta)}{\beta}>0 $: time preference
  \item $r_t$: the return on capital
 \end{enumerate}

\end{d1}

\textbf{Assumptions}

\begin{enumerate}
 \item[(A1).] 
We we model labour endowment shocks by  allowing  $l_t$  to be a solution of an SDE: 
  
  $$ dl^i_t=b(l^i_t)dt+\sigma(l^i_t)dW_t$$

  We can specify as: 
  \begin{itemize}
   \item[(B).] Geometric Brownian Motion
   \item[(OU).] Orstein - Uhlenbeck process. 
  \end{itemize}

Also, $(l_t^i)_{i=1}^N$ is is iid with bounded support given by $[l_{min},l_{max}]$ with $l_{min}>0$

\item[(A2).]  Utility function is of Constant Relative Risk Aversion type (CRRA) given by:
\begin{equation}
 U(c)=\frac{c^{1-\gamma}-1}{1-\gamma}
\end{equation}
for $\gamma>0$ with $U(c)=ln(c)$ if $\gamma=1$

Also, $\beta$, the utility discount factor is constant and time invariant

\item[(A3).]  The production function is Cobb-Douglas type. 
\begin{equation}
 F(K_t,L_t)=K_t^aL_t^{1-a}
\end{equation}

\item[(A4).]  We normalize $L_t$ to 1 i.e.
$$\mathbb{E}[l^i_t]=1 $$ for all $t>0$ and $i\in\{1,2,..,N \}$

\end{enumerate}

\textbf{Individual agent's problem}

The households are interested in maximizing:
\begin{equation}
 \underset{c^i_t\in A_{adm}}{max}\mathbb{E}_0\big[\int_0^Te^{-\beta t}U(c^i_t)dt+ \tilde{U}(c^i_T)\big]
\end{equation}

subject to a budget constraint 

\begin{align}
 &dk_t=[-c_t+w_tl_t+r_tk_t]dt\\
 &c_t\geq 0\\
 &k_t\geq -b \text{ a.s.}
\end{align}

$c_t\geq 0$ seems logical since there is no meaning in negative consumption, 

$k_t\geq -b$ is a limit on borrowing since in real world there is no such thing as an infinite credit line.

The limit on borrowing is no enough to rule Ponzi strategies, we need a borrowing constraint:

\begin{equation}
 \left\{
 \begin{split}
  &k_t\geq -\phi\\
  &\phi=
  \begin{cases}
   min\{b, \frac{wl_{min}}{r} \} \text{ if } r>0\\
   b \text{ if } r \leq 0
  \end{cases}
 \end{split}
\right.
\end{equation}

To incorporate the borrowing constraint in the model  we need to define:

$$\tilde{k}_t=k_t+\phi $$ 

so the budget constraint becomes:

\begin{align}
 &d\tilde{k}_t=[-c_t+w_tl_t+r(\tilde{k}_t-\phi)]dt  \\
 &c_t\geq 0 \notag\\
 &\tilde{k}_t\geq 0 \text{ a.s.}\notag
\end{align}

\textbf{Representative's firm problem}

The firms has to solve 
\begin{equation}
 \underset{K^N_t,L^N_t}{max}\Pi(K^N_t,L^N_t)= f(K_t,L_t)-(r_t+\delta)K_t-w_tL_t 
\end{equation}
with 

 \begin{equation*}
  \left\{
  \begin{split}
  \tilde{K}^N_t=\frac{1}{N}\sum_{i=1}^N\tilde{k}_t^i=\int kd\mu_t^N(k,l) \\
  L^N_t=\frac{1}{N}\sum_{i=1}^Nl_t^i=\int ld\mu_t^N(k,l)
  \end{split}
  \right.
 \end{equation*}

\textbf{Mean Field Game set-up}\\

First we solve Firm's problem, which by taking first order conditions on (4.2.9) together with assumptions  ($A3$),($A4$) yield:

\begin{equation}
 \left\{
 \begin{split}
  &a(\tilde{K}^N_t)^{a-1}=r_t+\delta\\
  &(\tilde{K}^N_t)^a=w_t
 \end{split}
\right.
\end{equation}

then we inject (4.3.1) in (4.2.8) which gives us:

\begin{equation}
 \left\{
 \begin{split}
  &d\tilde{k}_t^i=[(\tilde{K}^N_t)^al_t^i+(a(\tilde{K}^N_t)^{a-1}-\delta)(\tilde{k}_t^i-\phi)-c_t^i]dt\\
  &dl_t^i=b(l^i_t)dt+\sigma(l_t^i)dW_t^i\\
  &\tilde{K}^N_t=\frac{1}{N}\sum_{i=1}^N\tilde{k}_t^i=\int kd\mu_t^N(k,l) \\
  &\mu_t^N=\frac{1}{N}\sum_{i=1}^N \delta_{X_t^i\leq x}
 \end{split}
\right.
\end{equation}

the mean capital $\tilde{K}^N_t$ in the above equation gives us the mean filed interactions. 

Now as explain already in previous chapters we send the number of agents to infinity and try to solve the representative's agent problem. 

\section{Solution of Aiyagari MFG model }

Once we have $N\to\infty$ we would like the flow of  empirical measures $\mu_t^N$ to converge to some flow $\mu_t$ by a law of large numbers.  

We restate the optimal control problem for the representative household. 

\begin{equation}
 \underset{c_t\in A_{adm}}{max}\mathbb{E}_0\big[\int_0^Te^{-\beta t}U(c_t)dt+ \tilde{U}(c_T)\big]
\end{equation}
subject to 

\begin{equation}
 \left\{
 \begin{split}
  &d\tilde{k}_t=[\bar{K}_t^al_t+(a\bar{K}_t^{a-1}-\delta)(\tilde{k}_t-\phi)-c_t]dt\\
  &dl_t=b(l_t)dt+\sigma(l_t)dW_t\\
  &\bar{K}_t=\int kd\mu_t(k,l) \\
  &\mu_t=\mathcal{L}(X_t)
 \end{split}
\right.
\end{equation}

Our strategy here is going to be the same as in chapter 3 we are going to solve the optimal control problem when $\mu_t$ is fixed, but here things are a little bit easier since the state dynamics depend only on the mean capital. 

\begin{description}
 \item[Minimization of the Hamiltonian]
 
 $\vspace{1mm}$
  First we define the Hamiltonian 
  $$H(t,c,k,l,y_k,y_l,q)=U(c)+\langle[\bar{K}_t^al+(a\bar{K}_t^{a-1}-\delta)(k-\phi)-c],y_k\rangle+\langle b(l,t),y_l\rangle +\langle\sigma(l,t),q \rangle $$
  
  we take 
  $$ \frac{\partial H}{\partial c}=0 $$
  \begin{align}
   &U'(c)-y_k=0 \notag\\
   &\hat{c}=y_k^{-\frac{1}{\gamma}}
  \end{align}

  $\hat{c}$ is the optimal control rule which is  independent of $l$ and $y_l$ and so we can simplify our formulas by removing the from the definition of the Hamiltonian and the adjoint process and drop the subscript $k$ from  the $y$ variable. This way we end up with one dimensional deterministic adjoint process since we ruled out the stochastic part.
 
 \item[Adjoint process and forward-backward system ]
  We define the first order adjoint process 
  \begin{equation*}
   \left\{
   \begin{split}
   &dY_t=-\frac{\partial H}{\partial x}dt=-(a\bar{K}_t^{1-a}-\delta)Y_tdt\\
   &Y_T=1
   \end{split}
   \right.
  \end{equation*}

  Together with the state equation (4.3.2) (but independent of $l_t$) we end up with the forward-backward system 
  
  \begin{equation}
 \left\{
 \begin{split}
 &dY_t=-(a\bar{K}_t^{1-a}-\delta)Y_tdt\\
  &d\tilde{k}_t=[\bar{K}_t^al_t+(a\bar{K}_t^{a-1}-\delta)(\tilde{k}_t-\phi)-\hat{c}_t]dt\\
  &Y_T=1\\
  &\bar{K}_t=\int kd\mu_t(k,l) \\
  &\mu_t=\mathcal{L}(X_t)
 \end{split}
\right.
\end{equation}
  
  in order to calculate the mean filed interactions we take the expectation in $\tilde{k}_t$,  remembering also assumption ($A4$) so we end up in a system of ODEs 
  
  \begin{align}
 &dY_t=-(a\bar{K}_t^{1-a}-\delta)Y_tdt\\
  &d\bar{K}_t=\bar{K}_t^a-a\phi\bar{K}_t^{a-1}-\delta\bar{K}_t+\delta\phi)-Y_t^{-\frac{1}{\gamma}}dt\\
  &Y_T=1
\end{align}
   
\end{description}

We have to work using numerical approximations to solve (4.3.5)(4.3.6)

\textbf{To be continued...}

\begin{re2}
 \begin{enumerate}
  \item System (4.3.4) is a direct consequence of assumption ($A4$) which restricted mean filed interactions to the capital and so we could use ODE  methods to explore it.
  \item Equation (4.3.5) can be solved explicitly and together with the numerical approximation for the system we can substitute in the optimal control rule (4.3.3) and get the optimal consumption rule for the economy. This is a very important variable for economists since they can study the growth path of the economy and decide about optimal  macroeconomic policies.  
 \end{enumerate}

\end{re2}

\appendix

\chapter{Optimal Control}

\section{Introduction}
We are interested in studying a phenomena that can be described by a set of variables called $\mathit{state \hspace{1mm} variables}$ and a system of differential equations ($\mathit{dynamical \hspace{1mm} system} $) which define the path in which the state variables evolve.\\
We are interested in answering the following questions:
\begin{itemize}
\item
What is the asymptotic behaviour of our system?
\item
Can we add specific variables which we have under our control to the system to steer it to a target set?(Controllability)
\item
Can we find a path ($\mathit{trajectory} $) which makes a certain pre-decided criterion optimal?(Existence of optimal control)
\item
How can we design the variables(controls) to achieve this optimum? (Approximation of the optimal controls)
\end{itemize}

Let's introduce some notation:
\begin{d1}{Some terminology}
 
\begin{enumerate}
 
\item
$x \in \mathbb{R}^n$ be the state variables
\item
$\Omega$ be the unit cube in  $\mathbb{R}^m$ i.e. $\Omega := \{c| c\in \mathbb{R}^m, |c|\leq 1, i=1,2,...,m \}  $
\item
 $u(\cdot) \in \mathcal{U}_m$ the controls belonging to a set of measurable functions i.e. $ \mathcal{U}_m[0,t_1]=\{ u(\cdot)|u(\cdot)\in \Omega \, and \,  u(\cdot)\, measurable\, on\, [0,t_1] \} $
\item
$\mathcal{T}(t)$ the target set 
\item 
$\dot{x}(t)=f(x(t),u(t))$ the dynamics of the system under control $u(t)$ , $x(t_0)=x_0$ 
\item 
$x[t]\equiv x(t;x_0,u(\cdot))$ the response, i.e. the solution id the dynamical system when using the control $u(\cdot)$ 
\item
$J[u(\cdot)]=\int_0^{t_1}f^0(x[t],u(t))dt$ the criterion or value or cost function under which we are interested in finding the optimal path  

\end{enumerate}

\end{d1}

\begin{re2}
 \begin{enumerate}
 \item In this appendix we are going to calligraphic capital letters for our sets, for emphasis. 
 
  \item $\mathcal{T}(t)$ represents a time varying set in which we would like our response $x[t]$ to be included 
  
 \end{enumerate}

\end{re2}

\begin{description}
 \item[Optimal control problem] 

We want to find the control $u(\cdot)$ (if there exist one) which steers the system 
\begin{equation}
\dot{x}(t)=f(x(t),u(t) ) \label{1}
\end{equation}
in a way that  $x[t_1]\in \mathcal{T}(t_1)$ with the minimum cost(or maximum value) 
\begin{equation}
 J[u(\cdot)]=\int_0^{t_1}f^0(x[t],u(t))dt  \label{2}
\end{equation}

\end{description}

\section{Controllability}

Now in order to solve our basic optimal control problem we turn to the controllability question.

\begin{d1}{Controllable set}
The set
$$\mathcal{C}(t)=\{x_0\in \mathbb{R}^n|\exists u(\cdot)\in \mathcal{U}_m \hspace{1mm} such \hspace{1mm}  that \hspace{1mm}  x(t;x_0,u(\cdot))\in \mathcal{T}(t)\}$$
contains all states which can be steered to the target at time t
\end{d1}

The basic questions that arise are from the above definition are:
\begin{enumerate}
 \item 
to describe $\mathcal{C}$
\item
to show how $\mathcal{C}$ changes if we use special classes of controls

\end{enumerate}

Two desirable properties of $\mathcal{C}$ are: 
\begin{itemize}
\item
$0 \in Int\mathcal{C}$
\item
$ \mathcal{C}=\mathbb{R}^n$ in this case the system is completely controllable
\end{itemize}

\begin{d1}{Reachable set and Reachable cone}
\begin{itemize}
\item
The set
$$\mathcal{K}(t;x_0)=\{x(t;x_0,u(\cdot))|u(\cdot)\in \mathcal{U}_m \}$$
contains all states which can be reached in $\mathbb{R}^n$ at time t, from initial point $x_0$ and is called $\mathit{reachable \hspace{1mm} set}$
\item
and the set
$$ \mathcal{RC}(x_0)= \{ t,x(t;x_0,u(\cdot))|t\leq0 ,u(\cdot)\in \mathcal{U}_m \}=\bigcup_{t\geq0}\{t\}\times \mathcal{K}(t;x_0)   $$
 is called $\mathit{reachable \hspace{1mm} cone}$ 

\end{itemize}
\end{d1}

There exists a connection between reachable sets and controllable sets via the time reversed dynamical system $x(t)$ solves \eqref{1} with $x(0)=x_0$ and $x(t_1)=x_1$ if and only if $z(t)=x(t_1-t)$ solves:
\begin{equation}
\dot{z}(t)=-f(z,\tilde{u})  
\end{equation}
$$z(0)=x_1, z(t_1)=x_0, \tilde{u}(t)=u(t_1-t)$$

The two systems have the same trajectories, traversed in opposite directions

\begin{t1}

For the system \eqref{1} $\mathcal{C} $ is arc-wise connected. $\mathcal{C} $ is open if and only if $0\in Int\mathcal{C} $
\end{t1}

\begin{re1}
 
A set (or a topological space) X is arc-wise connected if $\exists f:[0,1]\to X $ s.t. $f(0)=a$ and  $f(1)=b$ for $a,b \in X$ with continuous inverse

\end{re1}

Again we will investigate \eqref{1} with the  extra assumption that $f(x,u)$ is continuously differentiable in x,u and $f(0,0)=0\in \mathbb{R}$ Therefore expand $f(x,u)$  about $(0,0)$ 
$$f(x,u)=f_x(0,0)x+f_u(0,0)u +o(|x|+|u|)$$
with $f_x, f_u$ the appropriate Jacobian matrices.

We expect the controllability of the nonlinear \eqref{1} near $0\in \mathbb{R}$ to be determined by the controllability of the linearisation:
$$\dot{x}= f_x(0,0)x+f_u(0,0)u=A_fx+B_fu$$
and define the controllability matrix:
$$M_f=\{ B_f,A_fB_f,A_f^2B_f,\dots A_f^{n-1}B_f  \}$$

\begin{t1}
If $rankM_f=n$ then $ 0\in Int\mathcal{C}$ for \eqref{1}

\end{t1}

\begin{t1}
For \eqref{1} suppose $rankM_f=n$ if solution $x(t)=0$ of the free system $\dot{x}=f(x,0)$ is globally asymptotically stable then $\mathcal{C}=\mathbb{R}^n$ for \eqref{1}
\end{t1}

\begin{re1}
 
We can use the Hartman-Grobeman theorem for topological equivalence of the linearised and the nonlinear system 

\end{re1}

There are three subsets of $\mathcal{U}_m$ that have some interest:
\begin{itemize}
\item Piecwise Constant
$$\mathcal{U}_{PC}[0,t_1]=\{u(\cdot)\in\mathcal{U}_m[0,t_1]| u(\cdot) piecewise\hspace{1mm} constant \hspace{1mm} on[0,t_1]\}$$
\item Absolutely continuous
$$\mathcal{U}_{AC}[0,t_1]=\{u(\cdot)\in\mathcal{U}_m[0,t_1]|  u(\cdot) absolutely \hspace{1mm} continuous,$$  $u(0)=u(t_1)=0 |u(t)|\leq1 \hspace{1mm} and \hspace{1mm} |\dot{u}(t)\leq \epsilon \hspace{1mm} a.e.\hspace{1mm}  on [0,t_1]\}$
\item Bang-bang (uses full power)
$$\mathcal{U}_{BB}[0,t_1]=\{u(\cdot)\in\mathcal{U}_m[0,t_1]||u^i(t)|=1 [0,t_1],i=1,...m\}$$
\end{itemize}

\section{Existence of optimal controls}

 Existence theory is in a nutshell is a study of a continuous or lower semicontinuous function $C[u(\cdot)]$ on a compact( in some sense) set of controls $\mathcal{U}_m$

 The problem \eqref{1}-\eqref{2} is in essence a mapping 
$$ C(u(\cdot))\to C[u(\cdot)]   $$
from $\mathcal{U}_m$ into $\mathbb{R}$ \\
This mapping can be extremely complicated since the cost functional $C[u(\cdot)]$ usually involves the response $x[\cdot] $ \\
The general approach should be:
\begin{enumerate}

 \item Show that $C[u(\cdot)]$ is bounded below, hence there exists a minimizing sequence $\{u_k(\cdot)\}_{k\in\mathbb{N}}$ with associated responses $\{x_k[\cdot]\}_{k\in\mathbb{N}}$ 
 \item Show that $\{x_k\}_{k\in\mathbb{N}}$ to a limit $x_*[\cdot]$ (not necessarily a response)
 \item Show that there is a $u_*(\cdot)\in\mathcal{U}_m$ for which $x_*[\cdot]$ is a response

\end{enumerate}

\begin{t1}{Existance}
For the problem \eqref{1}-\eqref{2} on a fixed interval $[0,T]$ with: $x_0$ given, $\mathcal{T}(t)=0$, $f(t,x,u)$ and $f^0(t,x,u)$ continuous. Assume:
\begin{enumerate}

 \item that the class of admissible controls which steer $x_0$ to the target set in time $t_1$ is nonempty
 \item satisfy an a priori bound:
$$ |x(t;x_0,u(\cdot))|\leq a \hspace{1mm} \forall \hspace{1mm} u \hspace{1mm} admissible $$
\item the set of points $f^0(t,x,\Omega)=\{(f^0(t,x,v),f^T(t,x,v))^T|v\in \Omega\}$ is convex in $\mathbb{R}^{n+1}$ 

\end{enumerate}

Then there exists and optimal control
\end{t1}

\section{Pontryagin's Maximum Principle}

In the previous section we gave the sufficient conditions about the existence of at least one optimal control. Here we are interested in the necessary conditions, which collectively are known as the Potryangin Maximum Principle.  

In this section we suppose the target set is $\mathcal{T}(t)=x_1$ and the cost is $C[u(\cdot)]=\int_{t_0}^{t_1}f^0(x[t],u(t))dt $ where $t_1$ is $\mathbf{unspecified}$

\begin{d1}{Dynamic cost variable}
We define as 
$x^0[t]=\int_{t_0}^{t_1}f^0(x[s],u(s))ds$ the dynamic cost variable
\end{d1}

\begin{re1}

If $u(\cdot)$ is optimal, then $x^0[t_1]$ is as small as possible
 
\end{re1}

If we set $\hat{x}[t]=(x^0, x^T[t])^T$ and $\hat{f}(t,\hat{x})=(f^0,f^T)^T $ then our original problem can be restated as:

\begin{description}
 \item[Restatement of the original problem]
 
Find an admissible control $u(\cdot)$ such that the (n+1)-dim solution of
\begin{equation}
 \dot{\hat{x}}[t]=\hat{f}(\hat{x},u(t)) 
\end{equation}
terminates at 
$\bigl(\begin{matrix}
x^0[t_1]\\x_1 
\end{matrix} \bigr)$ with $x^0[t_1]$ as small as possible.
\end{description}

In the linear case $\dot{x}(t)=Ax(t)+Bu(t)$ with cost function $C[u(\cdot)]=\int_0^{t_1}dt=t_1$ we know that the optimal control is going to be extremal(there exists a supporting hyperplane). We would like to use the same mechanism in the general nonlinear case. So we make the following definitions

For a given constant control $u(\cdot)$ any solution $\hat{x}[\cdot]$ of $\dot{\hat{x}}[t]=\hat{f}(x[t],u(t))$ is a curve in $\mathbb{R}^{n+1}$ If $\hat{b}_0$ is a tangent vector to $\hat{x}[\cdot]$ at $\hat{x}[t_0]$ then the solution $\hat{b}(t)$ of the linearised equation:
$$\dot{\hat{b}}(t)=\hat{f}_x(x[t],u(t))\hat{b}(t) $$ 
will be tangent to this curve at $\hat{x}[t]$ for all t. 

Thus the linearised equation describes the evolution of tangent vectors along the solution curves of the resulting autonomous equation.

 \begin{d1}{The Ajoint system}
 
  For a given admissible control $u(\cdot)$ and associated response $\hat{x}[\cdot]$  we consider the (n+1)-dim linear system
  \begin{equation}
   \dot{\hat{p}}(t)=-\hat{f}_{\hat{x}}(x[t],u(t))^T\hat{p}(t)\label{3}
  \end{equation}
The solutions of this system are called $\mathit{extended \hspace{1mm} costates}$
 \end{d1}
with $\hat{f}_{\hat{x}}$ the usual Jacobian matrix of $\hat{f}$ with respect to $\hat{x}$ \\
Thus if $\hat{b}(t)$ is tangent to $\hat{x}[\cdot]$ at $\hat{x}[t] $ for all t and if $\hat{p}(t_0)$ is perpendicular to $b(t_0)$ then $p(t)$ will be perpendicular to $\hat{x}[\cdot]$ at $\hat{x}[t]$ for all t.
\begin{re1}

The Adoint describes the evolution of vectors lying in the n-dim hyperplane P(t) attached to the extended response curve $\hat{x}[\cdot]$
\end{re1}

\begin{d1}{Hamiltonian}

 For a given control and extended response $(\hat{x}[\cdot],u(\cdot))$  we take any costate $\hat{p}(\cdot)$ and define the Hamiltonian as a the real-valued function of time:  
 $$H(\hat{p}, \hat{x},u)=<\hat{p},\hat{f}>= \sum_{j=0}^n p^j(t)f^j(x[t],u(t)) $$
 
\end{d1}

 and for system \eqref{3} we have:
 \begin{equation}
  \check{x}=grad_{\hat{p}}H(\hat{p},\hat{x},u)=\big(\frac{\partial H}{\partial p^0},\frac{\partial H}{\partial p^1},\cdots , \frac{\partial H}{\partial p^n} \big)^T
 \end{equation}

 \begin{equation}
  \check{p}=-grad_{\hat{x}}H(\hat{p},\hat{x},u)=-\big(\frac{\partial H}{\partial x^0},\frac{\partial H}{\partial x^1},\cdots , \frac{\partial H}{\partial x^n} \big)^T
 \end{equation}

 \begin{d1}{Legendre transform}
 $$\mathcal{H}(\hat{p},x)=\underset{v\in \Psi}{sup}H(\hat{p},x,v) $$
 $\mathcal{H}$ is the largest value of $H$ we can get for the given vectors $(\hat{p},x)$ using admissible values for $\mathbf{v}$
\end{d1}

\begin{t1}{Pontryagin Maximum Principle}

Consider the extended control problem (4) with measurable controls $u(\cdot)$ taking values in a fixed bounded set $\Psi \subset \mathbb{R}^m$ Suppose $(u(\cdot),\hat{x}[\cdot])$ is an optimal control-response pair. Then there exists an absolutely continuous function $\hat{p}(\cdot)$ solving the adjoint system a.e. on $[t_0,t_1]$ with:
\begin{equation}
 H(\hat{p}(t),x[t],u(t))=\mathcal{H}(\hat{p}(t),x[t]) 
\end{equation}

\begin{equation}
 \mathcal{H}(\hat{p}(t),x[t])=0, 
\end{equation}

\begin{equation}
 p^0(t)=p^0(t_0)\leq 0 
\end{equation}
\end{t1}

$\mathbf{Remarks}$\\
\begin{itemize}
 \item If $u(\cdot)$ is optimal for (4), then there is an associated response-adjoint pair, $(\hat{x}[\cdot],\hat{p}[\cdot])$ such that  for each t $H(\hat{p}(t),x[t],v)\leq 0$ for any $v\in \Psi$ 
 \item The PMP assumes that an optimal control exists. There maybe be a non empty set of candidates and yet no optimal control for a given problem. 
\end{itemize}

\section{Hamilton-Jacobi Equation}

$\mathbf{Note}$\\
We will change the notation to be closer to the PDE literature. We will use u for the solution of the Hamilton Jacobi and other letters for controls whereas needed

 \begin{d1}{Hamilton Jacobi}
 
 The partial differential equation 
  \begin{equation}
   \frac{\partial u(x,t)}{\partial t}+ H(D_x u(t,x))=0 \text{ in $\mathbb{R}^n\times(0,\infty)$} 
\end{equation}
\begin{equation}  
 u=g \text{ on $\mathbb{R}^n \times\{t=0\}$ }
 \end{equation}

  is called $\mathit{Hamilton \hspace{1mm} Jacobi \hspace{1mm} equation}$ 
 \end{d1}

\subsection{Derivation of HJE using calculus of variations}

Let $L:\mathbb{R}^n \times \mathbb{R}^n \to \mathbb{R}^n \to \mathbb{R} $ be the Lagrangian with $L=L(q,x)$, $q,x \in \mathbb{R}^n$ (q represents velocity, x represents state)
$$
 \left\{
 \begin{split}
D_qL=L_{q_1},L_{q_2}...L_{q_n}\\
D_xL=L_{x_1},L_{x_2}...L_{x_n}
 \end{split}
\right.
$$

We introduce the action functional
$$I[w(\cdot)]=\int_0^t L\big(\dot{w}(s),w(s)\big)ds$$ for $w(\cdot)$ belonging to the admissible class $\mathcal{A}=\{w(\cdot)\in C^2|w(o)=y,w(t)=x\} $

After defining the action functional the basic problem in calculus of variations is to find a curve $x(\cdot)\in \mathcal{A}$ satisfying 
$$ I[x(\cdot)]=\underset{w(\cdot)\in \mathcal{A}}{min} I[w(\cdot)] $$
we are asking for a function $x(\cdot)$ which minimizes the functional $I(\cdot)$ among all admissible candidates 

We assume next that there exists a $x(\cdot)\in \mathcal{A}$ that sastisfy our calculus of variations problem and we will deduce some of its properties

\begin{t1}{Euler-Lagrange}

Given a minimizer $x(\cdot) \in \mathcal{A}$ it solves the Euler-Lagrange equations

\begin{align}
 -\frac{d}{ds}\big(D_qL(\dot{x}(s),x(s)\big)+D_xL(\dot{x}(s),x(s))=0 \text{ $0\leq s\leq t$} \tag{E-L} 
\end{align}

\end{t1}

\begin{proof}
 
Choose $v\in \mathcal{A}$ it follows that $v(0)=v(t)=0$ and for $\tau \in \mathbb{R}$ we set $w(\cdot)=x(\cdot)+\tau v(\cdot)$ \\
$w(\cdot)\text{ in } C^2$ and $w(0)=w(t)=0$ so $w \in \mathcal{A}$ and $I[x(\cdot)]\leq I[w(\cdot)]$
We set also $i(\tau)=I[x(\cdot)+\tau v(\cdot)$ we differentiate with respect to $\tau$ noticing that $i^{'}(0)=0$ and we get the result 
\end{proof}

 \begin{re1}
  
 Any minimizer solves the E-L equations but it is possible that a curve $x(\cdot)\in\mathcal{A}$ may also solve E-L without being a minimizer. In this case $x(\cdot)$ is a critical point of $I$ 
 
 \end{re1}

 We now assume that $x(\cdot)$ is a critical point of $I$ and thus solves the $E-L$. We set 
$$ p(s)=D_qL(\dot{x}(s),x(s))\text{ for $0 \leq s \leq t$}$$
$p(s)$ is called the generalized momentum 

\begin{description}
 \item[Assumption]
Suppose for all $x,p \in \mathbf{R}^n$ that the equation $p=D_qL(q,x)$ can be uniquely solved for q as a smooth function of $p$ and $x$ $q=q(x,p)$
 
\end{description}

 \begin{d1}{Hamiltonian}
 The Hamiltonian associated with the Lagrangian $L$ is 
 $$H(p,x)=pq(p,x)-L(q(p,x),x) \text{ $p,x\in\mathbb{R}^n$} $$
\end{d1}

\begin{ex1}
 
Let $L(q,x)=\frac{1}{2} m q^2-\phi(x)$ be the Lagrangian. The corresponding E-L is 
$$m\ddot{x}(s)=f(x(s)) $$ 
for $f=-D\phi$ this is Newton's Law with the force field generated by the potential $\phi$ \\
Setting $p=D_qL=m|q|$ the Hamiltonian is 
$$H(p,x)=p\frac{p}{m}-\frac{1}{2}m{\frac{p}{m}}^2+\phi(x)=\frac{1}{2}mp^2+\phi(x)  $$
The sum of kinetic and potential energies 

\end{ex1}

If we rewrite $E-L$ in terms of $p(\cdot),x(\cdot)$ we arrive in the next theorem
\begin{t1}{Hamilton's ODE}
 The functions $x(\cdot),p(\cdot)$ satisfy Hamilton's equations
 $$\dot{x}(s)=D_pH(p(s),x(s))$$
 $$\dot{p}(s)=-D_xH(p(s),x(s)) \text{ for $0\leq s \leq t$}$$
and the mapping 
$$s \to H(p(s),x(s)) \text{ is constant}$$
 \end{t1}
(The sum of kinetic and potential energy is constant and these systems are called conservative)\\

\begin{proof}{(only the third statement)}
 
$\frac{d}{ds}\big(H(p(s),x(s))\big)=\sum^n_{i=1}\frac{\partial H}{\partial p_i}\dot{p}_i+\frac{\partial H}{\partial x_i}\dot{x}_i\underset{Hamilton's ODE}{=}\sum^n_{i=1}\frac{\partial H}{\partial p_i}-\frac{\partial H}{\partial x_i}+\frac{\partial H}{\partial x_i}\frac{\partial H}{\partial p_i}=0$
 
\end{proof}

\subsection{A candidate for the HJE}

Retuning to  the initial-value problem, we will investigate a connection between the PDE and the calculus of variations.\\
If $x \in \mathbb{R}^n$ is given and $g$ appropriate initial data we should presumably try to minimize the action functional, taking into account the initial condition for the PDE. 
$$\int_0^tL(\dot{w}(s))ds+g(w(0))$$ 

Finally we can now construct a candidate for  the initial-value problem in terms of the variational principle. 

\begin{equation}
 u(x,t)=inf\{\int_0^tL(\dot{w}(s))ds+g(w(0))|w(0)=y, w(t)=x \}
\end{equation}

With $w(\cdot) \in C^2,\hspace{1mm} w(t)=x$. 

We are going to investigate the sense in which $u$ solves the initial-value problem.

\begin{description}
 \item[Assumptions] 
 (They come naturally given the previous discussion but they are not sufficient to guarantee uniqueness) 
 \begin{itemize}
  \item H is smooth, convex and $\underset{|p|\to \infty}{lim} \frac{H(p)}{|p|}=\infty $
  \item $g:\mathbb{R}^n\to \mathbb{R}$ is Lipschitz continuous
 \end{itemize}
 
\end{description}

\begin{d1}{Hopf-Lax formula}
 $$u(x,t)=\underset{y\in\mathbb{R}^n}{min}\{tL\big(\frac{x-y}{t} \big)+g(y)  \} $$
 is the so called Hopf-Lax formula
\end{d1}

\begin{t1}

  The Hopf-Lax formula solves the minimization problem (14) 
 \end{t1}

 \begin{proof}
  
 $u(x,t)\leq \int_0^t L(\dot{w}(s))ds+g(y)$ and let us define $w(s)=\frac{s}{t}x+(1-\frac{s}{t})y =y+\frac{s}{t}(x-y)$ the convex combination of x,y with
 $\dot{w}(s)=\frac{y-x}{t}$
 $$\int_0^t L(\frac{x}{t}+(1-\frac{1}{t})y)ds=L(\frac{y-x}{t})\int_0^tds=tL(\frac{y-x}{t}) $$
 
 For the other hand-side by Jensen's inequality we get
 $$L\big(\frac{1}{t}\int_0^t\dot{w}(s)ds\big) \leq \frac{1}{t}\int_0^tL(\dot{w}(s))ds $$
 adding $g(y)$ to both sides and taking $inf$ over all $y \in \mathbb{R}^n$ we get the result. We have also to show that the $inf$ belongs to the set so it is actually a minimum.

 \end{proof}

 \begin{re1}{Convex duality of the Lagrangian and the Hamiltonian}   

 We hereafter suppose the Lagrangian $L:\mathbb{R}^n \to \mathbb{R}$ satisfies:
 \begin{itemize}
  \item $q \to L(q)$ is convex
  \item $\underset{q\to \infty}{lim}=\frac{L(q)}{|q|}=\infty $
 \end{itemize}
   
 \end{re1}

\begin{d1}

 The Legendre transform of L is: 
 $$L^*(p)=\underset{q \in \mathbb{R}^n}{sup}\{pq-L(q)\} $$
\end{d1}
$H=L^*$ the Hamiltonian is the Legendre transform of the Lagrangian and vice versa $L=H^*$ (Theorem)

 \begin{t1}{{Solution of HJ PDE}}
 
  The function $u$ defined by the Hopf-Lax formula is Lipschitz, differentiable a.e. in $\mathbb{R}^n\times(0,\infty)$ and solves the initial value problem
 \end{t1}
 
\begin{proof}
 
First we prove the theorem for a point $(x,t)$ where $u$ is differentiable by constructing the PDE. After we use Rodemacher's theorem to extend the result a.e. We will use the following lemma and double nesting.

\begin{l1}{Lemma}
 $$u(x,t)=\underset{y\in\mathbb{R}^n}{min}\{(t-s)L\big(\frac{x-y}{t-s}\big)+u(y,s)\}$$
\end{l1}

In other words to compute $u(\cdot,t)$ we calculate $u$ at the time $s$ then we use $u(\cdot,s)$ as the initial condition for the remaining time interval
  $$u(x+hq,t+h)=\underset{y\in \mathbb{R}^n}{min}\big\{hL\big(\frac{x+hq-y}{h}\big)+ u(x,t)\big\}\leq hL(q)+u(x,t) $$
 hence
 $$\frac{u(x+hq,t+h)-u(x,t)}{h} \leq L(q) $$
 $$\frac{u(x+hq,t+h)-u(x+hq,t)+u(x+hq,t)-u(x,t)}{h} \leq L(q) $$
 and for $h\to 0^+$ we get $(x\in \mathbb{R}^n)$
 $$qDu(x,t)+u_t(x,t)\leq L(q) \text{ for all $q \in \mathbb{R}^n$}$$
 $$\underset{q \in \mathbb{R}^n}{max}\{qDu(x,t)-L(q)\}+u_t(x,t)\leq 0 $$
 and finally by remembering the Legendre transform we arrive in
 $$H(Du(x,t))+u_t(x,t)\leq 0 $$
 
 For the other hand-side we need to consider the differences $u(x,t)-u(y,s)$
with $s=t-h$, $y=\frac{s}{t}x+(1-\frac{s}{t}z)$ \\and $z$  s.t. $u(x,t)\geq tL(\frac{x-z}{t})+g(z)$ $\frac{x-z}{t}=\frac{y-z}{t}$
$$u(x,t)-u(y,s) \geq tL(\frac{x-z}{t})+g(z)-[sL(\frac{y-z}{t})+g(z)]=(t-s)L(\frac{x-z}{t}) $$
$$\frac{u(x,t)-u(y,s)}{t-s} \geq L(\frac{x-z}{t}) $$ 
we change the $y,s$ variables so
$$\frac{u(x,t)-u((1-\frac{h}{t})x+\frac{h}{t}z, t-h)}{h}  \geq L(\frac{x-z}{t})$$
we add and subtract $u((1-\frac{h}{t})x+\frac{h}{t}z, t) $ to form the derivative and sending $h \to 0^+$ we get

$$\frac{x-z}{t} Du(x,t)+u_t\geq L(\frac{x-z}{t})$$
and again using the Legendre transform we get the PDE\\
So $$Du(x,t)+u_t(x,t)=0$$ for fixed $(x,t)$ and using Rodemacher's theorem we extend the differentiability of $u$ a.e.\\
$\mathbf{Reminder}$\\  
$(\mathbf{Rodemacher})$
Let $u$ be locally Lipschitz continuous in $U\subset \mathbb{R}^n$. Then u is differentiable almost everywhere in U

\end{proof}

\begin{ex1}{Counter Example of uniqueness}

$$u_t +|u_x|^2=0 \text{ in $\mathbb{R}^n\times(0,\infty)$}$$
$$u=0 \text{ in $\mathbb{R}^n\times \{t=0 \}$}$$
This initial value problem admits more than one solution i.e.
$$u_1(x,t)=0 $$
and
$$u_2(x,t)=
\begin{cases}
0 & \text{ if $|x|\geq t$}\\
x-t & \text{ if $0\leq x\leq t$}\\
-x-t & \text{ if $-t\leq x \leq 0$}
\end{cases}
$$

\end{ex1}

We need stronger assumptions to get uniqueness of the weak solution as the next theorem proposes
 
\begin{d1}{Semiconcavity and Uniform convexity}

We define the following notions:
\begin{itemize}
\item
$\mathbf{Semiconcavity}$\\
A function $u$ is called semiconcave if there exists a $C\in \mathbb{R}$ s.t.
$$ g(x+z)-2g(x)+g(x-z)\leq C|z|^2  \text{ for all $x,z\in \mathbb{R}^n$} $$
 \item
$\mathbf{Uniform\hspace{1mm} convexity}$\\
A $C^2$ function $H:\mathbb{R}^n\to\mathbb{R}$ is called uniformly convex (with constant $\theta\geq 0$) if
   $$ \sum_{i,j=1}^n H_{p_i,p_j}(p)\xi_i \xi_j\geq \theta |\xi|^2 \text{ for all $p,\xi \in \mathbb{R}^n$} $$

\end{itemize}
\end{d1}

\begin{t1}{Uniqueness HJE}

  Suppose $H$ is $C^2$ and satisfies the assumptions made earlier along with $g$. If either $g$ is semiconcave or $H$ is uniformly convex the $u$ defined by the Hopf-Lax formula is the only weak solution of the initial-value problem   
 \end{t1}

\section{Dynamic Programming Principle}

 Here we will derive a connection between the HJE and control problems
 In the rest of the presentation we will use $\alpha$ for for the controls and $\mathcal{A}$ for the class of the admissible controls. We define the value function as:
 $$u(x,t)=\underset{\alpha(\cdot)\in \mathcal{A}}{inf}I[\alpha(\cdot)]$$
 The least cost $\mathbf{given}$ we start at $x$ at time $t$.\\
 In essence we are embedding our given control problem into a larger class of problems.\\
 The idea is to show that $u$ solves a certain HJE and conversely that a solution of this PDE helps synthesize an optimal (feedback) control.

 \begin{t1}
 
  The value function $u$ is the unique (viscosity) solution of the $\mathit{terminal}$ value problem for the Hamilton-Jacobi equation:
  \begin{equation}
  \left\{
  \begin{split}
   & u_t+H(Du,x)=0 \text{ in $\mathbb{R}^n\times(0,\infty)$} \\
   &u=g \text{ on $\mathbb{R}^n \times\{t=T\}$ }
 \end{split}
 \right.
\end{equation}
with $H(p,x)=\underset{\alpha \in \mathcal{A}}{min}\{f(x,\alpha)p+f^0(x,\alpha) \}$ ($p,x\in\mathbb{R}^n$)
 \end{t1}
 
 \begin{re2}
  
\begin{enumerate}
 \item 

 If $u$ is the (viscosity) solution of the above problem then $w(x,t)=u(x,T-t)$ is the (viscosity) solution of the initial-value problem
\begin{equation}
  \left\{
  \begin{split}
&w_t-H(Dw,x)=0 \text{ in $\mathbb{R}^n\times(0,\infty)$} \\
&w=g \text{ on $\mathbb{R}^n \times\{t=0\}$ }
\end{split}
 \right.
\end{equation}
\item 
$H(x,p)=\underset{\alpha \in \mathcal{A}}{max}H(x,p,\alpha)=\underset{\alpha \in \mathcal{A}}{max} \{ f(x,\alpha)p+f^0(x,\alpha)   \} $ 

\end{enumerate}
 
 \end{re2}

 \subsection{Dynamic Programming Principle}
 
  Here with the DPP we use HJE to solve the control problem. 
 \begin{enumerate}
  \item 
  We solve the HJE and thereby compute the value function $u$.
  \item
  We define for each point $x\in \mathbb{R}^n$ and each time $0\leq t\leq T$ 
  $$\alpha^*(s)=\alpha\in\mathcal{A}$$
  $$\alpha=argmax \{u_t(x,t)+f(x,\alpha)D_xu(x,t)+f^0(x,\alpha)\} $$
  \item 
  Next(assuming $\alpha(\cdot,t)$ is sufficiently regular) we solve the ODE:
  $$\dot{x}^*(s)=f(x^*(s,\alpha(x^*(s),s)) \text{ $t\leq s \leq T$} $$
  $$x(t)=x$$
  and define the feedback control 
  $$\alpha^*(s)=\alpha(x^*(s),s)$$
 \end{enumerate}

We need also a so called verification theorem to prove that $\alpha^*(s)$ is indeed an optimal control.

 \begin{t1}{Verification Theorem}
 
  The control $\alpha^*$ defined by the DPP is optimal.
 \end{t1}

\chapter{Stochastic Optimal Control}

\section{Introduction}

We are interested in the stochastic version of the control problem we discussed in previous chapter and this reads as follows

\begin{d1}
	 Given an SDE 
	 \begin{equation}
	 \left\{
	 \begin{split}
	  &dX_t=b(t,X_t,u_t)dt+\sigma(t,X_t,u_t)dW_t\\
	  &X_0=x \in \mathbb{R}^n
	 \end{split}
\right.
	 \end{equation}

	 and a payoff functional 
	 \begin{equation}
	  J[u_{(\cdot)}]=\mathbb{E}\big[\int_0^Tf(t,X_t,u_t)dt+g(X_T\big]
	 \end{equation}
	we are interested in finding an optimal pair (if there is one) $X_t,u_t$ that makes the payoff functional optimal(max, or min).  
\end{d1}

The goal is to optimize the criterion by selecting a non-anticipative decision among the ones that satisfying all the constrains \\
But in this particular setting where the dynamics are described by an SDE and $Xt$ is a stochastic process (and probably $u_t$), we need to define it properly by a probability space $(\Omega,\mathcal{F},\mathbb{P}, \{\mathcal{F}_t\}_{t\geq0})$ where we can define m-dim Brownian motion $W_t$.   

We can make also the following remarks to motivate further definitions
\begin{itemize}
 \item 
At any time we need to determine which information is available to the controller,(the easy answer is at most $\mathcal{F}_t$, he should not be able to foretell what is going to happen afterwards) but we will see that the ``flow'' of information can be subject to modification. 
\item
The control can be either a deterministic function or a stochastic process. In the first case the control will not be of much use because the Ito integral of a deterministic function is a Gaussian random variable. In the second case is has to be non-anticipative because otherwise the integral will not be well defined. This non-anticipative nature of the control can be represented as ``$u_{(\cdot)}$ is $\mathcal{F}_t$ adapted'' 
 \end{itemize}
 
 \subsection{Formulation}

 \begin{d1}{The strong formulation}

 Let $(\Omega,\mathcal{F},\mathbb{P}, \{\mathcal{F}_t\}_{t\geq0})$ be a filtered probability space satisfying standard conditions, let $W(t)$ be a given m-dim Brownian motion. A control $u(\cdot)$ is called strongly admissible (s-adm) and $(x(\cdot),u(\cdot))$ a s-admissible pair if:
  \begin{enumerate}
   \item 
   $u_{(\cdot)}\in \mathcal{U}[0,T]$
   $$ \mathcal{U}[0,T]:=\{u:[0,T]\times \Omega \to U| U \text{time invariant metric space}, u_{(\cdot)} \text{ measurable} \}   $$ 
   \item
   $X_{(\cdot)}$ is the unique solution of the SDE on the given probability space. (In this sense we do not distinguish between the strong and the weak solution) 
  \item
  $X_t\in S(t)\hspace{1mm} \forall t\in [0,T]$ P-a.s. where $S(t)$ is a set that vary along time (state constrains)
  \item
  $f(\cdot,X_{(\cdot)},u_{(\cdot)})\in \L_\mathcal{F}^1(0,T;\mathbb{R})$ and $g(X_T)\in L_{\mathcal{F}_T}^1(\Omega,\mathbb{R})$
  \end{enumerate}

 \end{d1}

 The set of all s-adm controls is denoted by $\mathcal{A}^s$\\

 $\vspace{1mm}$

 \begin{description}
  \item[Problem $(\mathbf{P^s})$]
  \begin{equation*}
  \underset{u_{(\cdot)}\in \mathcal{A}^s }{min}J(u_{(\cdot)})=\underset{u_{(\cdot)}\in \mathcal{A}^s }{min}\mathbb{E}\bigg[\int_0^T f(X_t,u_t,t)dt+g(X_T\bigg]
  \end{equation*}

  Subject to 
  \begin{equation*}
	 \left\{
	 \begin{split}
	  &dX_t=b(t,X_t,u_t)dt+\sigma(t,X_t,u_t)dW_t\\
	  &X_0=x \in \mathbb{R}^n
	 \end{split}
\right.
	 \end{equation*}
 \end{description}

 In certain situations it will be more convenient or necessary to vary $(\Omega,\mathcal{F},\mathbb{P}, \{\mathcal{F}_t\}_{t\geq0})$ as well as $W(\cdot)$ and consider them as part of the control 
 
 \begin{d1}{Weak formulation}
 
  A 6-tuple $\pi=(\Omega,\mathcal{F},\mathbb{P}, \{\mathcal{F}_t\}_{t\geq0},W_{(\cdot)},u_{(\cdot)})$ is called a weakly admissible control system and $(X_{(\cdot)},u_{(\cdot)})$ a w-adm pair if 
  \begin{enumerate}
   \item 
   $(\Omega,\mathcal{F},\mathbb{P}, \{\mathcal{F}_t\}_{t\geq0})$ is a filtered probability space satisfying standard conditions 
   \item
   $W_{(\cdot)}$ is an B.M. on the probability space 
   \item
   $u_{(\cdot)}$ is $\mathcal{F}_t$-adapted on $(\Omega,\mathcal{F},\mathbb{P})$ taking values in $U$, with $U$ being a time invariant metric space
   \item
   $X_{(\cdot)}$ is the unique solution of the sde on the given probability space under $u_{(\cdot)}$. (In this sense we do not distinguish between the strong and the weak solution) 
   \item
  $f(\cdot,X_{(\cdot)},u_{(\cdot)})\in \L_\mathcal{F}^1(0,T;\mathbb{R})$ and $g(X_T\in L_{\mathcal{F}_T}^1(\Omega,\mathbb{R})$
  \end{enumerate}

 \end{d1}

 Symmetrical the set of all w-adm control systems is denoted by $\mathcal{A}^s$\\

 $\vspace{1mm}$

 \begin{description}
  \item [Problem $(\mathbf{P^w})$]
 
  \begin{equation*}
  \underset{u_{(\cdot)}\in \mathcal{A}^w }{min}J(u_{(\cdot)})=\underset{u_{(\cdot)}\in \mathcal{A}^w }{min}\mathbb{E}\bigg[\int_0^T f(X_t,u_t,t)dt+g(X_T\bigg]
  \end{equation*}

  Subject to (1)
\end{description}

 $\vspace{3mm}$
 
\begin{re2}
 
\begin{itemize}
 \item 
The strong formulation stems form the practical world while weak formulation sometimes serves as an auxiliary but effective model aiming at solving problems with the strong formulation. $\mathbf{Intuition}$ The objective of a stochastic control problem is to optimize the expectation of a certain random variable that depends only on the distribution of the processes involved. Therefore if the solutions in different probability spaces have the same pdf then one has more freedom in choosing a convenient probability space to work with. 
\item
We shall make a distinction between the information available to the controller and the information about the system. We denote $\mathcal{G}_t\subset \mathcal{F}_t $ the sub-filtration of the information available to the controller i.e. $\mathcal{F}_t$ is the information of the system. The idea is that only the specific path $X_{(\cdot,\omega)}$ might be seen by the controller  
\item
It was clear relatively early in the research of stochastic control systems that in the case where we have no control over the volatility the results are parallel with those in the deterministic case  
 \end{itemize}

\end{re2}

\section{An existence result}

We will present a simplified existence proof according to   Benes \cite{Benes1970}. It has very strong and restrictive assumptions that limit lot the applicability of the result but it is relatively straightforward to follow and focuses on the important issue of the availability of information for the controller and the system. All of our work will happen under weak formulation as we are going to start from a general space of continuous functions and then change the probability measure using an extension of Girsanov's theorem to translate the canonical process of the space i.e. the Wiener process into an equivalent that would be useful for our control problem. Also we will depart slightly from our notation and use small letters for Stochastic processes and to stress the dependences. 

$\vspace{1mm}$

\begin{description}
 \item[Assumptions]
 $\vspace{1mm}$
\begin{enumerate}
\item[(A0)]
$\sigma=1$ we have no control over volatility
 \item [(A1)] 
 $b(y,u,t)$ the drift part of the SDE grows with y either slower than linearly or linearly at a slow enough rate 
 $$|b(y,u,t)|^2 \leq k(1+|y(t)|^{2a}) \text{ a<1} $$
 \item [(A2)]
 $\mathcal{G}_t=\mathcal{F}_t$ the system depends on no more than what the controller knows.
\end{enumerate}
\end{description}

\subsection{Construction of the state process}

Let $\Gamma$ be a compact metric space of control points and  $C=C[0,1]$ the space of continuous functions $y(\cdot)$ with $y:[0,1]\to \mathbb{R}^n$ . For $0\leq s \leq  t \leq 1$ we introduce a filtration $S_t$ of  $\sigma$-algebras of C-subsets generated by the sets $\{y(s)\in A|\text{ A Borel, } y(\cdot) \in C \}$. This filtration represents the knowledge of the past from 0 to t. We suppose also that the dynamics are given by a function $b:[0.1]\times C \times \Gamma \to \mathbb{R}^n$ satisfying usual assumptions. We introduce an admissible control as a function $u:[0,1]\times C \to \Gamma$ Lebesgue for $x$ and $G_t$-adapted ($G_t$ represents the information available to the controller, $G_t\subset S_t$) for t, $\mathcal{U}$ the set of admissible controls.

We assume as given a probability space $(\Omega,\mathcal{B},\mathbb{P})$ on this space is defined a n-dim Brownian motion w with continuous sample paths. There is a set $\Omega_0 \in \mathcal{B}$ of full measure such that $w(\cdot,\omega)\in C$ for $\omega \in \Omega_0$ and we define $w(\omega)=w(\cdot,\omega)$. So 
$$U:=\{y(t_1)\in A|\text{ A Borel, } y\in C \}, \hspace{1mm} U\in S_1 $$ 
$$\mathcal{W}:=\{\omega|w(t_1,\omega)\in A,\text{ A Borel} \} \hspace{1mm} \mathcal{W}\in \mathcal{B} $$
but 
$$ \mathcal{W}\cap \Omega_0=w^{-1}U $$

and so $w^{-1}S_1 \subset \mathcal{B}$ 
The classes $\mathcal{G}_t:= w^{-1}G_t$ and $\mathcal{F}_t:=w^{-1}S_t$ are filtrations and they will provide us with a way of doing all of our work in the probability space and then return for our controls to the space $C$.  

In order to construct the SDE for the dynamics of the stochastic control problem, with translation of the canonical process of $(\Omega,\mathcal{B},\mathbb{P})$ we will need the following:

\begin{d1}{Admissible drifts}

$\mathcal{A}:=\{g:[0,1]\times\Omega\to \mathbb{R}^n|g(t,\omega)=b(t,w(\omega),u(t,w(\omega))),\hspace{1mm} u\in \mathcal{U} \}$
\end{d1}

\begin{d1}{attainable densities}

$\mathcal{D}:=\{\zeta:[0,1]\times\Omega\to \mathbb{R}^n|\zeta(\omega)=e^{\zeta(g)},g\in\mathcal{A}\}$
\end{d1}

Admissible drifts are random processes while attainable densities are random variables

$\vspace{1mm}$

We will introduce the new measure 
$$ d\tilde{\mathbb{P}}=e^{\zeta(g)}d\mathbb{P} \hspace{5mm} g=b(w(t),u(t,w),t) \text{ and $\tilde{\mathbb{P}}(\Omega)=1$} $$
where 
$$ \zeta(g)_\omega=\int_0^1b(w(\omega),u(t,w(\omega)),t)dw(t)-\int_0^1|b(w(\omega),u(t,w(\omega)),t)|^2dt$$
this procedure provides a solution of (1) in a sense that under $\tilde{\mathbb{P}}$ 
$$w(t,\omega)-\int_0^tb(w(\omega),u(s,w(\omega)),s)ds=W(t,\omega) \text{ is a Wiener process}$$

If we change name $x(t,\omega)$ to $w(t,\omega)$ we have 
$$x(t,\omega)=\int_0^tb(x(\omega),u(s,\omega),s)ds + W(t,\omega) $$

The above result is based on:

$\hspace{2mm}$

\begin{t1}{Girsanov}

 Let $\phi$ be a non anticipative Brownian functional with $\phi \in L_2 $ a.s. the following are equivalent:
 \begin{enumerate}
  \item
  $w(t)-\int_0^t\phi ds$ is a Wiener process under $d\tilde{\mathbb{P}}=e^{\zeta(g)}d\mathbb{P}$
  \item
  $E[e^{\zeta(\phi + \theta)}]=1 \hspace{1mm} \forall \theta \in \mathbb{R}^n$
 \end{enumerate}

\end{t1}

$\vspace{5mm}$

Proportional to the deterministic case we will introduce the dynamic cost variable to eliminate the dependence of the criterion on the control.

We replace n-dim vector $b$ by n+1-dim vector $f,b$ and we add another 1 dim  Brownian motion $w_0$ independent of $w$ to get:
$$z=(w_0,w)=(w_0,w_1,...,w_n),\text{ $h=(f,b)$ }  $$
$$\xi_h=\int_0^1h(t,z)dz(t)-\frac{1}{2}\int_0^1|h(t,z)|^2dt $$
then under $\tilde{\mathbb{P}}$ if $E[e^\xi]=1$ then 
$$z(t)-\int_0^th(s,z)ds $$
is a n+1 dim Wiener process.

We can cover also with similar arguments the case where $x(0)=a$ the initial data is non-zero.

The following statement can give us a hint of how we can restate our problem in a more friendly form.

$\mathbf{Statement}$\\
$$\mathbb{E}[\int_0^1f(w(\omega),u(t,w(\omega),t))dt e^{\zeta(b)}]=\mathbb{E}[w_0(1)e^\xi]$$

In this manner we can restate the minimization problem as:
\begin{align}
 &min\mathbb{E}[w_0(1)e^\xi]\\
 &\text{ subject to } \notag \\
 &g(t,\omega)=b(w(\omega),u(t,\omega),t) \text{ being an admissible drift} 
\end{align}

In this form of the problem we minimize the average of the value of $x_0(\cdot)$ at the endpoint 1, the functional $e^\xi$ determines what this averaging is.

\subsection{Optimal controls}

In the deterministic control theory it was enough to assume convexity of $b(t,y,\Gamma)$ (in the case of a system $\dot{y}=b(y,u,t)$) and show that a certain function obtained as a weak limit by a compactness argument was indeed an admissible optimal control. 

In the stochastic case things are much more complicated because control can depend on available information. We have already described the structure of the available information by the appropriate $\sigma$-algebras, the problem is that the information ($\mathcal{G}_t$) which is available to the controller may differ from that on which the system depends($\mathcal{F}_t$). Unfortunately the only case that can be solved by our approach is the case $\mathcal{G}_t=\mathcal{F}_t$

Leaving out technical results we will present the main propositions for the existence of optimal control in the stochastic case. 

\begin{t1}

The following hold for problem (B.2.1)-(B.2.2):

\begin{enumerate}
 \item 
 If for each $t,u$, $b(\cdot,t,u)$ is $G_t$-measurable and if for each $t,y$, $b(y,t,\Gamma)$ is convex then $\mathcal{A}$ is convex.  
\item
If $G_t=S_t$ and if $b(y,t,\Gamma)$ is convex for $t,y \in[0,1]\times C$  then $\mathcal{D}$ is convex
\item
If $|b(y,u,t)|^2 \leq k(1+|y(t)|^{2a}) \text{ a<1}$ then $\mathcal{D}$ is a bounded set of $L_2$
 \end{enumerate}

\end{t1}

The previous theorem stems directly from our assumptions

When $\mathcal{D}$ is a bounded subset of $L_2$ the following closure and existence results are proved in a natural way using strong and weak $L_2$-topologies. 

\begin{t1}
 
 $L_2 \cap \mathcal{D}$ is closed in $L_2$-norm topology
\end{t1}

\begin{t1}{Existence of an optimal control}
 
 If $G_t=S_t$, $b(y,t,\Gamma)$ is convex and $\mathcal{D}$ is $L_2$-bounded, then an optimal control exists.
\end{t1}

\subsection{Reachable set of stochastic control systems}

\section{Stochastic Maximum Principle}

We come now to the necessary conditions for an optimal control, which collectively are known as the stochastic maximum principle. Unlike the previous section where we limited ourselves under strong assumptions for educational purposes and simplicity, here we will treat a more general case applicable to a large class of problems.  

We consider the stochastic control system:
 \begin{equation}
	 \left\{
	 \begin{split}
	  &dX_t=b(X_t,u_t,t)dt+\sigma(X_t,u_t,t)dW_t\\
	  &X_0=x \in \mathbb{R}^n
	 \end{split}
\right.
	 \end{equation}
and cost
\begin{equation}
 J[u_{(\cdot)}]=\mathbb{E}\big[\int_0^T f(t,X_t,u_t)dt+g(X_T)\big]
\end{equation}

We will make the following assumptions

\begin{description}
 \item[Assumptions]
 $\vspace{1mm}$
\begin{enumerate}
 \item[(S0)]  
 $\{\mathcal{F}_t\}_{t\leq0}$ is the natural filtration generated by $W(t)$ augmented by all the $\mathbb{P}$-null sets in $\mathcal{F}$
 \item[(S1)]
 $(U,d)$ is a separable metric space and $T\leq0$
 \item[(S2)]
 The maps $b,\sigma,f,h$ are measurable, $\exists L>0$ and a modulus of continuity $\bar{\omega}:[0,\infty] \to [0,\infty]$ such that $b,\sigma,f,h$ satisfy Lipschitz type conditions
 \item[(S3)]
 The maps $b,\sigma,f,h$ are $C^2$ and satisfy growth conditions
 \end{enumerate}
 
\end{description}

 $$ \mathcal{U}[0,T]:=\{u:[0,T]\times \Omega \to U| \text{u is $\mathcal{F}_t-adapted$}  \}   $$
 
 Given $u_{(\cdot)}\in \mathcal{U}[0,T]$ the SDE (1)  has random coefficients
 
 \subsection{Adjoint equations}
 
 In the deterministic case we had the adjoint system that described the evolution of vectors lying in the n-dim hyperplane attached to the extended response curve. Here we will use the same mechanism introducing a pair of stochastic processes instead.

 We introduce the $\mathit{terminal}$ value problem for an SDE:
 \begin{equation}
  dp_t=-\bigg[b_x(t,\bar{X}_t,\bar{u}_t)^Tp_t+\sum_{j=1}^m\sigma_x^j(t,\bar{X}_t,\bar{u}_t)^Tq^j_t-f_x(t,\bar{x}_t,\bar{u}_t)\bigg]dt+q_tdW_t
 \end{equation}
 \begin{equation}
  p_T=-h_x(\bar{X}_T)
 \end{equation}

 This is a Backward Stochastic Differential Equation (BSDE) of first order. Any pair $(p_{(\cdot)},q_{(\cdot)})\in L_{\mathcal{F}}^2(0,T;\mathbb{R}^n)\times (L_{\mathcal{F}}^2(0,T;\mathbb{R}^n))^m$ satisfying the BSDE is $\mathcal{F}_t$-adapted. Under our assumptions (Adj) admits a unique solution. The existence theorem is in Appendix C. 
 
 $\mathbf{Interpretation}$
 In the deterministic case $p_{(\cdot)}$ (the adjoint variable) satisfies a Backward ODE, that is equivalent to a forward equation if we reverse time   however in the stochastic case this cannot happen. In addition $p_{(\cdot)}$ corresponds to the shadow price of the resource represented by the state variable. On the other hand in the stochastic case the controller has to balance carefully the scale of the control and the impact of it to the uncertainty. If a control is going to affect the volatility of the system $p_{(\cdot)}$ does not characterize completely the trade-off between cost and control gain in an uncertain environment. Things can very quickly turn ugly in partially observed systems, or when the whole path of the state process is not available to the controller$(\mathcal{G}_t \subset \mathcal{F}_t)$.

 One has to introduce another variable to reflect the uncertainty or risk factor of the system. 
 
 \begin{equation}
 \begin{split}
&dP_t=-\bigg[b_x^TP_t+P_tb_x+\sum_{j=1}^m(\sigma_x^j)^TP_t\sigma_x^j\\
&+\sum_{j=1}^m(\sigma_x^j)^TQ^j_t+Q^j_t\sigma_x^j+H_{xx}(t,\bar{X}_t,\bar{u},p_t,q_t)\bigg]dt\\
&+\sum_{j=1}^mQ^j_tdW^j_t  
 \end{split}
 \end{equation}
\begin{equation}
 P_T=-h_{xx}(\bar{X}_T)
\end{equation}
where the Hamiltonian $H$ is defined by:
\begin{equation}
\begin{split}
 &H(t,x,u,p,q)=<p,b>+tr[q^T\sigma]-f,\\
 & \text{ $(t,x,u,p,q)\in [0,T]\times\mathbb{R}^n\times U\times\mathbb{R}^n\times\mathbb{R}^{n\times m}$}
\end{split}
\end{equation}

The above equation is also a BSDE of second order in matrix form, the solution $ (P_{(\cdot)},Q_{(\cdot)}\in L_{\mathcal{F}}^2(0,T;\mathbb{R}^{n,n})\times (L_{\mathcal{F}}^2(0,T;\mathbb{R}^{n,n}))^m $ and 
$(\bar{X}_t,\bar{u}_t,p_{(\cdot)},q_{(\cdot)},P_{(\cdot)},Q_{(\cdot)})$ is called an optimal 6-tuple (admissible 6-tuple)  

Where $\mathbb{R}^{n,n}$ is the space of all $n\times n$ real symmetric matrices with the scalar product: $<A_1,A_2>_*=tr(A_1,A_2) \forall A_1,A_2\in \mathbb{R}^{n,n}$

To get formal motivation for the first and second order adjoint processes we refer to the original proof of the SMP by S. Peng 1990 \cite{Peng1990}

The last ingredient before we state the Maximum Principle for stochastic systems is the so-called Generalized Hamiltonian. 

\begin{d1}{Generalized Hamiltonian}

 Let $H(x,p,u)$ be the classical Hamiltonian with $p_{(\cdot)}$ the adjoint process satisfying the first order (adj) we call Generalized Hamiltonian the function:
 \begin{equation}
  G(t,x,u,p,P)=H(x,p,u)+\frac{1}{2}tr\{\sigma(t,x,u)^TP_{(\sigma(t,x,u))}\}
 \end{equation}
with $P$ given by (10),(11)
\end{d1}

The term $\frac{1}{2}tr\{\sigma(t,x,u)^TP_{(\sigma(t,x,u))}\}$ reflects the risk adjustment , which must be present when the volatility depends on the control.

\begin{t1}{Stochastic Maximum Principle}

We assume (S0-S3) and $(\bar{X}_t,\bar{u}_t)$ be an optimal pair then there are pairs of processes 
\begin{equation}
 \begin{split}
  &(p_{(\cdot)},q_{(\cdot)})\in L_{\mathcal{F}}^2(0,T;\mathbb{R}^n)\times (L_{\mathcal{F}}^2(0,T;\mathbb{R}^n))^m\\
  &(P_{(\cdot)},Q_{(\cdot)})\in L_{\mathcal{F}}^2(0,T;\mathbb{R}^{n,n})\times (L_{\mathcal{F}}^2(0,T;\mathbb{R}^{n,n}))^m 
 \end{split}
\end{equation}
as defined before, satisfying the first and second order adjoint equations such that the variational inequality:
\begin{equation}
\begin{split}
 &H(t,\bar{X}_t,\bar{u}_t,p_t,q_t)-H(t,\bar{X}_t,u_t,p_t,q_t)\\
 &-\frac{1}{2}tr\{[\sigma(t,\bar{X},\bar{u})-\sigma(t,\bar{X},u)]^TP_t[\sigma(t,\bar{X},\bar{u})-\sigma(t,\bar{X},u)]\}\geq 0 
\end{split}
\end{equation}
holds
\end{t1}

 \section{Dynamic Programming}
 
 With the dynamic programming principle we are trying to solve our stochastic control problem by embedding our problem into a larger class of problems which we solve collectively. We are going to define the value function of the control problem and with it form a second order nonlinear PDE the famous Hamilton-Jacobi-Bellman equation. Under assumptions the solution of the PDE problem helps us synthesize an optimal control in feedback form.  
 
 \subsection{Principle of optimality}
 \subsubsection{Introduction}
 
 We are going to make the same assumptions as in section 3 with the addition that the U is complete and the functions involved to be continuous in (t,x,u). We are going to refer to them as (S1'-S3') for the needs of this section. Also we are going to use the weak formulation 
 
 As in the deterministic case we are going to define the value function as:
 \begin{equation}
  \left\{
  \begin{split}
  &V(t,x)=\underset{u(\cdot)\in \mathcal{U}^w[t,T]}{inf}J(t,x;u(\cdot)) \text{    $\forall (t,x)\in [0,T]\times \mathbb{R}^n$}\\
  &V(T,x)=g(x) \text{      $\forall x \in \mathbb{R}^n$ }
  \end{split}
  \right.
 \end{equation}

 $\vspace{3mm}$
 
\begin{re2}
 
\begin{enumerate}
 \item 
$V(t,x)$ exhibits continuous dependence on the parameters under proper conditions. Such dependence will be useful for approximations in cases of degenerate parabolic problems.
\item
 If we assume the the function $g$ on the boundary, along with $f$ are semiconcave they ``push'' $V(t,x)$ to be semiconcave.
\end{enumerate}
\end{re2}

\subsubsection{Dynamic Programming Equation}

We will state the Bellman's principle of dynamic programming. We begin from:
\begin{equation}
\begin{split}
 V(t,x)=\underset{u_{(\cdot)}\in \mathcal{U}^w[t,T]}{inf}\mathbb{E}\big[\int_t^{t+h}f(s,x(s,u(s)),u(s))ds+V(t+h,x(t+h))|\mathcal{F}_t\big]\\
\text{ for $t\leq t+h \leq T$}
\end{split}
\end{equation}

which is simply the sum of the running cost on $[t,t+h]$ and the minimum expected cost obtained by proceeding optimally on $[t+h,T]$ with $(t+h,x(t+h))$ as initial data.   

Also the Legendre transform of $f$ gives:
\begin{equation}
 H(t,x,p)=\underset{u_{(\cdot)}\in U}{sup}[<b,p>-f]
\end{equation}

With these remarks in mind we can prove the following theorem:

\begin{t1}{Hamilton-Jacobi-Bellman Equation}

 Assume (S1')-(S3') and $V\in C^{1,2}([0,T]\times \mathbb{R}^n)$. Then $V$ is a solution of the terminal value problem of a (possibly degenerate) second-order partial differential equation:
 \begin{equation}
  \left\{
  \begin{split}
   &-V_t+\underset{u\in U}{sup}G(t,x,-V_x,-V_{xx})=0 \text{ $(t,x)\in [0,T]\times \mathbb{R}^n $}\\
   &u|_{t=T}=g(x) \text{ $x\in \mathbb{R}^n$}
  \end{split}
  \right.
 \end{equation}

 where $G(t,x,p,P)$ is the Generalized Hamiltonian defined in the previous section. 
\end{t1}

\subsubsection{Optimal control in feedback form}

 Here with the DPP we use HJB to solve the control problem. 
 \begin{enumerate}
  \item 
  We solve the HJB and thereby compute the value function $V$.
  \item
  We define for each point $x\in \mathbb{R}^n$ and each time $0\leq t\leq T$ 
  $$u^*(s)=u\in\mathcal{U}^w$$
  $$u=argmax \{-V_t(x,t)+G(t,x,-V_x,-V_{xx})\} $$
  \item 
  Next(assuming $u(\cdot,t)$ is sufficiently regular) we 'solve' the SDE:
  $$dX^*_s=b(X^*(s,u(X^*_s,s))+\sigma(X^*(s,u(X^*_s,s) \text{ $t\leq s \leq T$} $$
  $$x_t=x$$
  and define the feedback control 
  $$u^*_s=u(x^*_s,s)$$
 \end{enumerate}

 $\mathbf{Alternative}$\\
 
 If we let $V_s(t,x)=\underset{u(\cdot)\in \mathcal{U}^s}{inf}J(u(x);\pi)$  then $V_s=W$ a natural way to proceed is to select a Markov control $\bar{u}$ s.t. for each $(t,x)$ in the corresponding sets. 
 
 $\bar{u}(x,t)\in argmax \{-V_t(x,t)+G(t,x,-V_x,-V_{xx})$ if $\bar{u}$ together with any initial data determine a process $\bar{x}(s)$ that satisfy (1) then
 
 $$\bar{u}(s)=\bar{u}(\bar{x}(s),s)$$ 
 
 Once the corresponding control system $\bar{\pi}$ is verified to be admissible, is also optimal. 
 
 The main difficulty is to show existence of $\bar{\pi}$ with the required property.

\subsection{The verification theorem}
 
 Solving an optimal control requires finding an optimal control and the corresponding state trajectory. The main motivation of introducing dynamic programming is that one might be able to construct an optimal control in feedback form via the value function. 
 
 \subsubsection{Connection between SMP and DP }
 
 In the case where $V(x,t)$ is sufficiently smooth. 
 
 \begin{t1}
 
  Let (S0'-S2') hold and $(x,s)\in [0,T)\times \mathbb{R}^n$ be fixed, $(\bar{x}(\cdot),\bar{u}(\cdot),p(\cdot),q(\cdot))$ be an optimal 4-tuple for $P^s$ and the value function $V\in C^{1,2}([0,T])\times \mathbb{R}^n$ then
  \begin{equation}
   \begin{split}
   V_t(t,\bar{X}_t)&=G\big(t,\bar{X}_t,\bar{u}_t,-V_x(t,\bar{X}_t),-V_{xx}(t,\bar{X}(t))\big)\\
    &=\underset{u\in U}{max}\{G\big(t,\bar{X}_t,u_t,-V_x(t,\bar{X}_t),-V_{xx}(t,\bar{X}_t)\big)\}\\
    &\text{a.e. $t\in [s,T]$, P-a.s.}
   \end{split}
   \end{equation}
   
Furthermore if $V\in C^{1,3}([0,T]\times \mathbb{R}^n$ and $V_{tx}$ is also continuous then 

\begin{equation}
 \left\{
 \begin{split}
  &V_x(t,\bar{X}_t)=-p_t, \text{ $\forall t\in [s,T]$ P-a.s.}\\
  &V_{xx}=(t,\bar{X}_t)\sigma(t,\bar{X}_t,\bar{u}_t)=-q_t \text{ $\forall t\in [s,T]$ P-a.s.}
 \end{split}
\right.
\end{equation}
 \end{t1}

 $\vspace{3mm}$
 
 $\mathbf{Corollary}$\\
 Along the optimal trajectory $\bar{x}(t)$ the map 
 $$t\to V(t,\bar{x}_t)+\int_s^tf(r,\bar{x}_r,\bar{u}_r)dr $$ 
 is a martingale

\chapter{Backward Stochastic Differential Equations}

\section{Introduction}

In the classical stochastic analysis we are interested in modelling the dynamics of a phenomena that is evolving in time and is subject to random perturbations. This gave birth to the classical SDEs which represent the dynamics as a sum of the deterministic part called drift term and the random part called diffusion term.
$$dX_t=\underbrace{b(t,X_t)dt}_\text{drift} + \underbrace{\sigma(t,X_t)dW_t}_\text{diffusion}  $$
Usually we start the system from a specific point $X_0=x $ and we allow the time to move forward. However, here we are interested in asking the opposite question i.e. How can we describe the dynamics if we start from a given point and start moving backwards in time? 

A crucial point is the availability of information. In the ODE and PDE world it is very easy to answer the above question we can make the transformation $t \to T-t$ and we have reversed the time (we can move across a smooth, or not so smooth curve in one direction or in the opposite without any problem). On the other hand in the SDE world when the SDEs are in Ito sense we demand the solutions to be adapted to some filtration generated by the driving process of the SDE and so if we just reverse time we would destroy the adaptability of the process. To elaborate more on the concept of adaptability we will use an example taken from Yong and Zhou "Stochastic Controls" \cite{YZ1999}.

\subsection{An illustrative example}

To begin with we assume $(\Omega,\mathcal{F},\{\mathcal{F}_t\}_{t\geq 0},\mathbb{P}) $ to be a standard filtered probability space on which we can defined an m-dim Brownian Motion $W_t$ such that $\{\mathcal{F}_t\}_{t\geq 0}$ is generated by $W$ augmented with all the $\mathbb{P}$-null sets in $\mathcal{F}$. We will keep this setting for the rest of the notes but for the sake of our example we will assume that $m=1$.

Consider the following terminal value problem of the SDE:
\begin{equation}
 \left\{
	 \begin{split}
	  &dY_t=0, \text{ $t\in[0,T]$}\\
	  &Y_T=\xi
	 \end{split}
\right.
\end{equation}

Where $\xi$ is an $L^2$ random variable with values in $\mathbb{R}$ and $\mathcal{F}_T$ measurable, $T>0$ given (we will also keep the assumption that the terminal time T is deterministic and known a priori for the rest of the notes). We want to find an $\{\mathcal{F}_t\}_{t\geq 0}$-adapted solution $Y_{(\cdot)}$. However, this is impossible since the only solution of (1) is 

\begin{equation}
Y_t=\xi \hspace{1mm} \forall t \in [0,T]  
\end{equation}

Which is not necessarily adapted, the only option is$\xi$ to be $\mathcal{F}_0$ measurable and finally a constant. Thus if we expect any $\{\mathcal{F}_t\}_{t\geq 0}$-adapted solution, we have to reformulate (C.1.1), keeping in mind that new formulation should coincide with (C.1.2) in the case $\xi$ is a non-random constant. 

We start with (C.2.2). A natural way to to make $Y_(\cdot)$ adapted is to redefine it as:
\begin{equation}
 Y_t=\mathbb{E}[\xi|\mathcal{F}_t], \text{ $t\in[0,T]$}
\end{equation}
 Then $Y_{(\cdot)}$ is adapted and satisfies the terminal condition $Y_T$ since $\xi$ is $\mathcal{F}_T$ measurable, but no longer satisfies (C.1.1). So the next step is to find a new equation to describe $Y_{(\cdot)}$ and this will come from the martingale representation theorem since $Y_t=\mathbb{E}[\xi|\mathcal{F}_t]$ is a martingale. So the theorem states that:

 \begin{t1}
  Under the above setting the $\{\mathcal{F}_t\}_{t\geq 0}$-martingale Y can be written as:
  \begin{equation}
   Y_t=Y_0+\int_0^tZ_sdW_s \hspace{1mm} \forall t \in [0,T], \mathbb{P}-a.s.
  \end{equation}
where $Z_{(\cdot)}$ is a predictable, W-integrable process. 
 \end{t1}

Then

\begin{equation}
 \xi=Y_0+\int_0^TZ_sdW_s
\end{equation}

and eliminating $Y_0$ from (C.1.4) and (C.1.5) we get
$$ Y_t-\xi=\int_0^tZ_sdW_s-\int_0^TZ_sdW_s$$
\begin{equation}
Y_t=\xi - \int_t^TZ_sdW_s  
\end{equation}

This is the so called BSDE. The process $Z_{(\cdot)}$is not a priori known and is a part of the solution. As a matter of fact the term $Z_tdW_t$ accounts for the non-adaptiveness of the original $Y_t=\xi$. And the pair $(Y_{(\cdot)},Z_{(\cdot)})$ is called an $\{\mathcal{F}_t\}_{t\geq 0}$-adapted solution. 

Also in this particular example the solution is unique. (The proof is rather straightforward we apply Ito's formula to $|Y_t|^2$ take expectations, then assume a second pair satisfies (C.1.6) and we have to show that $\mathbb{P}\{Y_t=Y'_t,\forall t\in [0,T]\text{ and} Z(t)=Z'_t\text{ a.e. } t\in [0,T]  \}=1$ )

\section{Linear and nonlinear BSDEs}

 Here we will state an existence theorem for the general linear case and for the nonlinear case with a Lipschitz condition. To save some time and space with notation, we introduce the following definition:

 \begin{d1}
  
 $L^2_\mathcal{F}(\Omega;C([0,T]);\mathbb{R}^k)$:= The set of all $L^2$, $\{\mathcal{F}_t\}_{t\geq 0}$-adapted processes with continuous paths \\
 In the rest when we use the subscript $\mathcal{F}$ we mean $\{\mathcal{F}_t\}_{t\geq 0}$-adapted and when we use $\mathcal{F}_T$ we mean only  $\mathcal{F}_T$-measurable
  
 \end{d1}

For the general linear case we study the problem:

 \begin{description}
  \item[The general linear problem]  In k dimensions
  
  \begin{equation}
  \left\{
  \begin{split}
  &dY_t=\{A(t)Y_t+\sum_{j=1}^mB_j(t)Z^j_t+f(t)\}dt+Z_tdW_t, \hspace{1mm} t\in [0,T] \\
  &Y_T=\xi
  \end{split}
\right.
  \end{equation}
where $A(\cdot),B_1(\cdot),...,B_m(\cdot):[0,T]\to \mathbb{R}^{k\times k}$ bounded, $\{\mathcal{F}_t\}_{t\geq 0}$-adapted processes, and $f\in L^2_\mathcal{F}([0,T];\mathbb{R}^k)$, $\xi \in L^2_{\mathcal{F}_T}(\Omega,\mathbb{R}^k)$($\xi$ is only $\mathcal{F}_T$-measurable by our notation)

 \begin{t1}{Existence}

 Let $A(\cdot),B_1(\cdot),...,B_m(\cdot)\in L^\infty_\mathcal{F}([0,T];\mathbb{R}^{k\times k}) $ Then for any $f\in L^2_\mathcal{F}([0,T];\mathbb{R}^k)$ and $\xi \in L^2_{\mathcal{F}_T}(\Omega,\mathbb{R}^k)$, the BSDE (C.2.1) admits a unique adapted solution $(Y_{(\cdot)},Z_{(\cdot)})\in L^2_\mathcal{F}(\Omega;C([0,T]);\mathbb{R}^k)\times L^2_\mathcal{F}([0,T];\mathbb{R}^{k\times m}) $
 \end{t1}

 $\mathbf{Proof} $\\
 .........\\

 \item[The general nonlinear problem]

  \begin{equation}
  \left\{
  \begin{split}
  &dY_t=h(t,Y_t,Z_tdt+Z_tdW_t, \hspace{1mm} t\in [0,T] \hspace{1mm} a.s. \\
  &Y_T=\xi
  \end{split}
\right.
  \end{equation}
  Where $h:[0,T] \times \mathbb{R}^k \times \mathbb{R}^{k\times m}\times \Omega \to \mathbb{R}^k$ and $\xi \in  L^2_{\mathcal{F}_T}(\Omega,\mathbb{R}^k)$

\begin{t1}{Existence}

 If for any $(y,z)\in \mathbb{R}^k \times \mathbb{R}^{k\times m}$ and $h(t,y,z) \hspace{1mm} \{\mathcal{F}_t\}_{t\geq 0}$-adapted with $h(\cdot, 0,0)\in L^2_\mathcal{F}([0,T];\mathbb{R}^k)$ there exists a $L>0$ such that:
 \begin{equation}
 \begin{split}
  &|h(t,y,z)-h(t,\bar{y},\bar{z}|\leq L \{|y-\bar{y}|+|z-\bar{z}| \} \\
  &\forall t\in [0,T], y,\bar{y}\in \mathbb{R}^k, z,\bar{z}\in \mathbb{R}^{k\times m} a.s.
 \end{split}
\end{equation}
Then the BSDE (C.2.2) admits a unique adapted solution\\ $(Y_{(\cdot)},Z_{(\cdot)})\in L^2_\mathcal{F}(\Omega;C([0,T]);\mathbb{R}^k)\times L^2_\mathcal{F}([0,T];\mathbb{R}^{k\times m}) $
\end{t1}

 \end{description}

\subsection{The Stochastic Maximum Principle and Duality of BSDEs and SDEs}

Here we will try to motivate a connection of SDEs and BSDEs as it appeared in the proof of the SMP. 

For starters assume we have the following stochastic control problem.

 \begin{description}
  \item[Problem]
  
  \begin{equation}
  \underset{u_{(\cdot)}\in \mathcal{A}_{adm} }{min}J(u_{(\cdot)})=\underset{u_{(\cdot)}\in \mathcal{A}_{adm} }{min}E\bigg[\int_0^T f(X_t,u_t,t)dt+g(X_T)\bigg]
  \end{equation}

  Subject to 
  \begin{equation}
	 \left\{
	 \begin{split}
	  &dX_t=b(X_t,u_t, t)dt+\sigma(X_t,t)dW_t\\
	  &X_0=x \in \mathbb{R}^n
	 \end{split}
\right.
	 \end{equation}

 \end{description}

 First, we assume $(Y_{(\cdot)},u_{(\cdot)})$ to be an optimal pair, then we introduce the so called spike variation of the control $u^\epsilon_{(\cdot)}$ and $Y^\epsilon_{(\cdot)}$ the corresponding trajectory.
 
 \begin{equation}
  u^\epsilon_t=
  \begin{cases}
   v &\text{ if  } t\in[\tau,\tau+\epsilon]\\
   u_t & \text{otherwise}
  \end{cases}
 \end{equation}
 
 Then with a little bit of effort we can get an estimate for $\Delta Y_{\tau}:= Y_\tau-Y^\epsilon_\tau $ and get the first order variational equation:
 
 \begin{equation}
 \left\{
 \begin{split}
  dY^1_t=&\{b_x(Y_t,u_t)Y^1_t+(b(Y_t,u^\epsilon_t)-b(Y_t,u_t)\}dt\\
  &+\{\sigma_x(Y_t,u_t)Y^1_t) \}dW_t\\
  Y^1_0=0
 \end{split}
\right.
  \end{equation}

using (C.2.4) we can get an estimate of the criterion using $u^\epsilon_{(\cdot)} $

 \begin{equation}
\begin{split}
J(u^\epsilon_{(\cdot)})=&\mathbb{E}[\int_0^Tf_x(Y_s,u_s)Y^1_sds]+\mathbb{E}[g_x(Y_T)]Y^1_T\\
&+\mathbb{E}[\int_0^Tf(Y_s,u^\epsilon_s)-f(Y_s,u_s)ds]+   o(\epsilon)
\end{split}
 \end{equation}

 we will use Riesz Representation theorem to exploit (C.2.8):

$\vspace{1mm}$

$\mathbf{Reminder}$

\begin{t2}{Riesz Representation Theorem}
 
 Let H be a Hilbert space, and let H* denote its dual space, consisting of all continuous linear functionals from H into the field $ \mathbb {R}$  or $\mathbb {C}$ . If ${\displaystyle x}$ is an element of H, then the function ${\displaystyle \varphi _{x}}$, for all ${\displaystyle y}$ in H defined by:

$$\displaystyle \varphi _{x}(y)=\left\langle y,x\right\rangle $$

where ${\displaystyle \langle \cdot ,\cdot \rangle }$ denotes the inner product of the Hilbert space, is an element of H*. 
\end{t2}

Here we will work with the functional:

$$I(\phi_{(\cdot)})=\mathbb{E}[\int_0^Tf_x(Y_s,u_s)y^1_sds]+\mathbb{E}[g_x(Y_T)]Y^1_T $$

and $\phi_t$ is $(b(Y_t,u^\epsilon_t)-b(Y_t,u_t)$ for notational economy. And so (since $I(\cdot)$ is linear continuous)

from Riesz there a unique $p_{(\cdot)}\in  L^2_\mathcal{F}([0,T];\mathbb{R}^k)  $ such that:
$$I(\phi(\cdot))=\mathbb{E}\int_0^T \langle p_t,\phi_t\rangle dt $$

\begin{equation}
 \mathbb{E}[\int_0^Tf_x(Y_s,u_s)Y^1_s)ds]+\mathbb{E}[g_x(Y_T)]Y^1_T)=\mathbb{E}\int_0^T \langle p_s,(b(Y_t,u^\epsilon_t)-b(Y_t,u_t)\rangle ds
\end{equation}

and by defining the Hamiltonian:
\begin{equation}
 H(x,u,p)=f(x,u)+\langle p,b(t,x,u) \rangle
\end{equation}

we can get from (C.2.9):

\begin{align}
 &\mathbb{E}\int_0^T \langle p_t,(b(Y_t,u^\epsilon_t) \rangle+f(Y_t,u_t)- \langle p_t,(b(Y_t,u_t)>-f(Y_t,u_t)dt \notag \\
 &=\mathbb{E}\int_0^TH(Y_t,u^\epsilon_t,p_t)-H(Y_t,u_t,p_t)dt
\end{align}

Finally:

\begin{equation}
\begin{split}
 H(Y_\tau,v,p_\tau)-H(Y_\tau,u_\tau,p_\tau)\geq 0\\
 \forall v\in \mathcal{A} \text{ a.e. $\mathbb{P}$-a.s.}
\end{split}
\end{equation}

This proof even though it is simple and parallel to the deterministic case can give us the important hint about how to transform the criterion and form a BSDE from it. 

Now we can come back to our linear BSDE
\begin{equation}
  \left\{
  \begin{split}
  &dY_t=\{A(t)Y_t+\sum_{j=1}^mB_j(t)Z^j_t+f(t)\}dt+Z_tdW_t, \hspace{1mm} t\in [0,T] \\
  &Y_T=\xi
  \end{split}
\right.
  \end{equation}
  
and show how (C.2.13) is dual to an SDE similar to (C.2.5) in the Hilbert space $L^2_\mathcal{F}(\Omega;C([0,T]);\mathbb{R}^k)\times L^2_\mathcal{F}([0,T];\mathbb{R}^{k\times m})$ using Riesz Representation Theorem

\begin{equation}
 \begin{split}
  &I(\phi,\psi):=E\bigg[\int_0^T \langle X_t,-f(t) \rangle dt+ \langle X_T,\xi\rangle \bigg]\\
 &\forall (\phi,\psi)\in L^2_\mathcal{F}(\Omega;C([0,T]);\mathbb{R}^k)\times L^2_\mathcal{F}([0,T];\mathbb{R}^{k\times m})
 \end{split}
\end{equation}

where $X_{(\cdot)}$ is the solution of the SDE:
\begin{equation}
 \left\{
 \begin{split}
  &dX_t=(-A(t)^TX_t+\phi_t)dt+\sum_{j=1}^m (-B_j(t)^TX_t+\psi_t)dW_t\\
  &X_0=0
 \end{split}
\right.
\end{equation}

$I(\cdot,\cdot)$ is linear bounded, hence by Riesz there is a unique $(\bar{Y},\bar{Z})\in L^2_\mathcal{F}(\Omega;C([0,T]);\mathbb{R}^k)\times L^2_\mathcal{F}([0,T];\mathbb{R}^{k\times m})$ such that

\begin{equation}
 I(\phi,\psi)=\mathbb{E}\int_0^T\{\langle\phi(t),\bar{Y}_t\rangle+\langle\psi(t),\bar{Z}_t\rangle \}dt
\end{equation}

On the other hand, if $(Y_{(\cdot)},Z_{(\cdot)})$ is the solution of the BSDE (C.2.13), then by applying Ito's formula to $X_t\cdot Y_t$ and assuming the local martingale part of it is a true martingale we can derive that  (C.2.16) holds with $(\bar{Y}_{(\cdot)},\bar{Z}_{(\cdot)})$ replaced by $(Y_{(\cdot)},Z_{(\cdot)})$. Due to the uniqueness of the Riesz representation $(\bar{Y}_{(\cdot)},\bar{Z}_{(\cdot)})=(Y_{(\cdot)},Z_{(\cdot)})$.

$\vspace{1mm}$

$\mathbf{Remark} $\\
The SDE for $X_t$ appears in the proof of the SMP with control over volatility as the variational equation (in our case (C.2.5)) and the corresponding first order adjoint process reads as the following BSDE. 
\begin{equation}
 \begin{split}
 &-dp_t= H_x(X_s,u_t,p_t,q_t)ds+q_tdW_t\\
 &p_T=g_x(X_T)
 \end{split}
\end{equation}

\section{Systems of coupled Forward and Backward SDEs}

In the stochastic optimal control problems, the stochastic Hamiltonian system as introduced in the relative chapter is a system of coupled forward and backward SDEs (FBSDE) where the forward component $X_{(\cdot)}$ (state process) and the backward components $(p_{(\cdot)},q_{(\cdot)})$ (first-order adjoint process, shadow price) and $(P_{(\cdot)},Q_{(\cdot)})$ (second order adjoint process, risk sensitivity) are coupled through the maximum condition. 

So here we are interested in studying those systems and more specifically the case where the SDE is n-dim and the BSDE is 1-dim where we can prove existence and uniqueness for a fairly general system. In general coupled FBSDE are not necessarily solvable and there are very few classes of certainly solvable systems. 

First, following  \cite{YZ1999} or \cite{MPJ1994} or \cite{Delarue2002} we will introduce a heuristic derivation for the (n,1)-dim system and then state the existence and uniqueness theorem. We consider:

\begin{equation}
 \left\{
 \begin{split}
  &dX_t=b(t,X_t,Y_t,Z_t)dt+\sigma(t,X_t,Y_t,Z_t)dW_t\\
  &dY_t=h(t,X_t,Y_t,Z_t)dt+Z_tdW_t\\
  &X_0=x\in \mathbb{R}^n\\
  &Y_T=g(X_T)
 \end{split}
\right.
\end{equation}

Suppose $(X_t,Y_t,Z_t)$ is an adapted solution of (23) and Y,X are  related by:
$$ Y_t=\theta(t,X_t), \hspace{1mm} \forall t \in [0,T], \mathbb{P}-a.s.$$
where $\theta$ is deterministic and belongs to $C^{1,2}$. Then by Ito's formula we have: 
\begin{equation}
\begin{split}
 dY_t=d\theta(t,X_t)=&\{\theta_t(t,X_t)+\theta_x(t,X_t)b\big(t,X_t,\theta(t,X_t),Z_t\big)\\
 &+\frac{1}{2}tr\bigg[\theta_{xx}(t,X_t)\sigma \sigma^T\big(t,X_t,\theta(t,X_t),Z_t\big)\bigg]\}dt\\
 &+\{\theta_x(t,X_t)\sigma(t,X_t,\theta(t,X_t),Z_t) \}dWt
 \end{split}
 \end{equation}

by comparing (24) with (23) we get:
\begin{equation}
\left\{
\begin{split}
 &h(t,X_t,\theta(t,X_t)=\theta_t(t,X_t)+\theta_x(t,X_t)b\big(t,X_t,\theta(t,X_t),Z_t\big)+\frac{1}{2}tr\bigg[\theta_{xx}(t,X_t)\sigma \sigma^T\big(t,X_t,\theta(t,X_t),Z_t\big)\bigg]\\
 &\theta(T,X_T)=g(X_T)
 \end{split}
\right.
\end{equation}

\begin{equation}
 \theta_x(t,X_t)\sigma(t,X_t,\theta(t,X_t),Z_t)=Z_t
\end{equation}

The above argument suggests that we design the following four-step scheme:
\begin{enumerate}
 \item[$\mathbf{Step}$ 1] Find $z(t,x,y,p)$ satisfying the following:
 \begin{equation}
 \begin{split}
  z(t,x,y,p)=p\sigma(t,x,y,z(t,x,y,z(t,x,y,p))\\
  \forall (t,x,y,p)\in [0,T]\times \mathbb{R}^n\times \mathbb{R}\times \mathbb{R}^{1\times n}
 \end{split}
 \end{equation}

 \item[$\mathbf{Step}$ 2] Use $z$ obtained above to solve the parabolic problem for $\theta(t,x)$:
 \begin{equation}
 \left\{
\begin{split}
  &\theta_t(t,x)+\theta_x(t,x)b\big(t,x,\theta(t,x),z(t,x,y,p)\big)\\
  &+\frac{1}{2}tr\bigg[\theta_{xx}(t,x)\sigma\sigma^T\big(t,X,\theta(t,x),z(t,x,y,p)\big)\bigg]\\
  &-h(t,X,\theta(t,x),z(t,x,y,p))=0 \hspace{2mm} (t,x)\in[0,T]\times \mathbb{R}^n  \\
  &\theta(T,x)=g(x) \hspace{2mm} x\in\mathbb{R}^n
   \end{split}
\right.
 \end{equation}
\item[$\mathbf{Step}$ 3] Solve the SDE
\begin{equation}
 \left\{
 \begin{split}
  &dX_t=b\big(t,X_t,\theta(t,X_t),z(t,X_t,Y_t,p_t)\big)dt+ \sigma\big(t,X_t,\theta(t,X_t),z(t,X_t,Y_t,p_t)\big)dW_t\\
  &X_0=x
 \end{split}
\right.
\end{equation}

\item[$\mathbf{Step}$ 4] Set

\begin{equation}
 \left\{
 \begin{split}
  &Y_t:=\theta(t,X_t)\\
  &Z_t:=z(t,X_t,\theta(t,X_t),\theta_x(t,X_t))
 \end{split}
\right.
\end{equation}

\end{enumerate}

And this way the triple $(X_t,Y_t,Z_t)$ will provide an adapted solution to (23) 

\subsection{Implementation of the scheme}

The main challenge to implement the above scheme is the solution of the boundary value problem (28). For this we are going to use the results from the thoery of quasi linear parabolic equations and systems for the general case (where the BSDE is k dim, and the SDE n). We refer to the original work of Ladynzhenskaya Solonnikov and Ural'tseva 1968 \cite{LSU1968} and Edmunds and Peletier 1971 \cite{EP1971}  for a review. Ma et al. 1994 \cite{MPJ1994} were first to discuss the 4-step scheme and use the PDE approach to solve it for local times and Delarue 2002 \cite{Delarue2002} extended their result. 

The method to use the scheme in practice, in case the PDE (28) cannot be solved explicitly (which is the most probable scenario) is:

\begin{enumerate}
 \item Prove existence and uniqueness of (28) 
 \item Solve (28) numerically 
 \item Use a numerical scheme for the SDE (29)
 \item Set $Y_t,Z_t$ according to (30)
\end{enumerate}

For the sake of illustration we will give examples in the next section for the scheme's Implementation

\subsubsection{Existence and uniqueness of Quasi Linear Parabolic PDEs}

We will now discuss briefly the existence and uniqueness result for (28) without too much involvement with the PDE theory. 

The solvability of the boundary value problem is proved on the basis of the Leray-Schauder theorem and a priori estimates of the norms in the spaces involved in the general case in the original work form Ladynzhenskaya Solonnikov and Ural'tseva (1968). We will state the theorem as a lemma and use it to provide existence 

Now we have to make some assumptions to gain our result. 

$\vspace{1mm}$

needs revision!!!!!!!

\begin{description}
 \item[Assumptions]
 $\vspace{1mm}$
 \begin{enumerate}
 \item $m=n$ for (23) and $b$, $\sigma$, $h$, $g$ are smooth with uniformly bounded first-order derivative taking values in $\mathbb{R}^n$, $\mathbb{R}^{n\times n}$, $\mathbb{R}$, $\mathbb{R}$ respectively   
 \item The map $z \to -|\sigma(t,x,y,z)^T|^{-1}z $ is uniformly monotone. 
\end{enumerate}

\end{description}

$\vspace{1mm}$

From step 1 we get 

$$z= \sigma(t,x,y,z)^T p $$

which used in step 2 yields

\section{Examples}

\subsection{Application to Option pricing and alternative proof of the Black-Scholes formula}

Here we will apply the theory that was developed in the previous sections in pricing a European option. What follows is rather classical for the mathematical finance literature and can be found in several textbooks, we will follow El Karoui et al. (1997) \cite{EKPQ1997} and the book \cite{YZ1999}. We will mainly focus on the BSDEs and the mathematics rather than the finance theory with market's completeness etc for the rigorous formal approach we refer to \cite{EKPQ1997}

We will study a complete market we two assets one riskless $B_t$ called bond and and one risky asset $S_t$ called stock. Also we will assume an investor who has a total wealth $Y_t$ and invests $\pi_t$ in the risky asset. The dynamics are described by:

\begin{equation}
 \left\{
 \begin{split}
  &dB_t=r_tB_tdt \hspace{1mm} t\in[0,T]\\
  &B_0=b_0\in\mathbb{R}
 \end{split}
\right.
\end{equation}

\begin{equation}
 \left\{
 \begin{split}
  &dS_t=\mu(t)S_tdt+\sigma(t)S_tdW_t  \hspace{1mm} t\in[0,T]\\
  &S_0=x_0\in\mathbb{R}^+
 \end{split}
\right.
\end{equation}

\begin{equation}
 \left\{
 \begin{split}
&dY_t=N_S(t)dS_t+N_B(t)dB_t\\
&Y_0=y_0\in\mathbb{R}
 \end{split}
\right.
\end{equation}

We need to make some remarks here:
\begin{itemize}
 \item We assume the same probability space as it was introduced in the introduction
 \item (31) is an ODE while (32) is the familiar Geometric B.M. and (33) gives us the evolution of the wealth process $N_S(t)$ the number of shares of the stock and $N_B(t)$ the number of shares of the bond
 \item $r_t,\mu(t),\sigma(t)$ are predictable bounded processes for the sake of simplicity.
 \item $\pi_t=N_S(t)S_t$ and $Y(t)-\pi_t=N_B(t)B_t$
 \item We have control over the number of shares for both of them but because we can express the wealth process as a function of $\pi_t$ and we assume no risk preference we will use as control variable $\pi_{(\cdot)}$, we can always translate our strategy in terms of $N_S,N_B$. 
\end{itemize}

We further manipulate (33) and get:

$$
dY_t=\frac{\pi(t)}{S_t}dS_t+r_t(Y_t-\pi_t)dt
$$

\begin{equation}
dY_t=\{r_tY_t+[\mu(t)-r_t]\pi_t\}dt+\sigma(t)\pi_tdW_t
\end{equation}

Suppose now that the investor wants to sell a European option, the payoff of this option at maturity T is $\xi\in L^2$. The aim of the investor is to define the minimum initial amount of capital $y_0$ such that he can cover the payoff $\xi$ at time T.

So this is a BSDE problem and we can use the 4-stem scheme from section 3 to solve it. 

The FBSDE system reads as follows for $Z_t=\pi_t\sigma(t)$

\begin{equation}
 \left\{
 \begin{split}
  &dS_t=\mu(t)S_tdt+\sigma(t)S_tdW_t  \hspace{1mm} t\in[0,T]\\
  &dY_t=\{r_tY_t+[\mu(t)-r_t]\frac{Z_t}{\sigma(t)}\}dt+Z_tdW_t\\
  &Y_T=\xi\\
  &S_0=x_0
 \end{split}
\right.
\end{equation}

In this particular case the FBSDE is decoupled since $dS(t)$ involves no $Y(t)$ and $dY(t)$ involves no $S(t)$ 

\begin{itemize}
 \item [Step 1]
 Set 
 $$z(t,s,y,p)=\sigma(t)xp, \hspace{1mm} (t,s,y,p)\in [0,T]\times \mathbb{R}^3 $$
 \item [Step 2]
 Solve the PDE
 
\begin{equation}
 \left\{
 \begin{split}
&\theta_t + \frac{\sigma(t)^2s^2}{2}\theta_{ss}+r(t)s\theta_s-r(t)\theta=0 \hspace{2mm} (t,s)\in [0,T]\times\mathbb{R}\\
&\theta|_{t=T}=\xi
 \end{split}
\right.
\end{equation}

\item[Step 3]
Solve the SDE 

$$
 \left\{
 \begin{split}
  &dS_t=\mu(t)S_tdt+\sigma(t)S_tdW_t  \hspace{1mm} t\in[0,T]\\
  &S_0=x_0
 \end{split}
\right.
$$

\item[Step 4]
Set

\begin{equation}
 \left\{
 \begin{split}
  &Y_t=\theta(t,S_t)\\
  &Z_t=\sigma(t)S_t\theta_s(t,S_t)
 \end{split}
\right.
\end{equation}

\end{itemize}

Then the option price, at t=0 will be given by 

$$Y_0=y_0=\theta(0,s)$$

\subsubsection{An alternative proof of the Black Scholes formula}

To illustrate more on (35),(36) suppose we have a put option so $Y_T=(K-S_T)^+$ and $r_t=r,\mu(t)=\mu,\sigma(t)=\sigma$ are positive constants. Then (36) is the classical Black-Scholes PDE 

$$
 \left\{
 \begin{split}
&\theta_t + \frac{\sigma^2s^2}{2}\theta_{ss}+rs\theta_s-r\theta=0 \hspace{2mm} (t,s)\in [0,T]\times\mathbb{R}  \\
&\theta|_{t=T}=(K-s)^+
 \end{split}
\right.
$$

and at $s=0$ we have 

$$
\left\{
 \begin{split}
&\theta_t-r\theta=0\\
&\theta|_{t=T}=K
 \end{split}
\right.
$$

and so $\theta(t,0)=Ke^{r(t-T)}$. Therefore $\theta(t,s)$ for $s>0$(as stock prices can never be zero) solves:

\begin{equation}
 \left\{
 \begin{split}
&\theta_t + \frac{\sigma^2s^2}{2}\theta_{ss}+rs\theta_s-r\theta=0 \hspace{2mm} (t,s)\in [0,T]\times(0,\infty) \\
&\theta|_{s=0}=Ke^{r(t-T)} \hspace{2mm} t\in [0,T]\\
&\theta|_{t=T}=(K-s)^+ \hspace{2mm} s\in(0,\infty)
 \end{split}
\right.
\end{equation}

To solve (38) we can consider the successive changes of variables:
\begin{itemize}
 \item First the state $x=lns$ and $\phi(t,x)=\theta(t,e^s)$ satisfies
 
\begin{equation}
 \left\{
 \begin{split}
  &\phi_t+\frac{\sigma^2}{2}\phi_{xx}+(r-\frac{\sigma^2}{2})\phi_x-r\phi=0 \hspace{2mm} (t,x)\in [0,T]\times \mathbb{R}\\
  &\phi|_{t=T}=(K-e^x)^+ \hspace{2mm}  x\in \mathbb{R}
 \end{split}
\right.
\end{equation}

\item Then time $\tau=\gamma t $ and $\psi(\tau,x)=e^{-\frac{\alpha \tau}{\gamma}-\beta x}\phi(\frac{\tau}{\gamma},x)$ with 

$$\alpha=r+\frac{1}{2\sigma^2}(r-\frac{\sigma^2}{2})^2 $$
$$\beta=-\frac{1}{\sigma^2}(r-\frac{\sigma^2}{2}) $$
$$\gamma=\frac{\sigma^2}{2} $$

then $\psi(\tau,x)$ satisfies

\begin{equation}
 \left\{
 \begin{split}
  &\psi_\tau+\psi_{xx}=0 \hspace{2mm} (\tau,x)\in [0,\gamma T]\times \mathbb{R}\\
  &\psi|_\tau=e^{-\frac{\alpha T}{\gamma}-\beta x}(K-e^x)^+ \hspace{2mm}  x\in \mathbb{R}
 \end{split}
\right.
\end{equation}

\end{itemize}

Now we have transform (38) into (40), a simple heat equation which can be solved explicitly by common techniques (separation of variables etc) which in the end yields the familiar formula:

\begin{equation}
\left\{
 \begin{split}
 &\theta(t,s)=Ke^{-r(t-T)}N(-d_2)-N(-d_1)S_t\\
 &d_1=\frac{1}{\sigma\sqrt{T-t}}[\ln(\frac{S_t}{K})+(r+\frac{\sigma^2}{2})(T-t)]\\
 &d_2=d_1-\sigma\sqrt{T-t}\\
 &N(x)=\frac{1}{\sqrt{2\pi}}\int_{-\infty}^x e^{-\frac{z^2}{2}}dz
 \end{split}
\right.
\end{equation}

\subsection{A linear case of FBSDE}

Here we will study a linear one dimensional FBSDE to elaborate more on the 4 step scheme. The particular example is only pedagogical with no interpretation in finance or physics. We consider

\begin{equation}
 \left\{
 \begin{split}
  &dX(t)=\{X(t)+Y(t)\}dt+\{X(t)+Y(t)\}dW(t)\\
  &dY(t)=\{X(t)+Y(t)\}dt+Z(t)dW(t)\\
  &X(0)=x_0\in \mathbb{R}\\
  &Y(T)=g(X(T)
 \end{split}
\right.
\end{equation}

We will think about the terminal condition later to ensure the wellposedness of the problem. We apply the 4 step scheme. 

\begin{enumerate}
 \item [Step 1] 
 $$z(t,x,y,p)=p(x+y) \hspace{2mm} (t,x,y,p)\in [0,T]\times \mathbb{R}\times \mathbb{R}\times \mathbb{R}$$
 \item[Step 2]
 We will solve 
 $$
 \left\{
 \begin{split}
 &\theta_t+\frac{1}{2} (x+\theta)^2\theta_{xx}+(x+\theta)\theta_x+(x+\theta) \hspace{2mm} (t,x)\in [0,T]\times \mathbb{R}\\
 &\theta(T,x)=g(x) \hspace{2mm}\in\mathbb{R}
 \end{split}
\right.
 $$
\end{enumerate}

 \addcontentsline{toc}{chapter}{Bibliography}
 
\end{document}